\documentclass[11pt,reqno]{amsart}

\usepackage{amsmath,amssymb,amsthm,mathrsfs}
\usepackage[margin=0.9in]{geometry}

\usepackage{enumitem}
\usepackage{booktabs}
\usepackage{array}
\usepackage[dvipsnames]{xcolor}
\usepackage{tikz}
\usepackage{pgfplots}
\pgfplotsset{compat=1.18}
\usepgfplotslibrary{fillbetween}
\usetikzlibrary{arrows.meta,decorations.pathreplacing,patterns}

\usepackage{float}

\definecolor{TitleRed}{rgb}{0.60,0.00,0.00}
\definecolor{LinkBlue}{rgb}{0.00,0.25,0.70}
\definecolor{URLRed}{rgb}{0.55,0.00,0.10}

\definecolor{RevColor}{RGB}{0,120,60}
\definecolor{CheckColor}{RGB}{190,0,120}

\definecolor{lava}{rgb}{0.81, 0.06, 0.13}
\usepackage[colorlinks=true,linkcolor=lava,citecolor=blue,urlcolor=blue!60!black]{hyperref}

\theoremstyle{plain}
\newtheorem{theorem}{Theorem}[section]
\newtheorem{proposition}{Proposition}[section]
\newtheorem{lemma}{Lemma}[section]
\newtheorem{corollary}{Corollary}[section]
\newtheorem{conjecture}{Conjecture}[section]

\theoremstyle{definition}
\newtheorem{definition}{Definition}[section]
\newtheorem{assumption}{Assumption}[section]

\newtheorem{openproblem}{Open Problem}[section]

\theoremstyle{remark}
\newtheorem{remark}{Remark}[section]
\newcommand{\R}{\mathbb{R}}

\newcommand{\RN}{\R^N}
\newcommand{\Lp}[1]{L^{#1}(\RN)}
\newcommand{\Lpp}[1]{L^{#1}_{+}(\RN)}
\newcommand{\eps}{\varepsilon}
\newcommand{\ol}{\overline}
\newcommand{\pc}{p_{c}}
\newcommand{\pF}{p_{F}}
\newcommand{\pG}{p_{G}}
\newcommand{\pss}{p^{\ast\ast}}
\newcommand{\ac}{\alpha_{c}}
\newcommand{\Tmax}{T_{\max}}
\newcommand{\norm}[1]{\left\lVert #1\right\rVert}
\newcommand{\abs}[1]{\left\lvert #1\right\rvert}
\newcommand{\dx}{\,dx}
\newcommand{\ds}{\,ds}
\newcommand{\dy}{\,dy}
\newcommand{\eqdef}{:=}
\newcommand{\Dprime}[1]{\mathcal{D}'(#1)}
\newcommand{\Hmeas}{\mathcal{H}^{N-1}}

\DeclareMathOperator{\supp}{supp}
\DeclareMathOperator{\sgn}{sign}

\numberwithin{equation}{section}

\definecolor{orcidlime}{HTML}{A6CE39}

\newcommand{\orcidauthorB}{0000-0001-6038-1069}

\title[Critical decay and amplitude threshold]
{Critical decay and amplitude threshold for a porous medium equation with a singular weighted source}

\author[ R. Ben Belgacem and M. Majdoub]{Rihab Ben Belgacem and Mohamed Majdoub}
\address[R. Ben Belgacem]{Universit\'e de Tunis El Manar, Facult\'e des Sciences de Tunis, D\'epartement de math\'ematiques, Laboratoire \'equations aux d\'eriv\'ees partielles (LR03ES04), 2092 Tunis, Tunisie.}
\email{\tt rihabbenbelgacem662@gmail.com  }
\address[M. Majdoub]{Department of Mathematics, College of Science, Imam Abdulrahman Bin Faisal University, P. O. Box 1982, Dammam, Saudi Arabia.}
\address[M. Majdoub]{Basic and Applied Scientific Research Center, Imam Abdulrahman Bin Faisal University, P.O. Box 1982, 31441, Dammam, Saudi Arabia.}
\email{\tt mmajdoub@iau.edu.sa}
\email{\tt med.majdoub@gmail.com}
\email{\tt mohamed.majdoub@fst.rnu.tn}

\date{\today}

\begin{document}
\begin{abstract}
We study the Cauchy problem for the parabolic Hardy--H\'enon equation
$\partial_tu=\Delta u^m+|x|^\sigma u^p$ in $\mathbb{R}^N$, $N\ge3$, with slow
diffusion $m>1$ and a weight $|x|^\sigma$, $-2<\sigma<0$, which is singular at
the origin. Above the Fujita exponent $p_F=m+(\sigma+2)/N$, global and
blowing-up solutions coexist, and the existing theory does not determine how
the spatial decay of the initial data selects between them. We prove that the
second critical exponent is $\alpha_c=(\sigma+2)/(p-m)$, the decay rate of the
explicit singular steady state $S(x)=A|x|^{-\alpha_c}$. Data decaying more
slowly than $|x|^{-\alpha_c}$ blow up in finite time whatever their size,
while faster decaying data which are small in a weighted supremum norm give
rise to global solutions. At the critical decay the decay rate alone does
not decide, and we study the amplitude threshold along the family
$L\min\{1,|x|^{-\alpha_c}\}$.
The singular weight prevents a direct smoothing of $S$ into a bounded barrier,
and the data under consideration are not integrable. We therefore establish a
comparison principle for barriers which are unbounded at the origin and not
integrable at infinity, and we construct bounded stationary supersolutions
with exactly critical tail by gluing multiples of $S$ to regular steady
states. In the range where $S$ exists, it follows that the amplitude threshold
lies between $A$ and an explicit constant, that data lying below a proper
fraction of $S$ produce global bounded solutions, and that the amplitude $A$
is attained by global solutions whenever a Joseph--Lundgren type condition
fails. Where $S$ does not exist, forward self-similar solutions provide the
barriers, and the amplitude threshold is positive in the whole supercritical
range. We formulate conjectures on the exact value of the amplitude threshold.
\end{abstract}
\subjclass[2020]{{Primary 35K65, 35B33, 35B44; Secondary 35K67, 35K57,
35K59, 35K15, 35B51, 35A01, 35C06, 35J75}}

\keywords{{Parabolic Hardy--H\'enon equation, porous medium equation,
singular weight, second critical exponent, critical decay, amplitude threshold,
singular steady state, comparison principle, forward self-similar solutions,
Joseph--Lundgren exponent}}

\maketitle
%=====================================================================
\section{{Introduction and main results}}\label{sec:intro}
%=====================================================================

In this paper we study the Cauchy problem for the parabolic
Hardy--H\'enon equation with quasilinear degenerate diffusion,
\begin{equation}\label{eq:CP}
\left\{
\begin{aligned}
\partial_t u&=\Delta u^m+|x|^{\sigma}u^{p},
&&t>0,\ x\in\mathbb{R}^N,\\
u(0)&=u_0,&&x\in\mathbb{R}^N,
\end{aligned}
\right.
\end{equation}
where $N\ge 3$, $m>1$, $p>1$ and $\sigma\in(-2,0)$. The low-dimensional
cases $N\in\{1,2\}$ are discussed in Section~\ref{sec:discussion}.
Equation \eqref{eq:CP} couples a porous medium diffusion $\Delta u^m$ with
a superlinear source $|x|^{\sigma}u^{p}$, whose weight is singular at the
origin because $\sigma<0$. Equations of this type are called parabolic
Hardy--H\'enon equations. The name refers to H\'enon's model for rotating
stellar systems \cite{Henon73}, whose elliptic counterpart was studied by Ni
\cite{Ni82}, and to the work of Baras and Goldstein on the heat equation with
an inverse square potential. In the latter, the optimal constant of the Hardy
inequality separates existence from instantaneous blow-up \cite{BG84}.

In the semilinear case $m=1$, the Cauchy problem \eqref{eq:CP} has been
studied in depth. We refer to \cite{BTW17,CIT22,CITT24} for well-posedness,
critical weighted Lebesgue spaces, forward self-similar solutions and
uniqueness, to \cite{HisaTakahashi21,HisaSierzega24} for initial data with
singularities, to \cite{TayachiWeissler23} for life-span estimates, and to
\cite{Pinsky97} for the influence of a space-dependent weight on global
solvability. Related well-posedness, blow-up and Fujita-type results for
inhomogeneous and Hardy--H\'enon parabolic equations can be found in
\cite{MajdoubMliki21,Majdoub21,Majdoub23,MajdoubTorebek26a,MajdoubTorebek26b}.
A systematic account of the unweighted theory is given in \cite[Ch.~IV]{QS19}.
Reaction--diffusion equations of the form \eqref{eq:CP} also arise in models
of diffusing concentrations, population dynamics and heat conduction, see
\cite{Hu11,GalaktionovVazquez02,SGKM95} and the references therein. In all
these works, a central question is to determine the critical exponents which
separate global existence from finite-time blow-up.

This question has a long history. It starts with the work of Fujita
\cite{Fujita66} on the unweighted semilinear problem ($m=1$, $\sigma=0$),
where an exponent $p_F$ was shown to separate the subcritical, critical and
supercritical regimes. Since then, Fujita-type exponents have been determined
for a wide class of parabolic problems. The available techniques rely on
spectral arguments, on scaling, or on the construction of self-similar sub-
and supersolutions, see for instance
\cite{Levine90,DengLevine00,Qi98,MochizukiSuzuki97,Suzuki02,BLZ00,MajdoubMliki26}
and the references therein. Fujita-type thresholds in the presence of forcing terms
are studied in \cite{Majdoub23,AAAAM23,BBM26,MajdoubTorebek26a,KumarMajdoub26}.

For the quasilinear equation \eqref{eq:CP} with $\sigma<0$, two exponents
play a distinguished role, namely
\begin{equation}\label{eq:pGpF}
p_G\eqdef 1-\frac{\sigma(m-1)}{2}\in(1,m),\qquad
p_F\eqdef m+\frac{\sigma+2}{N}>m .
\end{equation}
The Fujita exponent $p_F$ was identified by Qi \cite{Qi98}, who proved that
every nontrivial nonnegative solution blows up when $m<p\le p_F$, whereas some
solutions are global when $p>p_F$. The exponent $p_G$ first appeared in the
classification of radially symmetric self-similar solutions for $1<p<m$, whose
form depends on the sign of $p-p_G$ \cite{IMS23,ILS24a,ILS24b}. Iagar and
Lauren\c{c}ot \cite{IL1} then showed that $p_G$ also governs the well-posedness
theory, and they carried out a detailed qualitative analysis in the whole range $m>1$, $p>1$, $\sigma\in\big(\max\{-2,-N\},0\big)$. We take
this work as our point of departure. In \cite{IL1} the authors establish the
local existence of weak solutions, together with uniqueness and a comparison
principle for data in $L^1\cap L^\infty$ (unconditionally when $p\ge p_G$ and
$N\ge3$). They prove global existence for $p\in(1,p_G]$, and finite-time
blow-up of every nontrivial nonnegative solution for $p_G<p\le p_F$. They also
obtain blow-up under a negative energy condition for $p\ge m$, global
existence for small data in the critical Lebesgue space $L^{r_0+1}$ when
$p>p_F$, and optimal grow-up rates for $p\in(1,p_G]$. Here the critical
Lebesgue exponent is governed by
\begin{equation}\label{eq:r0}
r_0\eqdef \frac{N}{\sigma+2}(p-m)-1=\frac{N}{\sigma+2}(p-p_F),
\qquad r_c\eqdef\max\{r_0,0\}.
\end{equation}

The H\'enon regime $\sigma>0$ was considered by the same authors in
\cite{IL2}, for $m>0$ and $p>\max\{1,m\}$. In \cite[Theorem~1.1]{IL2} they
determine the decay rates of bounded data $u_0\ge0$ with
$\lim_{|x|\to\infty}|x|^{\gamma}u_0(x)=L\in(0,\infty)$ for which no
nonnegative solution exists. Here $p_G=1+\sigma(1-m)/2$, which is larger than
$1$ only if $m<1$, and their result reads as follows:
\begin{itemize}
\item if $p\ge p_G$, non-existence holds for $\gamma<\sigma/(p-1)$;
\item if $1<p<p_G$ (a range which, for $\sigma>0$, requires $m<1$),
non-existence holds for $\gamma<(\sigma+2)/(p-m)$;
\item if $1<p<p_G$ and $\gamma=(\sigma+2)/(p-m)$, non-existence holds for
all amplitudes $L$ larger than some constant \cite[Eq.~(1.7)]{IL2}.
\end{itemize}
In the concluding discussion of \cite[\S3(c)]{IL2}, the authors also exhibit,
for $N\ge3$ and $p>p_c\eqdef m(\sigma+N)/(N-2)$, the explicit singular
stationary solution $S(x)=A|x|^{-(\sigma+2)/(p-m)}$. The same exponent
$p_c$ appears, for every $\sigma>\max\{-2,-N\}$, in the phase-space analysis
of global self-similar solutions carried out in \cite{ISS26}, where it
governs the existence of a critical point of the associated dynamical system
\cite[Section~2.1]{ISS26}. In the remark following \cite[Lemma~2.1]{ISS26},
the authors write $S$ explicitly, with the same amplitude, as the solution
corresponding to an orbit of that system, and they note that it exists for
$p>p_c$. {The singular stationary solution also plays a central role in
the recent work of Iagar, S\'anchez and Sarri\'on-Pedralva \cite{ISSpeak} on
peaking solutions, for $\sigma\in(-2,\infty)$ and $p>p_S(\sigma)$, an
exponent which coincides with $p^{\ast\ast}$ in \eqref{eq:hierarchy} below;
the constant $C(\sigma)$ of \cite[Eq.~(2.5)]{ISSpeak} is the amplitude $A$
of \eqref{eq:A}. There, $S$ is used as an upper barrier for the minimal
solution with the unbounded initial value $K|x|^{-(\sigma+2)/(p-m)}$,
$0<K<A$, and this solution is shown to be global, bounded for $t>0$ and of
forward self-similar form \cite[Proposition~3.1]{ISSpeak}. Backward
self-similar solutions which blow up at a finite time, and whose profiles
have exactly critical tail with amplitude in $(0,A)$, are also constructed
there. Their existence ranges involve a Joseph--Lundgren exponent
$p_{JL}(\sigma)$ \cite[Eq.~(2.4) and Theorems~2.2--2.4]{ISSpeak}, which we
relate to the present analysis in Section~\ref{sec:conjectures}.} The singular steady state is the object which
organises the whole of the present paper. The non-existence threshold
$\sigma/(p-1)$ of \cite{IL2} plays no role here, since it is negative, hence
vacuous, for $\sigma<0$.

Let us recall the classical picture that motivates the present work. For the
semilinear Fujita problem \(\partial_t u=\Delta u+u^{p}\) with \(p>1+2/N\),
global and blowing-up solutions coexist. Lee and Ni \cite{LeeNi92} showed that
the spatial decay of the initial datum decides between them. Data satisfying
\(\liminf_{|x|\to\infty}|x|^{a}u_0(x)>0\) for some \(a<2/(p-1)\) blow up in
finite time, whereas data bounded by a small multiple of
\((1+|x|^2)^{-1/(p-1)}\) give rise to global solutions; see also
\cite{GuiWang95} for the corresponding life-span estimates. The exponent
\(2/(p-1)\) is called the \emph{second critical exponent}. It is also the
decay rate of the singular steady state \(u_s\), which exists for
\(p>N/(N-2)\), and it is dual to the critical Lebesgue exponent \(N(p-1)/2\)
of the small-data theory. Wang \cite{Wang93} proved, moreover, that every
datum lying below \(\lambda u_s\), with \(\lambda\in(0,1)\), generates a
global classical solution. Second critical exponents have been obtained for
porous medium and fast diffusion with an unweighted source
\cite{MMH00,GuoGuo01}, for weights which are bounded near the
origin \cite{AndreucciDiBenedetto91,Suzuki02}, for other degenerate or doubly
nonlinear diffusions \cite{LiMu04,MuLiWang10,MMZ11,YangYangZheng12}, and for
weighted semilinear equations \cite{Pinsky97,BTW17}. Near the singular steady
state, finer phenomena appear. The ordering of the regular steady states and
the stability of the singular one change at the Joseph--Lundgren exponent
\cite{JosephLundgren73,GuiNiWang92,Wang93,PolacikYanagida03}. Moreover,
regular self-similar solutions may start from the singular steady state
\cite{GalaktionovVazquez97}, and even from multiples of it lying above it
\cite{SouWei03}. In the semilinear case, forward self-similar solutions with
algebraic tails also act as thresholds between global existence and blow-up
\cite{Naito12,Naito20}.

Despite this progress, an important part of the picture was missing for
\eqref{eq:CP}. In the Fujita range $1<p\le p_F$, the dichotomy between global
existence and blow-up is completely settled by \cite{IL1}. In the
supercritical range $p>p_F$, however, only two results were available.
Negative-energy data blow up \cite[Theorem~2.9]{IL1}, and small data in
$L^{r_0+1}$ are global \cite[Theorem~2.10]{IL1}. For the reference equation
$\partial_t u=\Delta u^m+K(x)u^{p}$, with a weight $K$ decaying like
$|x|^{\sigma}$ at infinity but regularized at the origin, the decay of
$u_0(x)$ as $|x|\to\infty$ which separates finite-time blow-up from global
existence is available \cite{AndreucciDiBenedetto91,Suzuki02}; see also
\cite{ZMHZ16,YangYangZheng12} for related degenerate and weighted models. As
observed in \cite{IL1}, the singularity of $|x|^{\sigma}$ at $x=0$ produces
results which differ from those for its regularized versions. {For
\eqref{eq:CP} itself, however, the influence of the spatial decay of the data
on the alternative between global existence and blow-up had not been studied
in the range $p>p_F$.} In the terminology of Lee and
Ni \cite{LeeNi92} and of its quasilinear analogues
\cite{Suzuki02,MMH00,GuoGuo01}, the missing object is
the second critical exponent of the singular quasilinear Hardy--H\'enon
equation.

The present paper fills this gap. {Its contributions are of three kinds.
For every $p>p_F$, we identify the second critical exponent
$\alpha_c=(\sigma+2)/(p-m)$: data with a slower algebraic tail blow up
whatever their amplitude, while data with a faster tail which are small in a
weighted supremum norm are global (Theorems~\ref{thm:slow} and
\ref{thm:sce}). At the critical decay, we show that the amplitude threshold
$L_\ast(\alpha_c)$ is positive and finite for every $p>p_F$. When the
singular steady state exists, that is, for $p>p_c$, it satisfies
$A\le L_\ast(\alpha_c)\le L_0$ with explicit constants, and the value $A$ is
attained by global solutions when a Joseph--Lundgren type condition fails
(Theorems~\ref{thm:threshold}, \ref{thm:allp} and \ref{thm:conj91},
Proposition~\ref{rem:LequalA}). On the technical side, we prove comparison
principles with stationary and self-similar barriers which are singular at the
origin or not integrable at infinity, and a positivity property of the
self-similar profiles (Lemmas~\ref{lem:comparison}, \ref{lem:comparisonU} and
\ref{lem:positivity}). For $m=1$ and $\sigma=0$, the second critical
exponent and the amplitude dichotomy at the critical decay go back to
\cite{LeeNi92,Wang93}, and second critical exponents for weights which are
bounded near the origin are obtained in
\cite{AndreucciDiBenedetto91,Suzuki02}. For \eqref{eq:CP} with $p>p^{\ast\ast}$,
global solutions with the unbounded initial values $K|x|^{-\alpha_c}$,
$0<K<A$, are constructed in \cite{ISSpeak}. To our knowledge, the results
listed above, which concern bounded data and the whole range $p>p_F$, were
not available for \eqref{eq:CP}.}

The passage from the regularized case to the
genuinely singular one is not merely technical. When the weight $K$ is bounded
near the origin, the standard comparison arguments with bounded stationary
supersolutions apply without obstruction. In \eqref{eq:CP}, on the
contrary, the singularity of $|x|^{\sigma}$ at $x=0$ rules out the naive
construction of barriers (Remark~\ref{rem:origin}). Moreover, the uniqueness
and comparison theory of \cite{IL1} is established for integrable data,
whereas we have to deal with data having slowly decaying or critical tails,
that is,
\[
u_0\in L^\infty(\mathbb{R}^N),\qquad u_0\geq 0,\qquad
\lim_{|x|\to\infty}|x|^{\gamma}u_0(x)=L\in(0,\infty).
\]
Such data are not integrable when $\gamma\le N$. For them we expect a critical
decay rate $\gamma=\alpha_c$ separating global existence for
$\gamma>\alpha_c$ (fast decay) from finite-time blow-up for $\gamma<\alpha_c$
(slow decay), and an amplitude threshold on $L$ at the borderline
$\gamma=\alpha_c$. The main difficulty is therefore twofold. First, one has to
identify the correct critical decay rate $\alpha_c$. It is dictated by the
singular steady state $S(x)=A|x|^{-\alpha_c}$, and it is the pointwise
counterpart of the critical Lebesgue exponent of \cite[Theorem~2.10]{IL1}.
Second, one has to construct bounded barriers at criticality, and to
establish comparison principles which remain valid for functions that are
unbounded at the origin and not integrable at infinity.

Our approach is organised around the singular steady state \(S\). Its decay
rate produces \(\alpha_c\), and, being a stationary supersolution, it traps the
flow from above. After a suitable modification near the origin, it provides
bounded barriers with exactly critical tail. Two technical ingredients make
this programme work in the singular setting. The first one is a comparison
lemma for barriers which are unbounded at the origin and not integrable at
infinity, which we obtain by comparing with the regularized problems of
\cite{IL1}. The second one is a scaling computation of the energy of the datum
on far annuli. Combined with the negative-energy criterion of \cite{IL1}, it
yields blow-up. Table~\ref{tab:literature} places our results in the context
of the existing literature.

\begin{table}[ht]
\centering
\small
\renewcommand{\arraystretch}{1.25}
\begin{tabular}{>{\raggedright\arraybackslash}p{3.3cm}>{\raggedright\arraybackslash}p{3.6cm}>{\raggedright\arraybackslash}p{7.0cm}}
\toprule
Reference & Setting & Results on the decay of $u_0$ and on critical exponents\\
\midrule
Lee--Ni \cite{LeeNi92},\newline Gui--Wang \cite{GuiWang95},\newline Wang \cite{Wang93}
& $\partial_tu=\Delta u+u^p$
& Second critical exponent $2/(p-1)$ for $p>1+2/N$; life span for slowly decaying data; global existence below $\lambda u_s$, $\lambda<1$.\\
Mukai--Mochizuki--Huang \cite{MMH00}, Guo--Guo \cite{GuoGuo01}
& $\partial_tu=\Delta u^m+u^p$
& Second critical exponent $2/(p-m)$, life span.\\
Andreucci--DiBenedetto \cite{AndreucciDiBenedetto91}, Suzuki \cite{Suzuki02}
& $\partial_tu=\Delta u^m+K(x)u^p$, $K\sim|x|^{\sigma}$ at infinity, $K$ bounded near $0$
& Fujita and second critical exponents for weights regularized at the origin.\\
\cite{LiMu04,MuLiWang10,MMZ11,YangYangZheng12}
& other degenerate, $p$-Laplacian or doubly nonlinear diffusions
& Second critical exponents and life span for slowly decaying data.\\
Qi \cite{Qi98}
& \eqref{eq:CP}, $m<p$
& Fujita exponent $p_F$: blow-up of all nontrivial solutions for $m<p\le p_F$.\\
Iagar--Lauren\c{c}ot \cite{IL1}
& \eqref{eq:CP}, $m>1$, $\sigma<0$
& Well-posedness, comparison for $L^1\cap L^\infty$ data, global existence for $p\le p_G$, blow-up for $p_G<p\le p_F$, small $L^{r_0+1}$ data global and negative energy blow-up for $p>p_F$.\\
Iagar--Lauren\c{c}ot \cite{IL2}
& \eqref{eq:CP}, $\sigma>0$
& Sharp decay thresholds for \emph{non-existence}: $\sigma/(p-1)$ if $p\ge p_G$, and $(\sigma+2)/(p-m)$ if $1<p<p_G$ (only for $m<1$).\\
Iagar--S\'anchez--Sarri\'on-Pedralva \cite{ISS26}
& \eqref{eq:CP}, $m\ge1$, $p>p_F$
& Global self-similar profiles; algebraic tails $\xi^{-(\sigma+2)/(p-m)}$; explicit singular stationary solution for $p>p_c$.\\
{Iagar--S\'anchez--Sarri\'on-Pedralva \cite{ISSpeak}}
& {\eqref{eq:CP}, $m>1$, $\sigma>-2$, $p>p^{\ast\ast}$}
& {Global forward self-similar solution from $K|x|^{-\alpha_c}$, $0<K<A$; backward self-similar solutions blowing up with critical tail $K<A$ and continued after blow-up (peaking solutions); Joseph--Lundgren exponent $p_{JL}(\sigma)$.}\\
\midrule
\textbf{This paper}
& \eqref{eq:CP}, $m>1$, $\sigma\in(-2,0)$, $N\ge3$, $p>p_F$
& Second critical exponent $\alpha_c=(\sigma+2)/(p-m)$ for every $p>p_F$; for $p>p_c$: trapping below $S$, global bounded solutions below $\vartheta S$, $\vartheta<1$, bounded barriers with critical tail, and $A\le L_\ast(\alpha_c)\le L_0$; for every $p>p_F$: $L_\ast(\alpha_c)>0$ via self-similar barriers.\\
\bottomrule
\end{tabular}
\caption{Position of the present results with respect to the literature.}
\label{tab:literature}
\end{table}

\subsection{{Main results}}\label{subsec:mainresults}\label{sec:main}

We now state our main results. Throughout, \(N\ge3\), \(m>1\),
\(\sigma\in(-2,0)\) and \(p>p_F\). Solutions are understood in the sense of
the Standing Assumptions~\ref{ass:IL}, and data belong to the class
\(\mathcal A_r\) defined in \eqref{eq:class}. The functions \(\phi_\gamma\)
and the weighted norm \(\norm{\cdot}_{[\gamma]}\) are defined in
\eqref{eq:phigamma}. The notion of stationary weak supersolution is recalled in
Definition~\ref{def:supersol}. The cut-off \(\zeta\) and the constants
\(I_1\), \(I_2\), \(M_0\), which enter the blow-up threshold \(L_0\), are
introduced in \eqref{eq:zeta}--\eqref{eq:M0}. The proofs are given in
Sections~\ref{sec:proofs} and \ref{sec:proofs-blowup}.

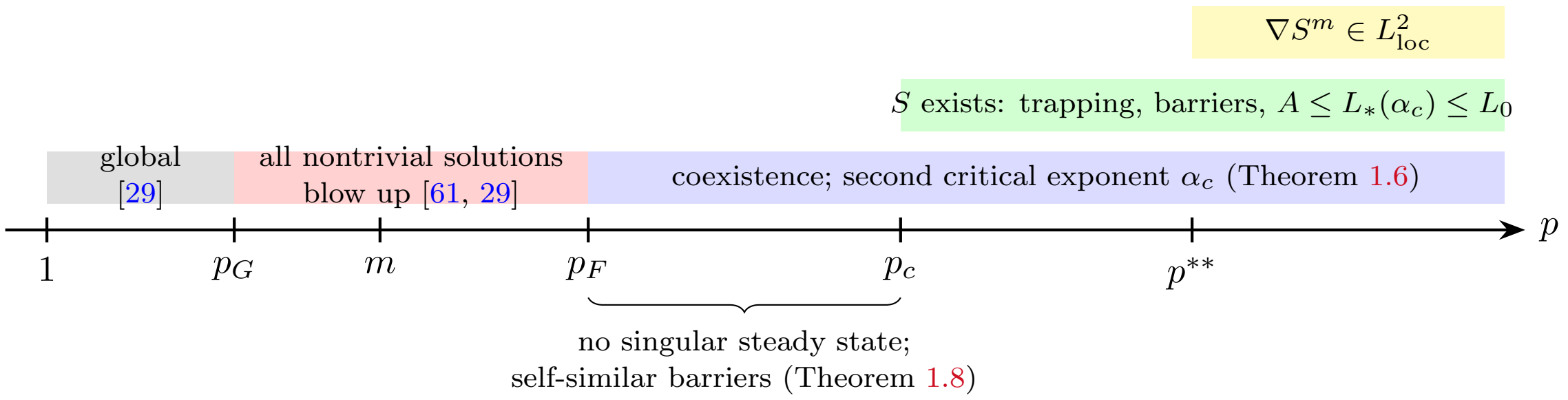
\begin{figure}[ht]
\centering
\begin{tikzpicture}[x=1cm,y=1cm,font=\small]
  \draw[thick,-{Stealth[length=2.5mm]}] (0,0) -- (14.6,0) node[right] {$p$};
  \foreach \x/\lab in {0.4/1,2.2/p_G,3.6/m,5.6/p_F,8.6/p_c,11.4/p^{\ast\ast}}
    {\draw[thick] (\x,0.12) -- (\x,-0.12) node[below] {$\lab$};}
  \fill[gray!25] (0.4,0.25) rectangle (2.2,0.75);
  \node[align=center,font=\scriptsize] at (1.3,0.5) {global\\\cite{IL1}};
  \fill[red!18] (2.2,0.25) rectangle (5.6,0.75);
  \node[align=center,font=\scriptsize] at (3.9,0.5) {all nontrivial solutions\\blow up \cite{Qi98,IL1}};
  \fill[blue!14] (5.6,0.25) rectangle (14.4,0.75);
  \node[align=center,font=\scriptsize] at (10.0,0.5) {coexistence; second critical exponent $\alpha_c$ (Theorem~\ref{thm:sce})};
  \fill[green!18] (8.6,0.95) rectangle (14.4,1.45);
  \node[align=center,font=\scriptsize] at (11.5,1.2) {$S$ exists: trapping, barriers, $A\le L_\ast(\alpha_c)\le L_0$};
  \fill[yellow!30] (11.4,1.65) rectangle (14.4,2.15);
  \node[align=center,font=\scriptsize] at (12.9,1.9) {$\nabla S^m\in L^2_{\mathrm{loc}}$};
  \draw[decorate,decoration={brace,amplitude=4pt,mirror}] (5.6,-0.65) -- (8.6,-0.65)
     node[midway,below=5pt,align=center,font=\scriptsize] {no singular steady state;\\ self-similar barriers (Theorem~\ref{thm:allp})};
\end{tikzpicture}
\caption{The scale of exponents for \eqref{eq:CP} with $m>1$,
$\sigma\in(-2,0)$ and $N\ge3$, and the ranges covered by \cite{IL1} and by the
present paper. The Joseph--Lundgren condition \eqref{eq:JL} can only hold for
$p>p^{\ast\ast}$, and it holds for all large $p$ when $N>10+4\sigma$.}
\label{fig:exponents}
\end{figure}

Our main results can be summarised as follows. For every $p>p_F$, the second
critical exponent is $\alpha_c=(\sigma+2)/(p-m)$
(Theorems~\ref{thm:slow}, \ref{thm:sce} and Corollary~\ref{cor:char}). For
$p>p_c$, every solution starting below $S$ stays below $S$
(Theorem~\ref{thm:trap}), and it is global and bounded if it starts below
$\vartheta S$ for some $\vartheta<1$ (Corollary~\ref{cor:belowS}). Moreover,
the amplitude threshold at the critical decay satisfies
$A\le L_\ast(\alpha_c)\le L_0$ (Theorems~\ref{thm:threshold} and
\ref{thm:conj91}, Corollary~\ref{cor:sandwich}). The lower bound $A$ is
obtained by gluing multiples of $S$ to regular radial steady states, and it is
optimal for barriers of this form. When a Joseph--Lundgren type condition
fails, the amplitude $A$ itself is attained by global solutions
(Proposition~\ref{rem:LequalA}). In the range $p_F<p\le p_c$, where $S$ does
not exist, forward self-similar solutions provide the barriers, and
$0<L_\ast(\alpha_c)\le L_0$ for every $p>p_F$ (Theorem~\ref{thm:allp}). To the
best of our knowledge, these results are new for the quasilinear
Hardy--H\'enon equation \eqref{eq:CP}. They are
summarised in Table~\ref{tab:summary} and Figure~\ref{fig:threshold} at the end
of this section, and the scale of exponents is displayed in
Figure~\ref{fig:exponents}.

The entire analysis revolves around the singular steady state and around a
single scaling computation on annuli. The candidate critical decay is
\begin{equation}\label{eq:ac}
\alpha_c\eqdef\dfrac{\sigma+2}{p-m}=\dfrac{N}{r_0+1}.
\end{equation}
The second identity in \eqref{eq:ac}, a consequence of \eqref{eq:r0}, is
structurally important. It expresses $\alpha_c$ as the pointwise dual of the
critical Lebesgue exponent $r_0+1$ of \cite[Theorem~2.10]{IL1}, exactly as in
the semilinear theory, where $2/(p-1)$ is dual to the Fujita-critical
Lebesgue exponent $N(p-1)/2$. Equivalently, $\alpha_c$ is the unique decay rate
which makes the profile $|x|^{-\alpha_c}$ invariant under the scaling
\begin{equation}\label{eq:scaling}
u_\lambda(t,x)\eqdef\lambda^{\alpha_c}\,u\bigl(\lambda^{\theta}t,\lambda x\bigr),
\qquad \theta\eqdef \alpha_c(m-1)+2>0,
\end{equation}
which leaves \eqref{eq:CP} invariant.

\subsubsection{The singular steady state and the exponents \texorpdfstring{$\pc$, $\pss$}{pc, pss}}\label{subsec:steady}
The critical decay $\alpha_c$ and the exponent $p_c$ come from two
requirements. An algebraic profile $A|x|^{-\alpha}$ must balance diffusion
against reaction in the stationary equation, which forces $\alpha=\alpha_c$,
and the resulting amplitude must be positive. The latter condition, namely
$\alpha_c m<N-2$, is also exactly the condition under which no Dirac mass
appears at the origin (Lemma~\ref{lem:dirac}). We therefore begin with the
singular steady state itself. We show that it exists precisely in the range
$p>p_c$, and that $p_c$ lies strictly above the Fujita exponent. The scale of
exponents reads
\begin{equation}\label{eq:hierarchy}
1<p_G<m<p_F<p_c<p^{**},\qquad
p^{**}\eqdef\dfrac{m(N+2\sigma+2)}{N-2}.
\end{equation}
The exponent $p^{**}$ is the Hardy--Sobolev critical exponent of the stationary
problem, and it is characterized by
\[
\nabla S^m\in L^2_{\mathrm{loc}}(\mathbb{R}^N)
\iff p>p^{**}
\]
(see Remarks \ref{rem:pss} and \ref{rem:pss-sobolev}). Moreover,
\[
\pc-m=\pss-\pc=\dfrac{m(\sigma+2)}{N-2},
\]
so that $p_c$ lies exactly midway between $m$ and $p^{**}$. The explicit pair
$(S,p_c)$ is noted for $\sigma>0$ in the concluding discussion of
\cite[\S3(c)]{IL2}. For $\sigma<0$ the same profile appears in the phase-space
analysis of \cite{ISS26}, where it is identified as the orbit leaving a
critical point of their dynamical system, with the amplitude \eqref{eq:A}, and
where it is observed that it exists for $p>p_c$ (see the remark following
\cite[Lemma~2.1]{ISS26}). The equivalence proved below and the distributional
identity on the whole of $\RN$ are therefore only a short complement to these
observations. {The singular steady state is used as an upper barrier in
\cite[proof of Proposition~3.1]{ISSpeak}, for the unbounded initial values
$K|x|^{-\ac}$ with $K<A$ and $p>\pss$. Its use for bounded data with
critical tail, within the solution framework of \cite{IL1} and in the whole
range $p>\pc$, together with the bounded barriers with critical tail which
are built from it, is the core of the analysis developed here.} We formalize these facts in the following
proposition.
\begin{proposition}\label{prop:steady}
Let $N\ge3$, $m>1$, $\sigma\in(-2,0)$ and $p>m$, and set
$\alpha_c=(\sigma+2)/(p-m)$.  The following are equivalent.
\begin{enumerate}[label=\textup{(\roman*)},leftmargin=2.2em]
\item There is $A>0$ such that $S(x)\eqdef A|x|^{-\alpha_c}$ solves
\begin{equation}\label{eq:stationary}
\Delta S^m+|x|^{\sigma}S^{p}=0\qquad\text{on }\mathbb{R}^N\setminus\{0\}.
\end{equation}
\item $p>p_c=m(\sigma+N)/(N-2)$.
\end{enumerate}
In that case $A$ is unique and given by
\begin{equation}\label{eq:A}
A^{\,p-m}=\alpha_c\,m\,(N-2-\alpha_c m)>0.
\end{equation}
In addition, $S$, $S^m$ and $|x|^{\sigma}S^{p}$ belong to
$L^1_{\mathrm{loc}}(\mathbb{R}^N)$, and \eqref{eq:stationary} holds in
$\mathcal{D}'(\mathbb{R}^N)$, that is, \emph{without} a Dirac mass at the
origin. Finally, the hierarchy \eqref{eq:hierarchy} holds:
$1<p_G<m<p_F<p_c<p^{**}$, with $p_c-m=p^{**}-p_c=m(\sigma+2)/(N-2)$.
\end{proposition}

The distributional identity in Proposition~\ref{prop:steady} holds for every
$p>\pc$. However, $S$ is a stationary weak supersolution on the whole of $\RN$
in the sense of Definition~\ref{def:supersol} only when $p>\pss$, since this
definition requires $S^m\in H^1_{\mathrm{loc}}(\RN)$ (Remark~\ref{rem:pss}).
On $\RN\setminus\{0\}$, $S$ is a classical solution for every $p>\pc$, and this
is all we shall use. Table~\ref{tab:exponents} collects the exponents of the
problem together with their classical semilinear counterparts.

\begin{table}[ht]
\centering
\small
\renewcommand{\arraystretch}{1.9}
\begin{tabular}{>{\raggedright\arraybackslash}p{3.9cm}>{\centering\arraybackslash}p{3.4cm}>{\centering\arraybackslash}p{2.2cm}>{\raggedright\arraybackslash}p{4.0cm}}
\toprule
Exponent & Problem \eqref{eq:CP} & $m=1$, $\sigma=0$ & Role\\
\midrule
$p_G$ & $1-\dfrac{\sigma(m-1)}{2}$ & $1$ & global existence for $p\le p_G$ \cite{IL1}\\
$p_F$ (Fujita) & $m+\dfrac{\sigma+2}{N}$ & $1+\frac2N$ & blow-up of all solutions for $p_G<p\le p_F$ \cite{Qi98,IL1}\\
$r_0+1$ (critical Lebesgue) & $\dfrac{N(p-m)}{\sigma+2}$ & $\dfrac{N(p-1)}{2}$ & small-data global existence \cite{IL1}\\
$\alpha_c$ (second critical) & $\dfrac{\sigma+2}{p-m}=\dfrac{N}{r_0+1}$ & $\dfrac{2}{p-1}$ & decay threshold (this paper)\\
$p_c$ (Serrin type) & $\dfrac{m(\sigma+N)}{N-2}$ & $\dfrac{N}{N-2}$ & existence of $S$\\
$p^{**}$ (Hardy--Sobolev) & $\dfrac{m(N+2\sigma+2)}{N-2}$ & $\dfrac{N+2}{N-2}$ & $\nabla S^m\in L^2_{\mathrm{loc}}$\\
Joseph--Lundgren condition & \eqref{eq:JL}{, i.e.\ $p\ge p_{JL}(\sigma)$} & $p\ge p_{JL}$ & stability of $S$ (Section~\ref{sec:conjectures})\\
\bottomrule
\end{tabular}
\caption{The exponents of problem \eqref{eq:CP} and their classical
semilinear analogues. They are ordered as $1<p_G<m<p_F<p_c<p^{**}$.}
\label{tab:exponents}
\end{table}

\begin{remark}\rm \label{rem:duality}
By \eqref{eq:r0} we have \(r_0+1=\tfrac{N}{\sigma+2}(p-m)=N/\ac\), that is,
\(\ac(r_0+1)=N\). Thus the pointwise critical decay \(\ac\) and the critical
Lebesgue exponent \(r_0+1\) are dual. The profile \(\abs{x}^{-\ac}\) is
precisely \emph{borderline} for membership in \(\Lp{r_0+1}\), since the
integral
\[
\int_{\{\abs{x}>1\}}\abs{x}^{-\ac(r_0+1)}\dx
=\int_{\{\abs{x}>1\}}\abs{x}^{-N}\dx=+\infty
\]
diverges logarithmically, while \(\abs{x}^{-\gamma}\in L^{r_0+1}(\{\abs{x}>1\})\)
for every \(\gamma>\ac\). Consequently, \cite[Theorem~2.10]{IL1}
(Theorem~\ref{thm:IL-super}(a) below) applies to \(\gamma>\ac\) but
\emph{just misses} the exact critical decay \(\gamma=\ac\). This is not a
technical accident. It is exactly the borderline at which the amplitude, and
not merely the decay rate, decides the outcome. The bounded barrier of
Theorem~\ref{thm:barrier} is the object which gives access to this borderline,
and it is the reason why Theorem~\ref{thm:global} is \emph{not} a consequence of
Theorem~\ref{thm:IL-super}(a). Quantitatively, for \(\gamma>\ac\) we have
\begin{equation}\label{eq:Cgamma}
C_\gamma\eqdef\norm{\phi_\gamma}_{r_0+1}
=\Bigl(\omega_N+\omega_N\,\frac{N}{\gamma(r_0+1)-N}\Bigr)^{\frac{1}{r_0+1}}<\infty,
\qquad C_\gamma\xrightarrow[\gamma\downarrow\ac]{}+\infty .
\end{equation}
This is the precise sense in which the small-data theory of \cite{IL1}
degenerates as the decay rate approaches \(\ac\).
\end{remark}

\begin{remark}\rm \label{rem:pss}
A direct computation gives \(\abs{\nabla S^m}^2=(A^m\ac m)^2\abs{x}^{-2\ac m-2}\).
We infer that
\[
\nabla S^m\in L^2\bigl(B(0,1)\bigr)\iff 2\ac m+2<N\iff \ac m<\tfrac{N-2}{2}
\iff p>\pss,
\]
while
\[
\nabla S^m\in L^2\bigl(\RN\setminus B(0,1)\bigr)\iff 2\ac m+2>N\iff p<\pss .
\]
Thus \(S^m\) never belongs to the global energy space, and \(\pss\) is exactly
the exponent at which the local energy of the singular profile becomes finite.
In the semilinear normalization \(m=1\), \(\sigma=0\), one recovers the Sobolev
exponent \(\pss=(N+2)/(N-2)\) and the Serrin exponent \(\pc=N/(N-2)\)
\cite{Serrin65}, while \(\ac=2/(p-1)\) is the critical decay of Lee and Ni
\cite{LeeNi92}. The formulas \eqref{eq:ac}--\eqref{eq:hierarchy} are therefore
the natural quasilinear and weighted extensions of the classical exponents.
The stability of \(S\) is not governed by \(\pss\), but by a larger exponent of
Joseph--Lundgren type \cite{JosephLundgren73}, see \eqref{eq:JL} in
Section~\ref{sec:conjectures}.
\end{remark}

\begin{remark}\rm \label{rem:pss-sobolev}
The exponent \(\pss\) admits a purely variational description, which explains
the coincidence observed in Remark~\ref{rem:pss}. Setting \(v\eqdef u^m\), the
stationary equation \(\Delta u^m+\abs{x}^{\sigma}u^{p}=0\) becomes the
semilinear Hardy--H\'enon equation
\[
\Delta v+\abs{x}^{\sigma}v^{q}=0,\qquad q\eqdef\frac{p}{m}.
\]
Its critical exponent is dictated by the Hardy--Sobolev inequality
\(\bigl(\int_{\RN}\abs{x}^{\sigma}\abs{v}^{q+1}\dx\bigr)^{1/(q+1)}
\le C\norm{\nabla v}_2\), which holds for \(-\sigma\in[0,2)\) with the
scale-invariant exponent \(q+1=2(N+\sigma)/(N-2)\) \cite{CKN84}. This gives
\[
q_S(\sigma)=\frac{N+2+2\sigma}{N-2},
\qquad\text{so that}\qquad p=m\,q_S(\sigma)=\pss .
\]
For \(\sigma=0\) this gives \(\pss=m(N+2)/(N-2)\), the Sobolev exponent of the
porous medium stationary problem. The exponent \(\pss\) coincides with the
critical exponent \(p_S(\sigma)\) which appears in the classification of global
self-similar solutions of \eqref{eq:CP} in \cite{ISS26}, see
Remark~\ref{rem:selfsimilar}.
\end{remark}

\subsubsection{Trapping below the singular steady state}\label{subsec:trapping}
The singular steady state $S$ provides a natural barrier for the flow. Indeed,
$S$ solves \eqref{eq:stationary} classically on $\mathbb{R}^N\setminus\{0\}$
by Proposition~\ref{prop:steady}. It follows that every sub-multiple
$\vartheta S$, $\vartheta\in(0,1]$, is a stationary weak supersolution of
\eqref{eq:CP} in $\RN\setminus\{0\}$ (Lemma~\ref{lem:submultiple}). Since
$\vartheta S$ blows up at the origin, Lemma~\ref{lem:comparison} applies to this
barrier, and we obtain the following result.

\begin{theorem}\label{thm:trap}
Let \(N\ge3\), \(m>1\), \(\sigma\in(-2,0)\), \(p>\pc\), let
\(\vartheta\in(0,1]\), and let \(u_0\in\mathcal A_r\), \(r>r_c\), satisfy
\begin{equation}\label{eq:below}
u_0(x)\le \vartheta\,S(x)=\vartheta A\abs{x}^{-\ac}\qquad\text{for a.e. }x\in\RN.
\end{equation}
Let \(u\) be the solution of \eqref{eq:CP} on \([0,\Tmax(u_0))\). Then
\begin{equation}\label{eq:trapped}
u(t,x)\le \vartheta\,S(x)=\vartheta A\abs{x}^{-\ac}
\qquad\text{for a.e. }x\in\RN\text{ and all }t\in[0,\Tmax(u_0)).
\end{equation}
In particular, for every \(\eps>0\),
\[
\sup_{t\in[0,\Tmax(u_0))}\ \operatorname*{ess\,sup}_{\abs{x}\ge\eps}u(t,x)
\ \le\ \vartheta A\,\eps^{-\ac}<\infty .
\]
Hence \(u\) cannot blow up outside any neighbourhood of the origin, and its
tail is controlled by the critical profile uniformly in time.
\end{theorem}

\begin{remark}\rm\label{rem:origin}
Estimate \eqref{eq:trapped} controls \(u\) everywhere except at the origin,
where \(S=+\infty\). Since \(\Tmax(u_0)<\infty\) forces
\(\limsup_{t\uparrow\Tmax}\norm{u(t)}_{\infty}=\infty\)
(Assumption~\ref{ass:IL}(A)), the trapping bound rules out blow-up away from
\(x=0\). By itself, however, it does not exclude a single-point blow-up at the
location of the singular potential. Indeed, under \eqref{eq:trapped} with
\(\vartheta=1\), the source obeys
\[
\abs{x}^{\sigma}u^{p}\le A^{p}\abs{x}^{\sigma-\ac p}=A^{p}\abs{x}^{-\ac m-2},
\qquad \ac m+2\in(2,N).
\]
Therefore this bound only yields \(\abs{x}^{\sigma}u^p\in L^{q}_{\mathrm{loc}}\)
for \(q<N/(\ac m+2)<N/2\), which is below the threshold \(q>N/2\) needed to
run an \(L^\infty\) smoothing argument \cite{Serrin64}. For
\(\vartheta\in(0,1)\), this difficulty disappears: combining
Theorem~\ref{thm:trap} with the bounded barriers of Theorem~\ref{thm:conj91},
we show in Corollary~\ref{cor:belowS} that such solutions are global and
uniformly bounded. Whether a solution starting below \(S\) itself may blow up
at the origin is open when the Joseph--Lundgren condition \eqref{eq:JL} holds
(Conjecture~\ref{conj:attract}). On the other hand, the naive regularization
of \(S\),
\[
S_\eps(x)\eqdef A\bigl(\abs{x}^2+\eps^{2}\bigr)^{-\ac/2},
\]
is \emph{not} a supersolution. Indeed, a computation based on
\eqref{eq:laplacian} gives
\[
\Delta S_\eps^m+\abs{x}^{\sigma}S_\eps^{p}
=A^{m}\bigl(\abs{x}^{2}+\eps^{2}\bigr)^{-\frac{\ac m}{2}-2}
\Bigl[A^{p-m}\abs{x}^{\sigma}\bigl(\abs{x}^2+\eps^2\bigr)^{\frac{2-\sigma}{2}}
-\ac m\bigl((N-2-\ac m)\abs{x}^{2}+N\eps^{2}\bigr)\Bigr].
\]
The bracket tends to \(+\infty\) as \(x\to0\), because \(\sigma<0\). Hence
\(S_\eps\) is a strict \emph{sub}solution near the origin, for every
\(\eps>0\). We conclude that a bounded barrier with critical tail cannot be
produced by this naive smoothing of \(S\), and that it has to be built. This is
done in Theorem~\ref{thm:barrier} below, by matching \(S\) to an inner profile
whose Laplacian reproduces the Hardy weight \(\abs{x}^{\sigma}\) exactly. It is
also done, {with a larger tail amplitude}, in Theorem~\ref{thm:conj91}, by matching
\(\vartheta S\) to regular radial steady states.
\end{remark}
\subsubsection{A bounded barrier with critical tail, and global existence}
\label{subsec:global}
We now construct a bounded stationary weak supersolution of \eqref{eq:CP}
with exactly critical tail. As explained above, the idea is to glue the
singular steady state $\mu S$ to the inner profile $\Psi$ of
Proposition~\ref{prop:inner}, whose Laplacian reproduces the Hardy weight
$|x|^{\sigma}$ exactly. The resulting barrier is described in the following
theorem.
\begin{theorem}
\label{thm:barrier}
Let \(N\ge3\), \(m>1\), \(\sigma\in(-2,0)\) and \(p>\pc\). Let \(A\), \(\ac\)
be as in Proposition~\ref{prop:steady}, let \(c\), \(w\), \(\Psi\) be as in
Proposition~\ref{prop:inner}, and set
\begin{equation}\label{eq:mu}
\mu\eqdef\min\Bigl\{1,\ \frac{(c/2)^{1/m}}{A\,2^{\ac}}\Bigr\}\in(0,1],
\qquad
\rho\eqdef\inf\Bigl\{\,\varrho\in(0,1]:\ \mu^{m}A^{m}\varrho^{-\ac m}\le w(\varrho)\Bigr\}.
\end{equation}
Then \(\rho\in(0,\tfrac12]\), and the function
\begin{equation}\label{eq:barrier-def}
\ol u(x)\eqdef
\begin{cases}
\Psi(x)=\Bigl[c\bigl(1-\tfrac12\abs{x}^{\sigma+2}\bigr)\Bigr]^{1/m},
&\abs{x}<\rho,\\[2mm]
\mu\,S(x)=\mu A\abs{x}^{-\ac}, &\abs{x}\ge\rho,
\end{cases}
\end{equation}
is continuous on \(\RN\), positive and radially nonincreasing. It is bounded,
with \(\norm{\ol u}_{\infty}=\ol u(0)=c^{1/m}\), and it satisfies
\begin{equation}\label{eq:barrier}
\overline{u}(x)=\mu A|x|^{-\alpha_c}\quad\text{for }|x|\ge\tfrac12,
\qquad \overline{u}\le\|\overline{u}\|_{\infty}<\infty .
\end{equation}
Moreover, \(\ol u\) belongs to \(\Lp{r+1}\) for every \(r>r_0\), and it is a
stationary weak supersolution of \eqref{eq:CP} in \(\RN\) in the sense of
Definition~\ref{def:supersol}.
\end{theorem}

\begin{theorem}
\label{thm:global}
Let \(N\ge3\), \(m>1\), \(\sigma\in(-2,0)\), \(p>\pc\), and let \(\ol u\),
\(\mu\), \(c\) be as in Theorem~\ref{thm:barrier}. Let \(\vartheta\in(0,1]\)
and let \(u_0\in\mathcal A_r\) for some \(r>r_c\) with
\begin{equation}\label{eq:below-barrier}
u_0(x)\le\vartheta\,\ol u(x)\qquad\text{for a.e. }x\in\RN.
\end{equation}
Then the solution of \eqref{eq:CP} is \emph{global}, \(\Tmax(u_0)=\infty\),
and
\begin{equation}\label{eq:global-bound}
u(t,x)\le\vartheta\,\ol u(x)\le\vartheta\, c^{1/m}
\qquad\text{for all }t\ge0\text{ and a.e. }x\in\RN .
\end{equation}
In particular, let \(L\in(0,\mu A]\) and let \(u_0\) satisfy
\(0\le u_0\le L\phi\), where \(\phi=\phi_{\ac}=\min\{1,\abs{x}^{-\ac}\}\).
(The datum \(u_0=L\phi\) has \emph{exactly critical decay}, namely
\(\lim_{\abs{x}\to\infty}\abs{x}^{\ac}u_0(x)=L\).) Then the solution is global
and obeys
\begin{equation}\label{eq:global-bound-L}
u(t,x)\le\frac{L}{\mu A}\,\ol u(x)\qquad\text{for all }t\ge0
\text{ and a.e. }x\in\RN .
\end{equation}
The same conclusion holds for every \(\gamma\ge\ac\) and every \(u_0\) with
\(0\le u_0\le L\phi_\gamma\), \(L\le\mu A\).
\end{theorem}

\begin{remark}\rm
\label{rem:compare-IL}
Theorem~\ref{thm:global} is not contained in Theorem~\ref{thm:IL-super}(a).
Indeed, for \(\gamma=\ac\) the datum \(L\phi\) does \emph{not} belong to
\(\Lp{r_0+1}\) (Remark~\ref{rem:duality}). Hence the critical Lebesgue
smallness condition \(\norm{u_0}_{r_0+1}\le C_0\) is not merely violated, it is
meaningless. Conversely, for \(\gamma>\ac\), Theorem~\ref{thm:IL-super}(a)
gives global existence under an integral smallness condition, whereas
Theorem~\ref{thm:global} gives it under a pointwise domination by \(\ol u\).
The latter also yields the uniform-in-time pointwise decay estimate
\eqref{eq:global-bound}, which the \(L^{r}\) framework does not provide. The
two statements are therefore complementary, and only the barrier reaches the
exact borderline.
\end{remark}

\begin{remark}\rm \label{rem:scaling}
The barrier is not unique. By \eqref{eq:scaling}, if \(\ol u\) is a stationary
supersolution, then so is \(\ol u_\lambda(x)\eqdef\lambda^{\ac}\ol u(\lambda x)\)
for every \(\lambda>0\). Indeed, since \(\ac p-\sigma=\ac m+2\) by \eqref{eq:ac},
we have
\[
\Delta\ol u_\lambda^{\,m}(x)+\abs{x}^{\sigma}\ol u_\lambda^{\,p}(x)
=\lambda^{\ac m+2}\Bigl[\Delta\ol u^{\,m}+\abs{y}^{\sigma}\ol u^{\,p}\Bigr]_{y=\lambda x}
\le0 .
\]
This family has the \emph{same} tail amplitude \(\mu A\), because the critical
profile is scale invariant, while
\(\norm{\ol u_\lambda}_{\infty}=\lambda^{\ac}c^{1/m}\) and the matching radius
is \(\rho/\lambda\). Thus the amplitude \(\mu A\) is intrinsic to the
construction, and it cannot be improved by rescaling. The loss of the factor
\(\mu\) comes from the explicit inner profile \(w\). It is removed in
Theorem~\ref{thm:conj91}, where \(w\) is replaced by regular radial steady
states, and this yields \(L_\ast(\ac)\ge A\). The identification of
\(L_\ast(\ac)\) itself is the subject of Conjecture~\ref{conj:threshold}. The
same scaling, combined with the large values of the barrier near the origin,
is the key to Corollary~\ref{cor:belowS} below.
\end{remark}
%======================================================
\subsubsection{Blow-up for slow decay, at any amplitude}
We now turn to the blow-up side. The key observation is that, for data
decaying more slowly than the critical rate $\alpha_c$, the energy of the
piece of $u_0$ sitting on a far annulus becomes negative, whatever the
amplitude. Combined with the negative energy criterion of
\cite[Theorem~2.9]{IL1} and with the comparison principle, this forces
finite-time blow-up.
\begin{theorem}\label{thm:slow}
Let \(N\ge3\), \(m>1\), \(\sigma\in(-2,0)\) and \(p>\pF\). Let
\(u_0\in\mathcal A_r\) for some \(r>r_c\), and assume that, for some
\(\gamma\in(0,\ac)\),
\begin{equation}\label{eq:slowtail}
\liminf_{\abs{x}\to\infty}\ \abs{x}^{\gamma}u_0(x)\ =\ \kappa_0\ \in(0,\infty] .
\end{equation}
Then \(\Tmax(u_0)<\infty\): the solution blows up in finite time, \emph{no
matter how small \(\kappa_0\) is}.
\end{theorem}
We now consider the critical decay $\gamma=\alpha_c$. In this borderline case,
the decay rate alone does not determine the behaviour of the solution, and an
amplitude threshold appears. The following theorem shows that blow-up occurs
in finite time above the explicit threshold $L_0$.
\begin{theorem}
\label{thm:blowup-crit}
Let \(N\ge3\), \(m>1\), \(\sigma\in(-2,0)\) and \(p>\pF\). Define
\begin{equation}\label{eq:L0}
L_0\eqdef 2^{\ac}M_0
=2^{\ac}\left(\frac{(m+p)\,I_1}{2m\,I_2}\right)^{\frac{1}{p-m}}\in(0,\infty),
\end{equation}
with \(I_1,I_2\) as in \eqref{eq:I1I2}. If \(u_0\in\mathcal A_r\), \(r>r_c\),
satisfies
\begin{equation}\label{eq:crittail}
\liminf_{\abs{x}\to\infty}\ \abs{x}^{\ac}u_0(x)\ >\ L_0 ,
\end{equation}
then \(\Tmax(u_0)<\infty\). In particular, for \(\phi=\phi_{\ac}\), the
solution with datum \(u_0=L\phi\) blows up in finite time whenever \(L>L_0\).
\end{theorem}
\begin{remark}\rm
\label{rem:whycritical}
The three cases of Lemma~\ref{lem:energy-annulus} are transparent. For
\(\gamma<\ac\), the reactive term of the energy dominates the diffusive one on
far annuli by a positive power of \(\lambda\). Blow-up is then forced by
placing the mass far away, at any amplitude: the potential decays, but the
annuli grow faster. For \(\gamma>\ac\), the diffusive term dominates on far
annuli, and no blow-up can be extracted in this way, consistently with the
global existence of small data. For \(\gamma=\ac\), the two homogeneities
balance exactly, the energy of the annulus piece is scale invariant, and only
the amplitude decides. This is precisely the mechanism by which \(\ac\), and no
other exponent, is the second critical exponent. It also shows that the
constant \(L_0\) of \eqref{eq:L0} may be optimized over the profile \(\zeta\),
but not by rescaling. Since the conclusion of Theorem~\ref{thm:blowup-crit}
holds for each admissible \(\zeta\), one may replace \(L_0\) throughout by
\begin{equation}\label{eq:L0star}
L_0^{\ast}\eqdef 2^{\ac}\inf\left\{
\left(\frac{(m+p)I_1(\zeta)}{2mI_2(\zeta)}\right)^{\frac{1}{p-m}}:\
\zeta\ \text{as in \eqref{eq:zeta}}\right\}\le L_0 ,
\end{equation}
the infimum being a variational problem of Hardy--Sobolev type on an annulus.
The factor \(2^{\ac}\) can be improved as well. If \(\zeta\) is replaced by
\(\tilde\zeta(y)\eqdef\abs{y}^{-\ac}\zeta(y)\), which is again admissible in
\eqref{eq:zeta}, then
\(\kappa_1\lambda^{-\ac}\tilde\zeta(x/\lambda)\le\kappa_1\abs{x}^{-\ac}\), and
the computation of Lemma~\ref{lem:energy-annulus} applies verbatim to
\(\tilde\zeta\). We obtain blow-up whenever the left-hand side of
\eqref{eq:crittail} exceeds the infimum of \(M_0(\zeta)\) over all admissible
\(\zeta\) with \(\zeta(y)\le\abs{y}^{-\ac}\). Since \(M_0(t\zeta)=t^{-1}M_0(\zeta)\)
for \(t>0\), and since \(2^{-\ac}\zeta\) satisfies this constraint for every
admissible \(\zeta\), this infimum is at most \(L_0^\ast\).
\end{remark}

\begin{remark}\rm \label{rem:largedata}
Theorem~\ref{thm:IL-super}(b), applied to \(v_{M,1}\) with \(M>M_0\), shows via
\eqref{eq:Esign} that \emph{some} large data blow up regardless of their
decay. The content of Theorem~\ref{thm:blowup-crit} is therefore not the mere
existence of blowing-up data at \(\gamma=\ac\). It is the fact that, at the
critical decay, blow-up is triggered by the tail amplitude alone, whatever the
behaviour of \(u_0\) on bounded sets. Combined with
Theorems~\ref{thm:global} and \ref{thm:conj91}, which provide global solutions
with exactly critical tail, it isolates an amplitude dichotomy at the
borderline. Neither statement alone would suffice.
\end{remark}
\subsubsection{The second critical exponent and the amplitude threshold}\label{subsec:sce}
%======================================================
We now assume that $m>1$, $\sigma\in(-2,0)$, $N\ge3$ and $p>p_F$, and we
consider initial data with a prescribed algebraic tail,
\begin{equation}\label{eq:tail}
u_0\in L^\infty(\mathbb{R}^N),\qquad u_0\ge0,\qquad
\lim_{|x|\to\infty}|x|^{\gamma}u_0(x)=L\in(0,\infty).
\end{equation}
Our aim is to determine the critical decay rate $\gamma=\alpha_c$ which
separates the two possible behaviours of the solution to \eqref{eq:CP},
namely global existence for $\gamma>\alpha_c$ (fast decay) and finite-time
blow-up for $\gamma<\alpha_c$ (slow decay). We also want to describe the
borderline $\gamma=\alpha_c$, where an amplitude threshold on $L$ is expected.

The results of the previous subsections settle the two extreme regimes. On the
one hand, Theorem~\ref{thm:slow} shows that, for $\gamma<\alpha_c$, every datum
with a tail of size $|x|^{-\gamma}$ produces finite-time blow-up, whatever the
amplitude. On the other hand, the small-data theorem of \cite{IL1}
(Theorem~\ref{thm:IL-super}(a)) shows that, for $\gamma>\alpha_c$, data which
are small in the weighted norm $\sup_x\max\{1,|x|^{\gamma}\}u_0(x)$ are global
(see Theorem~\ref{thm:sce}(ii) below). Moreover, for $p>p_c$,
Theorem~\ref{thm:barrier} provides a bounded barrier at the exact critical
decay $\gamma=\alpha_c$, which yields global existence for amplitudes
$L\le\mu A$.

The borderline case $\gamma=\alpha_c$ is the most delicate part of the
analysis. At the critical decay, the diffusive and reactive homogeneities
balance exactly, and the outcome depends on the size of the amplitude rather
than on the decay rate. The next results describe this threshold.
\begin{theorem}\label{thm:sce}
Let \(N\ge3\), \(m>1\), \(\sigma\in(-2,0)\) and \(p>\pF\), and let
\(\ac=(\sigma+2)/(p-m)=N/(r_0+1)\).
\begin{enumerate}[label=\textup{(\roman*)},leftmargin=2.2em]
\item \emph{(Slow decay: universal blow-up.)} If \(0<\gamma<\ac\) and
\(u_0\in\mathcal A_r\) satisfies \eqref{eq:tail}, or merely
\(\liminf_{\abs{x}\to\infty}\abs{x}^{\gamma}u_0(x)>0\), then
\(\Tmax(u_0)<\infty\), whatever the amplitude \(L>0\).
\item \emph{(Fast decay: global existence for small data.)} Let \(\gamma>\ac\)
and set \(\eps_\gamma\eqdef C_0/C_\gamma>0\), with \(C_0\) as in
Theorem~\ref{thm:IL-super}(a) and \(C_\gamma\) as in \eqref{eq:Cgamma}. If
\(u_0\ge0\) satisfies
\begin{equation}\label{eq:smallweighted}
\norm{u_0}_{[\gamma]}=\sup_{x\in\RN}\max\{1,\abs{x}^{\gamma}\}\,u_0(x)
\ \le\ \eps_\gamma,
\end{equation}
that is, \(0\le u_0\le\eps_\gamma\phi_\gamma\), then the solution of
\eqref{eq:CP} is global and \(u\in L^{\infty}((0,\infty),\Lp{r})\) for every
\(r\in[r_0+1,\infty)\). If moreover \(p>\pc\) and \(u_0\le\vartheta\ol u\) for
some \(\vartheta\in(0,1]\), the solution obeys the uniform pointwise bound
\eqref{eq:global-bound}.
\end{enumerate}
Consequently, \(\ac\) is the second critical exponent for \eqref{eq:CP} in the
sense of Lee and Ni \cite{LeeNi92}. For \(\gamma<\ac\), every datum with tail
\(\abs{x}^{-\gamma}\) produces finite-time blow-up. For \(\gamma>\ac\), small
data with that tail are global, while large data still blow up
\textup{(}Remark~\textup{\ref{rem:largedata}}\textup{)}.
\end{theorem}
The smallness condition \eqref{eq:smallweighted} couples the amplitude and the
decay, and this coupling is necessary. Smallness of the tail amplitude \(L\) in
\eqref{eq:tail}, together with smallness of \(\norm{u_0}_{\infty}\), does
\emph{not} suffice, no matter how fast the prescribed decay is.

\begin{proposition}
\label{prop:nodecayonly}
Let \(N\ge3\), \(m>1\), \(\sigma\in(-2,0)\) and \(p>\pF\). For every
\(\gamma>0\), every \(\eps>0\) and every \(L\in(0,\eps]\), there exists
\(u_0\in\mathcal A_r\) \textup{(}for every \(r\) large enough\textup{)} with
\[
\norm{u_0}_{\infty}=\eps,\qquad
\lim_{\abs{x}\to\infty}\abs{x}^{\gamma}u_0(x)=L,\qquad
\Tmax(u_0)<\infty .
\]
\end{proposition}

At the borderline \(\gamma=\ac\), the outcome is decided by the amplitude, and,
for \(p>\pc\), a genuine threshold exists. It is convenient to record the whole
picture at once, by testing along the rays \(L\mapsto L\phi_\gamma\).

\begin{theorem}\label{thm:threshold}
Let \(N\ge3\), \(m>1\), \(\sigma\in(-2,0)\) and \(p>\pF\). For \(\gamma>0\) set
\[
\mathcal G_\gamma\eqdef\bigl\{L>0:\ \Tmax(L\phi_\gamma)=\infty\bigr\},
\qquad
L_\ast(\gamma)\eqdef\sup\mathcal G_\gamma\in[0,\infty),
\]
with the convention \(\sup\emptyset=0\). Then the following hold.
\begin{enumerate}[label=\textup{(\roman*)},leftmargin=2.2em]
\item \emph{(Dichotomy along each ray.)} If \(0<L<L_\ast(\gamma)\), then
\(\Tmax(L\phi_\gamma)=\infty\), and if \(L>L_\ast(\gamma)\), then
\(\Tmax(L\phi_\gamma)<\infty\).
\item \emph{(Finiteness and monotonicity.)} \(L_\ast(\gamma)\le 2^{\gamma}M_0<\infty\)
for every \(\gamma>0\), with \(M_0\) as in \eqref{eq:M0}, and
\(\gamma\mapsto L_\ast(\gamma)\) is nondecreasing on \((0,\infty)\).
\item \emph{(Slow decay.)} \(L_\ast(\gamma)=0\) for every \(\gamma\in(0,\ac)\).
\item \emph{(Fast decay.)} \(L_\ast(\gamma)\ge\eps_\gamma>0\) for every
\(\gamma>\ac\), with \(\eps_\gamma\) as in Theorem~\ref{thm:sce}(ii).
\item \emph{(Critical decay.)} If in addition \(p>\pc\), then
\begin{equation}\label{eq:Lstar}
L_\ast\eqdef L_\ast(\ac)\in\bigl[\mu A,\ L_0\bigr]\subset(0,\infty),
\end{equation}
with \(\mu A\) as in Theorem~\ref{thm:barrier} and \(L_0\) as in
\eqref{eq:L0}. Moreover, for \(0<L\le\mu A\), the global solution satisfies
\(u(t,\cdot)\le\frac{L}{\mu A}\ol u\) for all \(t\ge0\). In particular, it
retains the critical decay uniformly in time.
\end{enumerate}
\end{theorem}
The lower bound in \eqref{eq:Lstar} is improved to \(L_\ast(\ac)\ge A\) in
Theorem~\ref{thm:conj91} and Corollary~\ref{cor:sandwich} below.
\begin{corollary}\label{cor:char}
Under the hypotheses of Theorem~\ref{thm:threshold},
\[
\ac=\sup\bigl\{\gamma>0:\ L_\ast(\gamma)=0\bigr\}
=\inf\bigl\{\gamma>0:\ L_\ast(\gamma)>0\bigr\} .
\]
Thus \(\ac\) is exactly the decay rate at which the amplitude threshold
function ceases to vanish. When \(p>\pc\), the function \(L_\ast\) is
moreover positive and finite at \(\ac\) itself.
\end{corollary}

The positivity of \(L_\ast(\ac)\) actually holds in the whole supercritical
range, including \(\pF<p\le\pc\), where no singular steady state is available.
In that range the barriers are provided by the forward self-similar solutions
constructed by Iagar, S\'anchez and Sarri\'on-Pedralva \cite{ISS26}. For
\(b>0\), we denote by \(f(\cdot;b)\) the profile of \cite{ISS26}, namely the
solution of the profile equation \eqref{eq:profile} of
Section~\ref{sec:conjectures} with \(f(0;b)=b\) and the
local behaviour \cite[Eq.~(1.10)]{ISS26} at \(\xi=0\). {We write \(b\) for
the central value of the profile in order to avoid any confusion with the
amplitude \(A\) of the singular steady state.} By
\cite[Theorem~1.2]{ISS26}, if \(\pF<p<\pss\), there is \(b_\ast>0\) such that,
for every \(b\in(0,b_\ast)\),
\begin{equation}\label{eq:ISStail}
\lim_{\xi\to\infty}\xi^{\ac}f(\xi;b)=c_\infty(b)\in(0,\infty),
\end{equation}
and \eqref{eq:ISStail} holds for every \(b>0\) if \(p\ge\pss\). We set
\(\mathcal B_{\rm alg}\eqdef(0,b_\ast)\) if \(\pF<p<\pss\) and
\(\mathcal B_{\rm alg}\eqdef(0,\infty)\) if \(p\ge\pss\).

\begin{theorem}\label{thm:allp}
Let \(N\ge3\), \(m>1\), \(\sigma\in(-2,0)\) and \(p>\pF\). Then
\[
0<\sup_{b\in\mathcal B_{\rm alg}}c_\infty(b)\le L_\ast(\ac)\le L_0 .
\]
More precisely, for every \(b\in\mathcal B_{\rm alg}\), the profile
\(f(\cdot;b)\) is positive on \([0,\infty)\), and every datum
\(u_0\in\mathcal A_r\), \(r>r_c\), with \(u_0\le U(\tau,\cdot)\) for some
\(\tau>0\), where \(U(t,x)=t^{-\alpha}f(\abs{x}t^{-\beta};A)\), generates a
global solution satisfying \(u(t)\le U(t+\tau)\le\tau^{-\alpha}M_b\) for all
\(t\ge0\), where \(M_b\eqdef\sup_{\xi\ge0}f(\xi;b)<\infty\) and \(\alpha\),
\(\beta\) are given by \eqref{eq:U}. Consequently, for
every \(p>\pF\), the zero set of \(L_\ast\) is exactly \((0,\ac)\), that is,
\[
\ac=\min\bigl\{\gamma>0:\ L_\ast(\gamma)>0\bigr\}.
\]
\end{theorem}

The proof is given in Section~\ref{sec:conjectures}. Its new ingredient is a
positivity property of nonnegative self-similar profiles,
Lemma~\ref{lem:positivity}, which is the self-similar form of the persistence
of positivity for porous medium type equations. It shows that a profile which
is positive near the origin and near infinity cannot vanish in between. In the
range \(\pF<p\le\pc\), Theorem~\ref{thm:allp} replaces the singular barriers of
Theorems~\ref{thm:barrier} and \ref{thm:conj91}. It gives no explicit value of
the lower bound, and the identification of \(L_\ast(\ac)\) remains open there
(Open Problem~\ref{prob:pcrange}).

\begin{remark}\rm \label{rem:henon}
It is worth comparing Theorem~\ref{thm:threshold} with the H\'enon case
\(\sigma>0\) treated in \cite{IL2}. There, for \(p\ge\pG\), the sharp
threshold for \emph{non-existence} of any nonnegative solution is
\(\gamma<\sigma/(p-1)\). For \(\sigma>0\) and \(m>1\), this is the only
relevant range, since then \(\pG=1+\sigma(1-m)/2<1\). The scale-invariant
decay \((\sigma+2)/(p-m)=\ac\) appears as a non-existence threshold only in the
fast-diffusion range \(1<p<\pG\), which requires \(m<1\)
\cite[Eqs.~(1.6)--(1.7)]{IL2}. In the Hardy regime \(\sigma<0\) studied here,
the quantity \(\sigma/(p-1)\) is negative. Hence no decay rate produces
non-existence, and the local theory of \cite{IL1} applies to all bounded data
with algebraic decay. The borderline phenomenon migrates to the
scale-invariant decay \(\ac\), and the alternative it decides is global
existence versus finite-time blow-up, rather than existence versus
non-existence. What the two regimes do share is the singular steady state
\(S\) and the exponent \(\pc\) of Proposition~\ref{prop:steady}. The profile
\(A\abs{x}^{-\ac}\) is exhibited in \cite[\S3(c)]{IL2} with the amplitude
\(A^{p-m}=m(\sigma+2)(N-2)(p-\pc)/(p-m)^2\), which coincides with \eqref{eq:A}.
The same explicit formula is written for \(\sigma\in(-2,0)\) in the remark
following \cite[Lemma~2.1]{ISS26}.
\end{remark}

\subsubsection{\texorpdfstring{{An improved lower bound}}{An improved lower bound} at the critical decay}\label{subsec:sharp}

The lower bound \(\mu A\) in Theorem~\ref{thm:threshold}(v) is not optimal.
Replacing the explicit inner profile \(w\) of Proposition~\ref{prop:inner} by
regular radial steady states (Lemma~\ref{lem:A}), we obtain the following
{improved bound, with the amplitude $A$ of the singular steady state}.

\begin{theorem}\label{thm:conj91}
Let \(N\ge3\), \(m>1\), \(\sigma\in(-2,0)\) and \(p>\pc\). For every
\(\vartheta\in(0,1)\) there is a bounded, radially nonincreasing stationary
weak supersolution \(\ol u_\vartheta\) of \eqref{eq:CP} in \(\RN\), in the
sense of Definition~\ref{def:supersol}, such that
\begin{equation}\label{eq:barrier-theta}
\ol u_\vartheta(x)=\vartheta A\abs{x}^{-\ac}\quad\text{for }\abs{x}\ge1,
\qquad
\ol u_\vartheta(x)\ge\vartheta A\quad\text{for }\abs{x}\le1 .
\end{equation}
Consequently, if \(u_0\in\mathcal A_r\), \(r>r_c\), satisfies
\(u_0\le\ol u_\vartheta\) a.e.\ in \(\RN\), then \(\Tmax(u_0)=\infty\) and
\begin{equation}\label{B-Infty}
u(t,x)\le\ol u_\vartheta(x)\le\norm{\ol u_\vartheta}_\infty
\qquad\text{for all }t\ge0\text{ and a.e. }x\in\RN .
\end{equation}
In particular, for every \(L\in(0,A)\), the solution with datum
\(L\phi_{\ac}\) is global and uniformly bounded in time.
\end{theorem}

\begin{corollary}\label{cor:sandwich}
Let \(N\ge3\), \(m>1\), \(\sigma\in(-2,0)\) and \(p>\pc\). Then
\[
A\le L_\ast(\ac)\le L_0 .
\]
In particular, \(A\le L_0\) for every cut-off \(\zeta\) as in \eqref{eq:zeta},
that is, \(A\le L_0^{\ast}\) in the notation of \eqref{eq:L0star}. Moreover,
\(L_\ast(\gamma)\ge A\) for every \(\gamma\ge\ac\).
\end{corollary}

\begin{figure}[ht]
\centering
\begin{tikzpicture}
\begin{loglogaxis}[
  width=12.5cm, height=7.2cm,
  xmin=0.01, xmax=400, ymin=0.05, ymax=1.6,
  xlabel={$|x|$}, ylabel={}]
\addplot[very thick,blue!70!black] coordinates {(0.01,0.57002) (0.012115,0.56996) (0.014678,0.56989) (0.017783,0.56981) (0.021544,0.56971) (0.026102,0.56959) (0.031623,0.56944) (0.038312,0.56927) (0.046416,0.56906) (0.056234,0.5688) (0.068129,0.56848) (0.08254,0.56811) (0.1,0.56765) (0.12115,0.5671) (0.14678,0.56643) (0.17783,0.56562) (0.21544,0.56465) (0.26102,0.56347) (0.31623,0.56206) (0.38312,0.56035) (0.46416,0.55831) (0.56234,0.55585) (0.68129,0.5529) (0.8254,0.54938) (1,0.54519) (1,0.54514) (1.3707,0.49074) (1.8789,0.44178) (2.5754,0.3977) (3.5302,0.35802) (4.839,0.3223) (6.6329,0.29014) (9.0919,0.26119) (12.462,0.23513) (17.083,0.21167) (23.416,0.19055) (32.096,0.17154) (43.995,0.15442) (60.305,0.13901) (82.662,0.12514) (113.31,0.11266) (155.31,0.10142) (212.89,0.091296) (291.82,0.082187) (400,0.073986)};
\addplot[thick,dashed,orange!85!black] coordinates {(1,0.54519) (1.1561,0.5415) (1.3366,0.53731) (1.5453,0.53255) (1.7865,0.52718) (2.0654,0.52111) (2.3878,0.5143) (2.7606,0.50667) (3.1916,0.49818) (3.6898,0.48875) (4.2659,0.47834) (4.9318,0.46692) (5.7018,0.45446) (6.5919,0.44095) (7.621,0.42642) (8.8107,0.41088) (10.186,0.39441) (11.776,0.37707) (13.615,0.35896) (15.74,0.34021) (18.198,0.32094) (21.038,0.30127) (24.323,0.28136) (28.12,0.26132) (32.51,0.24126) (37.585,0.22129) (43.453,0.20147) (50.236,0.18182) (58.079,0.16231) (67.146,0.14287) (77.628,0.12328) (89.747,0.10318) (103.76,0.081779) (119.96,0.057059) (138.68,0.016894) (139.39,0.013759) (140.09,0.0097048) (140.51,0.0061286) (140.72,0.003062)};
\addplot[thick,densely dotted,black] coordinates {(0.01,2.5303) (0.015551,2.184) (0.024183,1.8851) (0.037606,1.6271) (0.05848,1.4044) (0.090942,1.2122) (0.14142,1.0463) (0.21992,0.90313) (0.342,0.77953) (0.53183,0.67284) (0.82704,0.58076) (1.2861,0.50128) (2,0.43267) (3.1102,0.37346) (4.8365,0.32235) (7.5212,0.27823) (11.696,0.24016) (18.188,0.20729) (28.284,0.17892) (43.984,0.15443) (68.399,0.1333) (106.37,0.11505) (165.41,0.099309) (257.22,0.085717) (400,0.073986)};
\addplot[thin,gray] coordinates {(0.01,2.8114) (0.015551,2.4267) (0.024183,2.0946) (0.037606,1.8079) (0.05848,1.5605) (0.090942,1.3469) (0.14142,1.1626) (0.21992,1.0035) (0.342,0.86614) (0.53183,0.7476) (0.82704,0.64529) (1.2861,0.55698) (2,0.48075) (3.1102,0.41496) (4.8365,0.35817) (7.5212,0.30915) (11.696,0.26684) (18.188,0.23032) (28.284,0.1988) (43.984,0.17159) (68.399,0.14811) (106.37,0.12784) (165.41,0.11034) (257.22,0.095242) (400,0.082207)};
\draw[gray,densely dashed] (axis cs:1,0.05) -- (axis cs:1,1.6);
\node[blue!70!black,anchor=south west] at (axis cs:0.012,0.585) {$\ol u_\vartheta$};
\node[gray!70!black,anchor=south west] at (axis cs:0.25,0.99) {$S$};
\node[black,anchor=north east] at (axis cs:0.045,1.3) {$\vartheta S$};
\node[orange!85!black,anchor=south west,font=\small] at (axis cs:9,0.43) {regular steady state};
\node[anchor=north west,font=\small,gray!70!black] at (axis cs:1.08,1.5) {matching radius $|x|=1$};
\end{loglogaxis}
\end{tikzpicture}
\caption{The barrier $\ol u_\vartheta$ of Theorem~\ref{thm:conj91}, computed
numerically for $N=3$, $m=2$, $\sigma=-1$, $p=5$ (so that $p_c=4<p<p^{**}=6$)
and $\vartheta=0.9$, in logarithmic scales: $\ol u_\vartheta$ (solid blue), the
continuation of the rescaled regular steady state (dashed orange),
$\vartheta S$ (dotted) and $S$ (thin gray). The barrier is the rescaled regular
steady state for $|x|<1$ and coincides with $\vartheta S$ for $|x|\ge1$. Beyond
the matching radius, the regular steady state lies above $\vartheta S$ and
vanishes at a finite radius, so that $\ol u_\vartheta$ is the minimum of the
two branches, as required by the sign condition on the jump of the normal
derivative.}
\label{fig:barrier}
\end{figure}
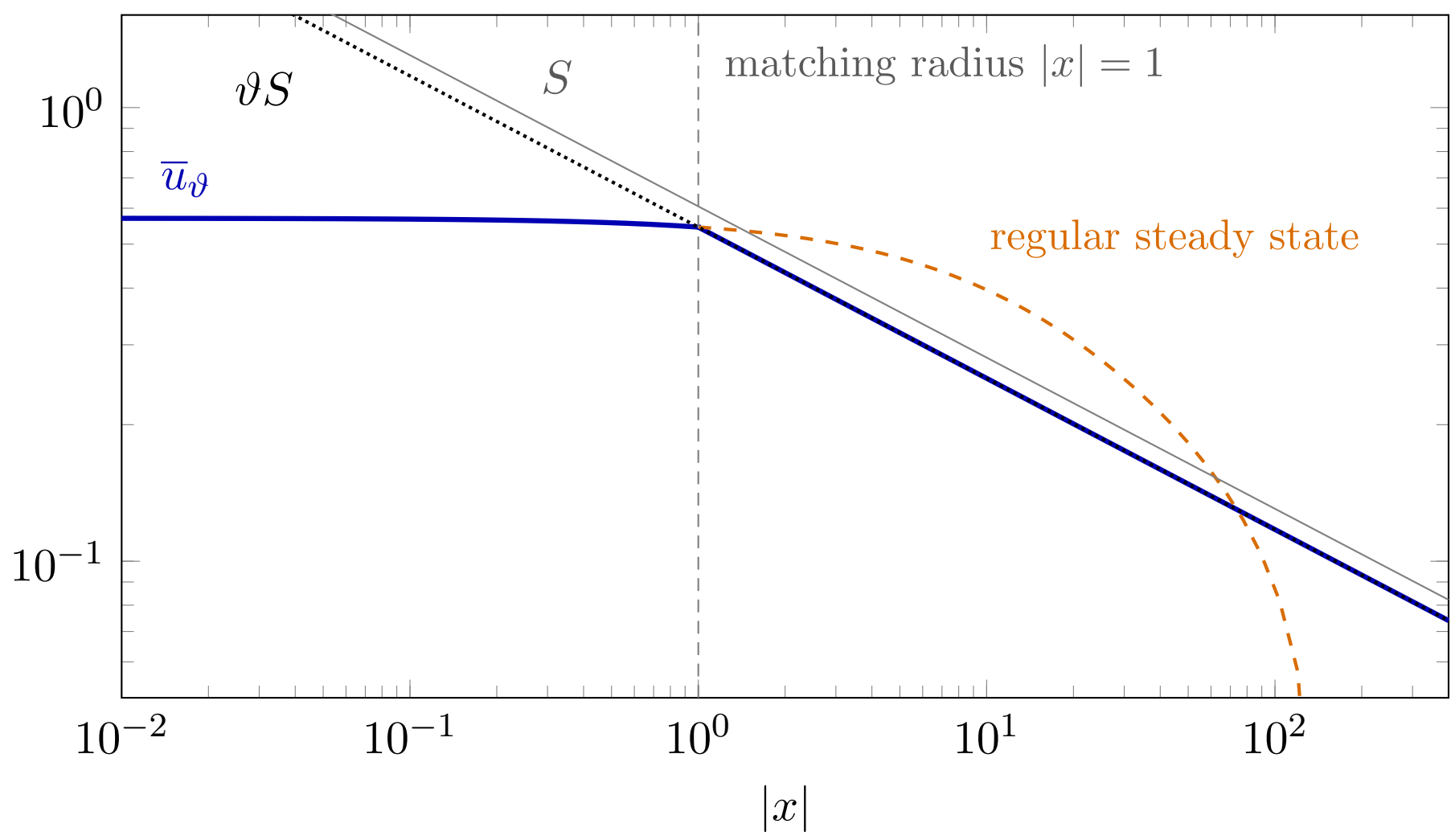

The barrier of Theorem~\ref{thm:conj91} is illustrated in
Figure~\ref{fig:barrier}. The construction also works with \(\vartheta=1\),
provided the Emden--Fowler function of Lemma~\ref{lem:A} reaches the level
\(A^m\). We show that this is the case whenever \(p\) does not satisfy the
Joseph--Lundgren condition, that is, whenever either \(p\le\pss\), or
\(p>\pss\) and
\begin{equation}\label{eq:JLfails}
\bigl(N-2-2\ac m\bigr)^2<4(\sigma+2)(N-2-\ac m)
\end{equation}
(compare with \eqref{eq:JL} in Section~\ref{sec:conjectures}).

\begin{proposition}
\label{rem:LequalA}
Let \(N\ge3\), \(m>1\), \(\sigma\in(-2,0)\) and \(p>\pc\), and assume that
either \(p\le\pss\), or \(p>\pss\) and \eqref{eq:JLfails} holds. Then there is
a bounded, radially nonincreasing stationary weak supersolution \(\ol u_1\) of
\eqref{eq:CP} in \(\RN\) such that
\[
\ol u_1(x)=A\abs{x}^{-\ac}\quad\text{for }\abs{x}\ge1,
\qquad
\ol u_1(x)\ge A\quad\text{for }\abs{x}\le1 .
\]
Consequently, every \(u_0\in\mathcal A_r\), \(r>r_c\), with
\(u_0\le\ol u_1\) generates a global solution bounded by
\(\norm{\ol u_1}_\infty\). In particular, \(A\in\mathcal G_{\ac}\).
\end{proposition}

Combining the trapping theorem with the bounded barriers, we finally obtain
global existence for all bounded data lying below a proper fraction of the
singular steady state. This is the quasilinear and singular-weight
counterpart of a result of Wang \cite{Wang93} for \(\partial_tu=\Delta u+u^p\).

\begin{corollary}\label{cor:belowS}
Let \(N\ge3\), \(m>1\), \(\sigma\in(-2,0)\), \(p>\pc\) and \(\vartheta\in(0,1)\).
Let \(u_0\in\mathcal A_r\), \(r>r_c\), satisfy \(u_0\not\equiv0\) and
\(u_0\le\vartheta S\) a.e.\ in \(\RN\), and set
\(\lambda\eqdef\bigl(\norm{u_0}_\infty/(\vartheta A)\bigr)^{1/\ac}\). Then
\(\Tmax(u_0)=\infty\) and, with \(\ol u_\vartheta\) as in
Theorem~\ref{thm:conj91},
\begin{equation}\label{eq:belowS}
u(t,x)\le\lambda^{\ac}\,\ol u_\vartheta(\lambda x)
\le\Bigl(\frac{\norm{u_0}_\infty}{\vartheta A}\Bigr)\norm{\ol u_\vartheta}_\infty
\qquad\text{for all }t\ge0\text{ and a.e. }x\in\RN .
\end{equation}
Under the assumptions of Proposition~\ref{rem:LequalA}, the same conclusion
holds with \(\vartheta=1\), the barrier \(\ol u_\vartheta\) being replaced by
\(\ol u_1\).
\end{corollary}

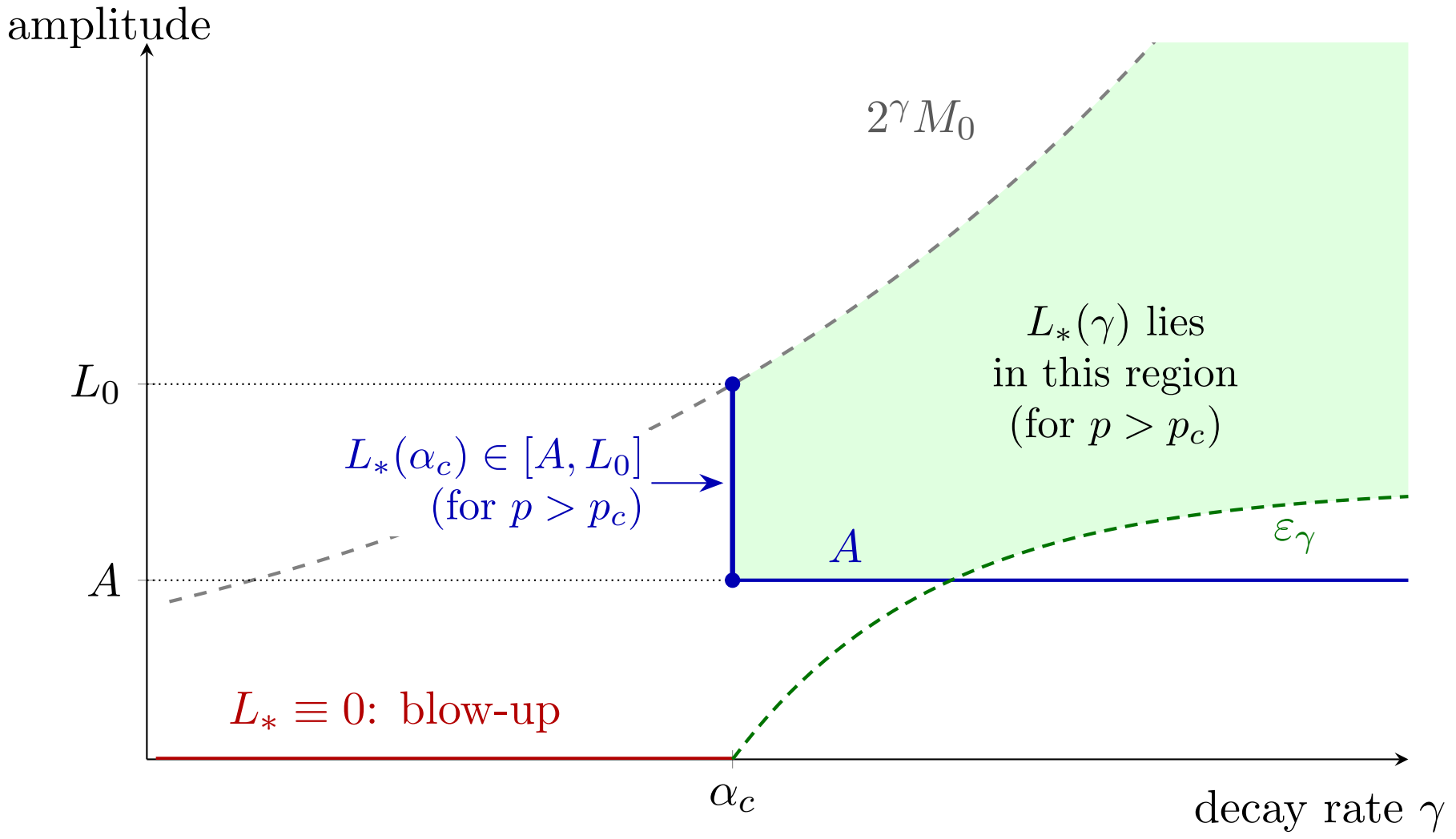
\begin{figure}[ht]
\centering
\begin{tikzpicture}
\begin{axis}[
  width=12cm, height=7.5cm,
  axis lines=left,
  xmin=0, xmax=2.8, ymin=0, ymax=4.2,
  xtick={1.3}, xticklabels={$\alpha_c$},
  ytick={1.05,2.2}, yticklabels={$A$,$L_0$},
  xlabel={decay rate $\gamma$}, ylabel={amplitude},
  x label style={at={(axis description cs:0.93,-0.02)},anchor=north},
  y label style={at={(axis description cs:-0.03,0.97)},anchor=south,rotate=-90}]
  \addplot[name path=up,draw=none,domain=1.3:2.8,samples=60] {0.893*2^(x)};
  \addplot[name path=lo,draw=none,domain=1.3:2.8,samples=80] {max(1.05,1.6*(1-exp(-2.2*(x-1.3))))};
  \addplot[fill=green!20,opacity=0.6] fill between[of=up and lo];
  \addplot[domain=0.05:2.8,samples=80,dashed,thick,gray] {0.893*2^(x)};
  \node[gray!70!black,anchor=east] at (axis cs:1.88,3.75) {$2^{\gamma}M_0$};
  \addplot[very thick,red!70!black] coordinates {(0.02,0) (1.3,0)};
  \node[red!70!black,anchor=south] at (axis cs:0.55,0.06) {$L_\ast\equiv0$: blow-up};
  \draw[densely dotted] (axis cs:0,1.05) -- (axis cs:1.3,1.05);
  \draw[densely dotted] (axis cs:0,2.2) -- (axis cs:1.3,2.2);
  \addplot[domain=1.3:2.8,samples=2,thick,blue!70!black] {1.05};
  \node[blue!70!black,anchor=south] at (axis cs:1.55,1.06) {$A$};
  \draw[very thick,blue!70!black] (axis cs:1.3,1.05) -- (axis cs:1.3,2.2);
  \fill[blue!70!black] (axis cs:1.3,1.05) circle (1.8pt);
  \fill[blue!70!black] (axis cs:1.3,2.2) circle (1.8pt);
  \node[blue!70!black,anchor=east,align=right,font=\small,fill=white,inner sep=1.5pt] (lab) at (axis cs:1.12,1.62) {$L_\ast(\alpha_c)\in[A,L_0]$\\ (for $p>p_c$)};
  \draw[-{Stealth[length=2mm]},blue!70!black] (lab.east) -- (axis cs:1.28,1.62);
  \addplot[domain=1.3:2.8,samples=60,thick,densely dashed,green!45!black] {1.6*(1-exp(-2.2*(x-1.3)))};
  \node[green!45!black,anchor=north] at (axis cs:2.55,1.5) {$\eps_\gamma$};
  \node[font=\small,align=center] at (axis cs:2.15,2.25) {$L_\ast(\gamma)$ lies\\ in this region\\ (for $p>p_c$)};
\end{axis}
\end{tikzpicture}
\caption{Schematic picture of the threshold function $\gamma\mapsto
L_\ast(\gamma)$ of Theorem~\ref{thm:threshold}. It vanishes on $(0,\alpha_c)$,
it is nondecreasing and bounded above by $2^{\gamma}M_0$, and it is bounded
below by $\eps_\gamma>0$ (dashed green) for $\gamma>\alpha_c$. For $p>p_c$, its
value at $\alpha_c$ lies in $[A,L_0]$, where $L_0=2^{\alpha_c}M_0$, and
$L_\ast(\gamma)\ge\max\{A,\eps_\gamma\}$ for $\gamma>\alpha_c$
(Corollary~\ref{cor:sandwich}), which gives the shaded region. For
$p_F<p\le p_c$, one has $0<L_\ast(\alpha_c)\le L_0$
(Theorem~\ref{thm:allp}), together with the bounds by $\eps_\gamma$ and
$2^\gamma M_0$. Only the qualitative features of the curves are meaningful.}
\label{fig:threshold}
\end{figure}

The main results are summarised in Table~\ref{tab:summary}, and the
threshold function is sketched in Figure~\ref{fig:threshold}.

\begin{table}[ht]
\centering
\small
\renewcommand{\arraystretch}{1.3}
\begin{tabular}{>{\raggedright\arraybackslash}p{2.6cm}>{\raggedright\arraybackslash}p{1.6cm}>{\raggedright\arraybackslash}p{6.4cm}>{\raggedright\arraybackslash}p{3.6cm}}
\toprule
Decay of $u_0$ & Range of $p$ & Outcome & Result\\
\midrule
$\gamma<\alpha_c$ & $p>p_F$ & blow-up for every amplitude; $L_\ast(\gamma)=0$ & Thm.~\ref{thm:slow}, \ref{thm:sce}(i), \ref{thm:threshold}(iii)\\
$\gamma=\alpha_c$, $L>L_0$ & $p>p_F$ & blow-up & Thm.~\ref{thm:blowup-crit}\\
$\gamma=\alpha_c$, $L<A$ & $p>p_c$ & global and bounded, below $\ol u_\vartheta$ & Thm.~\ref{thm:conj91}\\
$\gamma=\alpha_c$, $L=A$ & $p>p_c$, \eqref{eq:JL} fails & global and bounded & Prop.~\ref{rem:LequalA}\\
$\gamma=\alpha_c$ & $p>p_c$ & $A\le L_\ast(\alpha_c)\le L_0$ & Cor.~\ref{cor:sandwich}\\
$\gamma=\alpha_c$ & $p>p_F$ & $0<\sup_b c_\infty(b)\le L_\ast(\alpha_c)\le L_0$ & Thm.~\ref{thm:allp}\\
$\gamma>\alpha_c$, small in $\norm{\cdot}_{[\gamma]}$ & $p>p_F$ & global, $L_\ast(\gamma)\ge\eps_\gamma$ & Thm.~\ref{thm:sce}(ii), \ref{thm:threshold}(iv)\\
$\gamma\ge\alpha_c$, $L<A$ & $p>p_c$ & global and bounded, $L_\ast(\gamma)\ge A$ & Cor.~\ref{cor:sandwich}\\
$\gamma>0$, small $L$ and $\norm{u_0}_\infty$ & $p>p_F$ & blow-up may occur & Prop.~\ref{prop:nodecayonly}\\
$u_0\le\vartheta S$, $\vartheta\le1$ & $p>p_c$ & $u(t)\le\vartheta S$ as long as $u$ exists & Thm.~\ref{thm:trap}\\
$u_0\le\vartheta S$, $\vartheta<1$ & $p>p_c$ & global and uniformly bounded & Cor.~\ref{cor:belowS}\\
\bottomrule
\end{tabular}
\caption{Summary of the main results for $N\ge3$, $m>1$ and $\sigma\in(-2,0)$.}
\label{tab:summary}
\end{table}

The paper is organized as follows. Section~\ref{sec:prelim} collects the
functional framework, the comparison principles with barriers which are
unbounded at the origin or not integrable at infinity, the explicit and the
regular radial stationary profiles, and the energy computation on annuli.
Section~\ref{sec:proofs} is devoted to the singular steady state and to global
existence at the critical decay. It contains the proofs of
Proposition~\ref{prop:steady}, Theorems~\ref{thm:trap}, \ref{thm:barrier},
\ref{thm:global} and \ref{thm:conj91}, Corollaries~\ref{cor:sandwich} and
\ref{cor:belowS}, and Proposition~\ref{rem:LequalA}.
Section~\ref{sec:proofs-blowup} contains the proofs of the blow-up results and
of the characterization of the second critical exponent, namely
Theorems~\ref{thm:slow}, \ref{thm:blowup-crit}, \ref{thm:sce} and
\ref{thm:threshold}, Proposition~\ref{prop:nodecayonly} and
Corollary~\ref{cor:char}. Section~\ref{sec:conjectures} discusses the
identification of \(L_\ast(\ac)\), slowly varying corrections to the critical
tail and the range \(\pF<p\le\pc\). Finally, Section~\ref{sec:discussion}
places the results in a broader perspective.

%======================================================
\section{{Preliminaries and technical tools}}\label{sec:prelim}
%======================================================

In this section we introduce the notation and the functional framework used
throughout the paper, and we recall the results of \cite{IL1} on which our
analysis rests. We then collect the technical tools needed for the comparison
arguments: stationary weak supersolutions, comparison lemmas for barriers which
are unbounded at the origin or not integrable at infinity, two families of
stationary profiles, and an energy computation on annuli.

\subsection{{Notation}}
For \(r\ge0\), we write \(\Lpp{r}\) for the nonnegative cone of \(\Lp{r}\) and
\(\norm{\cdot}_r\) for the norm of \(\Lp{r}\). We denote by \(\omega_N\) the
volume of the unit ball of \(\RN\), by \(B(0,R)\) the open ball of radius \(R\)
centred at the origin, by \(\Hmeas\) the \((N-1)\)-dimensional Hausdorff
measure, and by \(\Dprime{\Omega}\) the space of distributions on an open set
\(\Omega\subseteq\RN\). For \(\gamma>0\) we set
\begin{equation}\label{eq:phigamma}
\phi_\gamma(x)\eqdef\min\{1,\abs{x}^{-\gamma}\},\qquad
\norm{u_0}_{[\gamma]}\eqdef\sup_{x\in\RN}\max\{1,\abs{x}^{\gamma}\}\,u_0(x),
\end{equation}
so that \(0\le u_0\le \eps\phi_\gamma\) if and only if \(u_0\ge0\) and
\(\norm{u_0}_{[\gamma]}\le\eps\). We abbreviate \(\phi\eqdef\phi_{\ac}\). Test
functions are denoted by \(\varphi\). Constants are allowed to depend on
\(N,m,p,\sigma\) without further mention.

\subsection{{Weak solutions and the well-posedness framework}}
We use the notion of weak solution and the well-posedness results of
\cite{IL1}.

\begin{definition}[Weak solution, {\cite[Definition~2.1]{IL1}}]\label{def:weak}
Let \(r\ge0\), \(T\in(0,\infty]\) and
\(u_0\in \Lpp{r+1}\cap \Lp{\infty}\). A nonnegative function \(u\) belonging to
\(L^{\infty}_{\mathrm{loc}}\bigl([0,T),\Lpp{r+1}\cap\Lp{\infty}\bigr)\), and
such that \(\nabla u^{(m+r)/2}\in L^2_{\mathrm{loc}}([0,T),\Lp2)\), is a weak
solution of \eqref{eq:CP} on \([0,T)\) if
\[
\int_0^T\!\!\int_{\RN}\Bigl[(u_0-u)\partial_t\varphi-u^m\Delta\varphi
-\abs{x}^{\sigma}u^{p}\varphi\Bigr]\dx\ds=0
\qquad\text{for all }\varphi\in C^2_c([0,T)\times\RN).
\]
For \(u_0\in\Lpp{1}\cap\Lp{\infty}\), the solution is an \emph{energy
solution} in the sense of \cite[Corollary~2.3]{IL1}, namely
\(\nabla u^m\in L^2_{\mathrm{loc}}([0,T),\Lp2)\) and the first-order weak
formulation holds.
\end{definition}

Data with the critical or with a slower algebraic tail are \emph{not}
integrable. It is therefore essential to record the exact class in which the
theory of \cite{IL1} is available. We shall systematically work in
\begin{equation}\label{eq:class}
\mathcal A_r\eqdef \Lpp{r+1}\cap\Lp{\infty},\qquad r>r_c=\max\{r_0,0\}.
\end{equation}
By interpolation with \(\Lp\infty\), we have \(\mathcal A_{r}\subset\mathcal
A_{r'}\) whenever \(r\le r'\). Hence the family \eqref{eq:class} is increasing
in \(r\), and it contains every datum considered in this paper. Indeed, a datum
\(u_0\in\Lp\infty\) with \(u_0(x)\lesssim\abs{x}^{-\gamma}\) at infinity
belongs to \(\mathcal A_r\) for every \(r>\max\{r_c,N/\gamma-1\}\). The
condition \(r>r_c\) is exactly the condition \(r_1>r_c\) of the existence
result \cite[Theorem~2.2]{IL1}.

\begin{assumption}[The framework of \cite{IL1}]\label{ass:IL}
In the whole paper we assume $m>1$, $\sigma\in(-2,0)$, $N\ge3$ and
$p>p_F$, so that $p>p_F>m>p_G$. Let \(u_0\in\mathcal A_r\) for some
\(r>r_c\). For \(\eta\in(0,1)\), let \(u_\eta\) denote the unique global
weak solution of
\begin{equation}\label{eq:Reta}
\partial_t u_\eta=\Delta u_\eta^m
+\bigl(\abs{x}^2+\eta^2\bigr)^{\sigma/2}\frac{u_\eta^{p}}{1+\eta u_\eta^{p-1}},
\qquad u_\eta(0)=u_0\,\mathbf 1_{B(0,1/\eta)},
\tag{$\mathrm R_\eta$}
\end{equation}
which belongs to \(C([0,\infty),\Lpp1)\cap
L^\infty_{\mathrm{loc}}([0,\infty),\Lp\infty)\) and satisfies the comparison
principle \cite[\S3.1, Eqs.~(3.22)--(3.23)]{IL1}. We call \emph{the solution}
of \eqref{eq:CP} the monotone limit \eqref{eq:monotone} below, and we set
\begin{equation*}
\Tmax(u_0)\eqdef\sup\Bigl\{T>0:\ \sup_{\eta\in(0,1)}
\norm{u_\eta}_{L^\infty((0,T)\times\RN)}<\infty\Bigr\}\in(0,\infty].
\end{equation*}
We use the following properties.
\begin{enumerate}[label=\textup{(\Alph*)},leftmargin=2.4em]
\item \textup{(Well-posedness and blow-up criterion.)}
The solution \(u\) is a nonnegative weak solution of \eqref{eq:CP} on
\([0,\Tmax(u_0))\) in the sense of Definition~\ref{def:weak}. If
\(\Tmax(u_0)<\infty\), then
\[
\limsup_{t\uparrow\Tmax(u_0)}\norm{u(t)}_{\infty}=\infty .
\]
If moreover \(u_0\in\Lpp1\cap\Lp\infty\), then \(u\) is the unique solution of
\eqref{eq:CP} given by \cite[Corollary~2.7]{IL1}, and \(\Tmax(u_0)\) is its
maximal existence time.
\item \textup{(Comparison.)} If \(u_{0,1}\le u_{0,2}\) with
\(u_{0,i}\in\mathcal A_{r}\), then the corresponding solutions satisfy
\(u_1\le u_2\) on the common existence interval. In particular,
\(\Tmax(u_{0,2})\le \Tmax(u_{0,1})\).
\item \textup{(Regularizing approximations, \cite[\S3.1, Eq.~(3.22)]{IL1}.)}
The map \(\eta\mapsto u_\eta\) is nonincreasing, and the solution of
\eqref{eq:CP} is the monotone limit
\begin{equation}\label{eq:monotone}
u(t,x)=\lim_{\eta\to0}u_\eta(t,x)=\sup_{\eta\in(0,1)}u_\eta(t,x).
\end{equation}
\end{enumerate}
\end{assumption}

Let us justify (A) and (B). Uniqueness and comparison are established in
\cite[Theorem~2.5]{IL1} only for data in \(\Lpp1\cap\Lp\infty\), whereas the
data considered here are, in general, not integrable. The monotonicity in (C)
is \cite[Eq.~(3.23)]{IL1}, and \eqref{eq:monotone} is the construction
\cite[Eq.~(3.26)]{IL1} in the proof of \cite[Theorem~2.2]{IL1}. We also note
that, for \(u_0\in\Lpp1\cap\Lp\infty\), the above definition of \(\Tmax(u_0)\)
is exactly the one used in the proof of \cite[Corollary~2.7]{IL1}.

For (B), the truncated data are ordered, \(u_{0,1}\mathbf 1_{B(0,1/\eta)}\le
u_{0,2}\mathbf 1_{B(0,1/\eta)}\). The comparison principle for \eqref{eq:Reta}
gives \(u_{\eta,1}\le u_{\eta,2}\), and letting \(\eta\to0\) in
\eqref{eq:monotone}, we obtain \(u_1\le u_2\). Moreover, by monotonicity,
\(\norm{u_i}_{L^\infty((0,T)\times\RN)}=\sup_{\eta}\norm{u_{\eta,i}}_{L^\infty((0,T)\times\RN)}\)
for every \(T>0\), and the inequality \(\Tmax(u_{0,2})\le\Tmax(u_{0,1})\)
follows from the definition of \(\Tmax\).

For (A), fix \(T<\Tmax(u_0)\) and let \(M\) be a bound for \(u_\eta\) on
\((0,T)\times\RN\), uniform in \(\eta\). Let \(s\ge0\) be such that
\(u_0\in\Lp{s+1}\), which is the case for every \(s\ge r\). We multiply
\eqref{eq:Reta} by \(u_\eta^{s}\) as in \cite[Eq.~(3.9)]{IL1}, and we use the
bound \(g_\eta(x,\tau)\le\abs{x}^{\sigma}\tau^{p}\) of
Lemma~\ref{lem:lipschitz} below. Splitting \(\RN\) into \(B(0,1)\), on which
\(\abs{x}^{\sigma}\) is integrable, and its complement, on which
\(\abs{x}^{\sigma}\le1\), we get
\[
\frac{d}{dt}\norm{u_\eta(t)}_{s+1}^{s+1}
+\frac{4ms(s+1)}{(m+s)^2}\norm{\nabla u_\eta^{(m+s)/2}(t)}_2^2
\le(s+1)\Bigl(\frac{N\omega_N}{N+\sigma}M^{p+s}
+M^{p-1}\norm{u_\eta(t)}_{s+1}^{s+1}\Bigr).
\]
Since \(\norm{u_\eta(0)}_{s+1}\le\norm{u_0}_{s+1}\), Gronwall's lemma bounds
\(\norm{u_\eta(t)}_{s+1}\) and \(\int_0^T\norm{\nabla u_\eta^{(m+s)/2}}_2^2\)
uniformly in \(\eta\) and \(t\in[0,T]\). Since \(0\le u_\eta\uparrow u\le M\),
dominated convergence allows us to pass to the limit in the weak formulation of
\eqref{eq:Reta}. Here the reaction terms converge in \(L^1_{\mathrm{loc}}\),
because they converge pointwise and are dominated by
\(M^p\abs{x}^{\sigma}\in L^1_{\mathrm{loc}}(\RN)\). Fatou's lemma and the weak
compactness of \(\nabla u_\eta^{(m+r)/2}\) in \(L^2\) give the regularity
required in Definition~\ref{def:weak}. Hence \(u\) is a weak solution on
\([0,T)\).

Assume next that \(\Tmax(u_0)<\infty\) and that \(\norm{u(t)}_\infty\le M\)
for \(t<\Tmax(u_0)\). Then, for every \(\eta\), we have \(u_\eta\le M\) on
\((0,\Tmax(u_0))\times\RN\), and the above inequality with \(s=r\) gives
\(K\eqdef\sup_{\eta}\sup_{t<\Tmax(u_0)}\norm{u_\eta(t)}_{r+1}<\infty\). We now
use the a priori estimates of \cite{IL1} for \eqref{eq:Reta}. As stated in the
proof of \cite[Theorem~2.2]{IL1}, see \cite[Eq.~(3.25)]{IL1}, the arguments of
\cite[Propositions~3.3 and 3.4]{IL1} apply verbatim to the regularized problem.
Indeed, they only rely on the differential identity \cite[Eq.~(3.9)]{IL1},
on the functional inequality \cite[Proposition~3.2]{IL1}, which concerns a
fixed function and not the equation, and on an \(L^\infty\) estimate of
Kawanago \cite[Th\'eor\`eme~1]{Kawanago93}. For \eqref{eq:Reta}, the identity
holds with \(\abs{x}^{\sigma}u^{p}\) replaced by \(g_\eta(x,u_\eta)\). Since
\(0\le g_\eta(x,\tau)\le\abs{x}^{\sigma}\tau^{p}\), it becomes an inequality in
the direction used in \cite{IL1}, and the source bound required by Kawanago's
estimate is unchanged.
Since \eqref{eq:Reta} is autonomous, we may apply these estimates to the
solution \(t\mapsto u_\eta(t_0+t)\), whose initial value satisfies
\(\norm{u_\eta(t_0)}_{r+1}\le K\) and \(\norm{u_\eta(t_0)}_\infty\le M\). In
the proof of \cite[Proposition~3.3]{IL1}, the existence time is
\(\norm{v_0}_{r+1}^{(r+1)(1-\nu_r)}/C\) with \(\nu_r>1\) (because \(p>p_G\)),
and all the constants of \cite[Propositions~3.3 and 3.4]{IL1} depend on the
initial value \(v_0\) only through upper bounds for \(\norm{v_0}_{r+1}\) and
\(\norm{v_0}_\infty\). We infer that there are \(T_1>0\) and \(C>0\), depending
only on \(N,m,p,\sigma,r,K,M\), such that
\[
\norm{u_\eta(t)}_\infty\le C\qquad\text{for all }t\in[t_0,t_0+T_1],\
t_0\in(0,\Tmax(u_0))\text{ and }\eta\in(0,1).
\]
Choosing \(t_0=\Tmax(u_0)-T_1/2\), we deduce that the quantity
\(\sup_\eta\norm{u_\eta}_{L^\infty((0,\Tmax(u_0)+T_1/2)\times\RN)}\) is
finite. This contradicts the definition of \(\Tmax(u_0)\).

Finally, let \(u_0\in\Lpp1\cap\Lp\infty\). Then the above estimates with
\(s=0\) and \(s=m\) show that \(u\) has the regularity of
\cite[Corollary~2.3]{IL1}. By \cite[Lemma~3.6]{IL1}, which relies on the
uniqueness part of \cite[Theorem~2.5]{IL1}, any weak solution in the sense of
\cite[Corollary~2.3]{IL1} coincides with \(\sup_\eta u_\eta\). Together with the
coincidence of the two definitions of \(\Tmax(u_0)\), this identifies \(u\) with
the solution of \cite[Corollary~2.7]{IL1}.

Two elementary properties of the regularized solutions will be used
repeatedly. For fixed \(\eta\in(0,1)\), the solution \(u_\eta\) belongs to
\(C([0,\infty),\Lpp1)\cap L^\infty_{\mathrm{loc}}([0,\infty),\Lp\infty)\)
\cite[\S3.1]{IL1}. Hence, for every \(T>0\),
\begin{equation}\label{eq:ueta-infty}
K_{\eta,T}\eqdef\sup_{t\in[0,T]}\norm{u_\eta(t)}_{\infty}<\infty ,
\end{equation}
and
\begin{equation}\label{eq:ueta-one}
M_{\eta,T}\eqdef\sup_{t\in[0,T]}\norm{u_\eta(t)}_{1}<\infty .
\end{equation}
The constants \(K_{\eta,T}\) and \(M_{\eta,T}\) may depend on \(\eta\), and
this is harmless in the comparison arguments below, where \(\eta\) is fixed
until the very last step. We emphasise that \eqref{eq:ueta-one} concerns the
solution with the truncated datum \(u_0\mathbf 1_{B(0,1/\eta)}\), which is
bounded with compact support. The data considered in this paper need not
belong to \(\Lp1\).

\begin{theorem}[Supercritical existence and blow-up, {\cite[Thms.~2.9, 2.10]{IL1}}]
\label{thm:IL-super}
Assume \(p>\pF\), so that \(r_0=r_c>0\).
\begin{enumerate}[label=\textup{(\alph*)},leftmargin=2.2em]
\item \emph{(Small data.)} There is \(C_0=C_0(N,m,p,\sigma)>0\) such that, if
\(u_0\in\Lpp{r_0+1}\cap\Lp{\infty}\) with \(\norm{u_0}_{r_0+1}\le C_0\), then
the solution of \eqref{eq:CP} is global and
\(u\in L^{\infty}\bigl((0,\infty),\Lp{r}\bigr)\) for all
\(r\in[r_0+1,\infty)\).
\item \emph{(Negative energy.)} Let \(p\ge m\), and let
\(u_0\in\Lpp1\cap\Lp{\infty}\) with
\(u_0\in L^{m+p}_+(\RN,\abs{x}^{\sigma}\dx)\) and \(\nabla u_0^m\in\Lp2\). If
\begin{equation}\label{eq:energy}
E(u_0)\eqdef\frac{\norm{\nabla u_0^m}_2^2}{2}
-\frac{m}{m+p}\int_{\RN}\abs{x}^{\sigma}u_0(x)^{m+p}\dx<0,
\end{equation}
then \(\Tmax(u_0)<\infty\).
\end{enumerate}
\end{theorem}

Let us comment on the way Theorem~\ref{thm:IL-super} is used. In part (a),
\cite[Theorem~2.10]{IL1} concerns the solution constructed in
\cite[Theorem~2.2]{IL1}, that is, the monotone limit \eqref{eq:monotone}. Its
proof in \cite[\S5]{IL1} consists of a priori estimates: a uniform
\(L^{r_0+1}\) bound, uniform \(L^{r+1}\) bounds for all \(r\ge r_0\), and an
\(L^\infty\) bound on every bounded time interval. They are derived from the
differential inequalities of \cite[\S5]{IL1} and from the \(L^\infty\)
estimate of \cite[Proposition~3.4]{IL1}. As in the proof of
\cite[Theorem~2.2]{IL1}, these computations apply verbatim to the
solutions \(u_\eta\) of \eqref{eq:Reta}. This is because
\(g_\eta(x,\tau)\le\abs{x}^{\sigma}\tau^{p}\) and
\(\norm{u_0\mathbf 1_{B(0,1/\eta)}}_q\le\norm{u_0}_q\) for every \(q\), so that
all the resulting bounds are uniform in \(\eta\). Consequently,
\(\sup_\eta\norm{u_\eta}_{L^\infty((0,T)\times\RN)}<\infty\) for every \(T>0\),
that is, \(\Tmax(u_0)=\infty\) in the sense of Assumption~\ref{ass:IL}, and
Fatou's lemma gives \(u\in L^{\infty}((0,\infty),\Lp{r})\). We also note that
the proof of \cite[Theorem~2.10]{IL1} uses the smallness condition in the form
\(\norm{u_0}_{r_0+1}^{r_0+1}<C_0^{\mathrm{IL}}\), see \cite[Eq.~(5.3)]{IL1}.
Part (a) therefore holds with
\(C_0\eqdef(C_0^{\mathrm{IL}}/2)^{1/(r_0+1)}\), which is the constant used
below. In part (b), the datum is integrable, so that the solution is the unique
one of Assumption~\ref{ass:IL}(A).

The energy \eqref{eq:energy} is the natural Lyapunov functional of
\eqref{eq:CP}. Multiplying \eqref{eq:CP} by \(\partial_t u^m\) and
integrating gives, formally,
\(m\int u^{m-1}(\partial_t u)^2\dx=-\tfrac{d}{dt}E(u(t))\), so that \(E\) is
nonincreasing along the flow. The proof of \cite[Theorem~2.9]{IL1} then
follows the argument of Berryman and Holland \cite{BerrymanHolland80}, and an
alternative route is the concavity method of Levine \cite{Levine73}. Both
convert \(E(u_0)<0\) into finite-time blow-up. We refer to
\cite{Levine90,DengLevine00} for surveys on the role of critical exponents in
blow-up theorems, and to
\cite{SGKM95,Galaktionov94,GalaktionovVazquez02,Vazquez07} for the quasilinear
background.

\subsection{{Two elliptic identities across the origin}}
Two elementary elliptic identities will be used repeatedly. For
\(b\in\R\setminus\{0\}\) we have
\begin{equation}\label{eq:laplacian}
\Delta\bigl(\abs{x}^{-b}\bigr)=b\,(b+2-N)\,\abs{x}^{-b-2},
\qquad
\Delta\bigl(\abs{x}^{\,b}\bigr)=b\,(b+N-2)\,\abs{x}^{\,b-2},
\qquad x\in\RN\setminus\{0\}.
\end{equation}
The following lemma records when these pointwise identities persist across
the origin in the sense of distributions. This is the point at which the
threshold \(\pc\) enters the comparison theory. In this singular setting, the
absence of a Dirac mass at the origin has to be checked whenever \(S\) is
compared with a solution.
\begin{lemma}\label{lem:dirac}
Let \(N\ge3\).
\begin{enumerate}[label=\textup{(\roman*)},leftmargin=2.2em]
\item If \(0<b<N-2\), then \(\abs{x}^{-b}\in L^1_{\mathrm{loc}}(\RN)\),
\(\abs{x}^{-b-2}\in L^1_{\mathrm{loc}}(\RN)\), and
\(\Delta\abs{x}^{-b}=b(b+2-N)\abs{x}^{-b-2}\) in \(\Dprime{\RN}\).
\item If \(b>0\), then \(\abs{x}^{b-2}\in L^1_{\mathrm{loc}}(\RN)\) and
\(\Delta\abs{x}^{b}=b(b+N-2)\abs{x}^{b-2}\) in \(\Dprime{\RN}\).
\end{enumerate}
\end{lemma}

\begin{proof}
Let \(\varphi\in C^\infty_c(\RN)\), and let \(f\) be smooth on
\(\RN\setminus\{0\}\). Green's formula on \(\{\abs{x}>\eps\}\) reads
\begin{equation}\label{eq:green}
\int_{\{\abs{x}>\eps\}}f\,\Delta\varphi\dx
=\int_{\{\abs{x}>\eps\}}\Delta f\,\varphi\dx
+\int_{\{\abs{x}=\eps\}}\Bigl(f\,\partial_\nu\varphi-\varphi\,\partial_\nu f\Bigr)
\,d\Hmeas,
\end{equation}
where \(\nu\) denotes the \emph{outer} normal to the domain
\(\{\abs{x}>\eps\}\), that is, \(\nu=-x/\abs{x}\) on \(\{\abs{x}=\eps\}\). No
boundary term appears at infinity, since \(\varphi\) has compact support. We
let \(\eps\downarrow0\) in \eqref{eq:green}, using \eqref{eq:laplacian}.

\smallskip
\noindent\emph{(i) The case \(f(x)=\abs{x}^{-b}\), \(0<b<N-2\).}
By \eqref{eq:laplacian}, \(\Delta f=b(b+2-N)\abs{x}^{-b-2}\) on
\(\RN\setminus\{0\}\). On the sphere \(\{\abs{x}=\eps\}\) we have
\(f=\eps^{-b}\), \(\partial_\nu f=b\,\eps^{-b-1}\), and
\(\Hmeas(\{\abs{x}=\eps\})=N\omega_N\eps^{N-1}\). Hence the two boundary
contributions in \eqref{eq:green} satisfy
\[
\int_{\{\abs{x}=\eps\}}f\,\partial_\nu\varphi\,d\Hmeas
=O\bigl(\eps^{N-1-b}\bigr),
\qquad
-\int_{\{\abs{x}=\eps\}}\varphi\,\partial_\nu f\,d\Hmeas
=O\bigl(\eps^{N-2-b}\bigr).
\]
The first term vanishes as \(\eps\downarrow0\) because \(b<N-1\). The second
one, which is the binding term, vanishes exactly when \(b<N-2\). If \(b=N-2\),
it converges to \(-(N-2)N\omega_N\varphi(0)\neq0\) in general, which produces a
Dirac mass at the origin, and if \(b>N-2\), it even blows up. The condition
\(b<N-2\) also guarantees that \(\abs{x}^{-b-2}\in L^1_{\mathrm{loc}}(\RN)\),
while \(b<N\) gives \(\abs{x}^{-b}\in L^1_{\mathrm{loc}}(\RN)\). Since
\(\varphi\) has compact support, both \(f\Delta\varphi\) and
\(\Delta f\,\varphi\) belong to \(L^1(\RN)\), and dominated convergence allows
us to pass to the limit in the volume integrals of \eqref{eq:green}. This
proves (i).

\smallskip
\noindent\emph{(ii) The case \(f(x)=\abs{x}^{b}\), \(b>0\).}
Here \(\Delta f=b(b+N-2)\abs{x}^{b-2}\) on \(\RN\setminus\{0\}\). On
\(\{\abs{x}=\eps\}\) we have \(f=\eps^{b}\) and \(\partial_\nu f=-b\,\eps^{b-1}\),
so that the two boundary contributions are of order
\(\eps^{N-1+b}\) and \(\eps^{N-2+b}\). Both vanish as \(\eps\downarrow0\),
because \(N\ge3\) and \(b>0\). Moreover \(\abs{x}^{b-2}\in L^1_{\mathrm{loc}}(\RN)\),
since \(b-2>-2>-N\). Hence no restriction on \(b\) beyond \(b>0\) is needed,
and the same dominated convergence argument as in (i) proves (ii).
\end{proof}
%=====================================================================
\subsection{Stationary weak supersolutions}\label{sec:background}
The barriers used in the comparison arguments fall outside the class of
comparison functions for which \cite[Theorem~2.5]{IL1} is stated. The singular
steady state $S$ is unbounded at the origin, and the bounded profile
$\overline{u}$ is not integrable at infinity. We therefore isolate, once and
for all, a comparison lemma adapted to them. The point is that it suffices to
compare with the \emph{regularized} problems \eqref{eq:Reta}, whose reaction
term is globally Lipschitz, and then to pass to the monotone limit
\eqref{eq:monotone}.
\begin{definition}[Stationary weak supersolution]\label{def:supersol}
Let $\Omega\subseteq\mathbb{R}^N$ be open. A function
$v\in L^1_{\mathrm{loc}}(\Omega)$ with $v>0$ a.e.,
$v^m\in H^1_{\mathrm{loc}}(\Omega)$ and
$|x|^{\sigma}v^{p}\in L^1_{\mathrm{loc}}(\Omega)$ is a \emph{stationary weak
supersolution} of \eqref{eq:CP} in $\Omega$ if
\begin{equation}\label{eq:supersol}
\int_{\Omega}\Bigl(\nabla v^m\cdot\nabla\varphi
-|x|^{\sigma}v^{p}\,\varphi\Bigr)\,dx\ \ge\ 0
\qquad\text{for all }0\le\varphi\in C^\infty_c(\Omega),
\end{equation}
that is, $\Delta v^m+|x|^{\sigma}v^{p}\le0$ in $\mathcal{D}'(\Omega)$.
\end{definition}
Let $\varphi\in H^1(\Omega)\cap L^\infty(\Omega)$ be nonnegative with compact
support in $\Omega$. Mollifying $\varphi$ produces nonnegative test functions
$\varphi_k\in C^\infty_c(\Omega)$, supported in a fixed compact subset of
$\Omega$ and uniformly bounded, such that $\varphi_k\to\varphi$ in $H^1$ and
a.e. Since $\nabla v^m\in L^2$ and $|x|^{\sigma}v^{p}\in L^1$ on that compact
set, we may pass to the limit in \eqref{eq:supersol}. Hence
\eqref{eq:supersol} extends to every such $\varphi$. This extension is used in
the proof of Lemma~\ref{lem:comparison}.

A first, elementary consequence of Definition~\ref{def:supersol} and of the
superlinearity \(p>m\) is that sub-multiples of a stationary supersolution
remain stationary supersolutions. This is used in the trapping argument of
Theorem~\ref{thm:trap} and in the construction of the bounded barrier of
Theorem~\ref{thm:barrier}.
\begin{lemma}\label{lem:submultiple}
Let $p>m$, let $\Omega\subseteq\mathbb{R}^N$ be open and let $v$ be a
stationary weak supersolution of \eqref{eq:CP} in $\Omega$. Then $\vartheta v$
is a stationary weak supersolution of \eqref{eq:CP} in $\Omega$ for every
$\vartheta\in(0,1]$.
\end{lemma}
\begin{proof}
The regularity requirements of Definition~\ref{def:supersol} are invariant
under multiplication by a positive constant. Indeed, since \(\vartheta>0\), we
have \(\vartheta v\in L^1_{\mathrm{loc}}(\Omega)\), \(\vartheta v>0\) a.e.,
\((\vartheta v)^m=\vartheta^m v^m\in H^1_{\mathrm{loc}}(\Omega)\) and
\(\abs{x}^{\sigma}(\vartheta v)^{p}=\vartheta^{p}\abs{x}^{\sigma}v^{p}
\in L^1_{\mathrm{loc}}(\Omega)\).

It remains to verify the distributional inequality \eqref{eq:supersol}. For
\(0\le\varphi\in C^\infty_c(\Omega)\), we have
\begin{align*}
\int_{\Omega}\Bigl(\nabla(\vartheta v)^m\cdot\nabla\varphi
-\abs{x}^{\sigma}(\vartheta v)^{p}\varphi\Bigr)\dx
&=\vartheta^{m}\int_{\Omega}\Bigl(\nabla v^m\cdot\nabla\varphi
-\abs{x}^{\sigma}v^{p}\varphi\Bigr)\dx\\
&\quad+\bigl(\vartheta^{m}-\vartheta^{p}\bigr)
\int_{\Omega}\abs{x}^{\sigma}v^{p}\varphi\dx
=:J_{1}+J_{2}.
\end{align*}
The term \(J_{1}\) is nonnegative, since \(\vartheta^{m}>0\) and the integral is
nonnegative by \eqref{eq:supersol}. The term \(J_{2}\) is also nonnegative,
since \(0<\vartheta\le1\) and \(p>m\) yield \(\vartheta^{m}-\vartheta^{p}\ge0\),
while \(\abs{x}^{\sigma}v^{p}\varphi\ge0\) a.e.\ on \(\Omega\). Hence
\(J_1+J_2\ge0\) for all \(0\le\varphi\in C^\infty_c(\Omega)\), which is
precisely \eqref{eq:supersol} for \(\vartheta v\).
\end{proof}
\subsection{A global Lipschitz truncation}

The following lemma shows that, in the regularized problem \eqref{eq:Reta},
the reaction term is a global Lipschitz truncation of the original reaction
\(\abs{x}^{\sigma}u^{p}\). This property is the key ingredient which allows us
to run the \(L^{1}\)-contraction argument in the proof of the comparison lemma
below, since it makes the Gronwall step independent of any local boundedness
assumption on the barrier.

\begin{lemma}
\label{lem:lipschitz}
For $\eta\in(0,1)$, set
$g_\eta(x,s)\eqdef\bigl(|x|^2+\eta^2\bigr)^{\sigma/2}\,
\dfrac{s^{p}}{1+\eta s^{p-1}}$ for $x\in\mathbb{R}^N$ and $s\ge0$. Then, for
all $x\in\mathbb{R}^N$ and $s,\hat{s}\ge0$,
\begin{equation}\label{eq:lip}
0\le g_\eta(x,s)\le |x|^{\sigma}s^{p},
\qquad
0\le g_\eta(x,s)\le \eta^{\sigma-1}s,
\qquad
|g_\eta(x,s)-g_\eta(x,\hat{s})|\le \Lambda_\eta\,|s-\hat{s}|,
\end{equation}
with $\Lambda_\eta\eqdef p\,\eta^{\sigma-1}$.
\end{lemma}

\begin{proof}
Since \(\sigma<0\), the map \(t\mapsto t^{\sigma/2}\) is nonincreasing on
\((0,\infty)\). The inequalities \(\eta^{2}\le\abs{x}^{2}+\eta^{2}\) and
\(\abs{x}^{2}\le\abs{x}^{2}+\eta^{2}\) then yield
\[
\bigl(\abs{x}^{2}+\eta^{2}\bigr)^{\sigma/2}
\le\min\{\eta^{\sigma},\abs{x}^{\sigma}\}.
\]
The first estimate in \eqref{eq:lip} follows from
\(\bigl(\abs{x}^{2}+\eta^{2}\bigr)^{\sigma/2}\le\abs{x}^{\sigma}\) and
\((1+\eta s^{p-1})^{-1}\le1\). The second one follows from
\(\bigl(\abs{x}^{2}+\eta^{2}\bigr)^{\sigma/2}\le\eta^{\sigma}\) and
\(s^{p}/(1+\eta s^{p-1})\le s^{p}/(\eta s^{p-1})=s/\eta\).

For the third estimate, set \(h(s)\eqdef s^{p}/(1+\eta s^{p-1})\) for \(s\ge0\).
For \(s>0\), we compute
\[
h'(s)=\frac{p\,s^{p-1}+\eta\,s^{2p-2}}{\bigl(1+\eta s^{p-1}\bigr)^{2}}
\le\frac{p\,s^{p-1}\bigl(1+\eta s^{p-1}\bigr)}{\bigl(1+\eta s^{p-1}\bigr)^{2}}
=\frac{p\,s^{p-1}}{1+\eta s^{p-1}}
\le\frac{p}{\eta},
\]
where we used \(p\ge1\). Hence \(h\) is globally Lipschitz on \([0,\infty)\),
with Lipschitz constant \(p/\eta\). Since
\(g_\eta(x,s)=\bigl(\abs{x}^{2}+\eta^{2}\bigr)^{\sigma/2}h(s)\), we deduce that,
for all \(x\in\mathbb{R}^N\) and all \(s,\hat{s}\ge0\),
\[
\abs{g_\eta(x,s)-g_\eta(x,\hat{s})}
\le\eta^{\sigma}\cdot\frac{p}{\eta}\abs{s-\hat{s}}
=p\,\eta^{\sigma-1}\abs{s-\hat{s}},
\]
which is the third estimate in \eqref{eq:lip}.
\end{proof}

\begin{remark}\rm\label{rem:lipconst}
The constant \(p\,\eta^{\sigma-1}\) is not sharp, but it is convenient. Writing
\(y=\eta s^{p-1}\), we have \(h'(s)=\eta^{-1}(py+y^2)(1+y)^{-2}\), and an
elementary computation gives
\[
\sup_{s>0}h'(s)=\frac1\eta\quad\text{if }1<p\le2,\qquad
\sup_{s>0}h'(s)=\frac{p^{2}}{4(p-1)\,\eta}\quad\text{if }p>2.
\]
Indeed, the derivative of \(y\mapsto(py+y^2)(1+y)^{-2}\) has the sign of
\(p+(2-p)y\). For \(p\le2\) this function increases to its limit \(1\), while
for \(p>2\) it attains its maximum at \(y=p/(p-2)\).
Hence the Lipschitz constant \(\eta^{\sigma-1}\) stated after
\cite[Eq.~(3.22)]{IL1} is valid for \(1<p\le2\), while for \(p>2\) it must be
multiplied by \(p^2/(4(p-1))>1\). This does not affect the arguments of
\cite{IL1}, which only use that the regularized reaction is globally
Lipschitz.
\end{remark}
%================================================
\subsection{Comparison with a stationary supersolution}
We now prove the comparison lemma which is used throughout the sequel. It
allows us to compare the weak solution of \eqref{eq:CP} with stationary
supersolutions which are unbounded at the origin and not integrable at
infinity. This situation is not covered by Assumption~\ref{ass:IL}(B), which
only compares two solutions with data in the class \eqref{eq:class}.
\begin{lemma}\label{lem:comparison}
Let $u_0\in\mathcal{A}_r$ for some $r>r_c$, and let $u$ be the solution of
\eqref{eq:CP} on $[0,T_{\max}(u_0))$. Let
$v:\mathbb{R}^N\setminus\{0\}\to(0,\infty)$ satisfy
\begin{enumerate}[label=\textup{(H\arabic*)},leftmargin=2.6em]
\item $v\in L^1_{\mathrm{loc}}(\mathbb{R}^N)$,
$v\in L^\infty(\{|x|\ge\delta\})$ for every $\delta>0$, and $v$ is a
stationary weak supersolution of \eqref{eq:CP} in
$\mathbb{R}^N\setminus\{0\}$ in the sense of Definition~\ref{def:supersol};
\item either $\displaystyle\lim_{x\to0}v(x)=+\infty$, or
$v\in L^\infty(B(0,1))$, $v^m\in H^1(B(0,1))$ and $v$ is a stationary weak
supersolution in all of $\mathbb{R}^N$;
\item $u_0\le v$ a.e.\ in $\mathbb{R}^N$.
\end{enumerate}
Then
\begin{equation}\label{eq:comparison}
u(t,x)\le v(x)\qquad\text{for a.e. }x\in\mathbb{R}^N\setminus\{0\}\text{ and all }
t\in[0,T_{\max}(u_0)).
\end{equation}
\end{lemma}
Before proving Lemma~\ref{lem:comparison}, we state the elementary Gronwall
lemma which will be used in the proof.

\begin{lemma}\label{lem:gronwall}
Let \(T>0\), let \(\Lambda\ge0\), and let \(F\in L^\infty(0,T)\) be a
nonnegative function satisfying
\[
F(t)\le\Lambda\int_0^tF(s)\,ds\qquad\text{for a.e. }t\in[0,T].
\]
Then \(F(t)=0\) for a.e.\ \(t\in[0,T]\).
\end{lemma}

\begin{proof}
Set \(G(t)\eqdef\int_0^tF(s)\,ds\) for \(t\in[0,T]\). Then \(G\) is absolutely
continuous, \(G'=F\) a.e., and the assumption gives \(G'\le\Lambda G\) a.e. We
infer that
\[
\frac{d}{dt}\bigl(e^{-\Lambda t}G(t)\bigr)
=e^{-\Lambda t}\bigl(G'(t)-\Lambda G(t)\bigr)\le0
\qquad\text{a.e. on }(0,T),
\]
so that \(t\mapsto e^{-\Lambda t}G(t)\) is nonincreasing. Since \(G(0)=0\) and
\(G\ge0\), we get \(G\equiv0\) on \([0,T]\), hence \(F=G'=0\) a.e.
\end{proof}

\begin{proof}[Proof of Lemma~\ref{lem:comparison}]
By Assumption~\ref{ass:IL}(C) and \eqref{eq:monotone}, it suffices to prove that
\[
u_\eta(t,\cdot)\le v\qquad\text{for every }\eta\in(0,1)\text{ and every }t>0,
\]
where \(u_\eta\) is the global solution of \eqref{eq:Reta}. Fix
\(\eta\in(0,1)\) and \(T>0\), and set \(K\eqdef K_{\eta,T}\) and
\(M_T\eqdef M_{\eta,T}\), which are finite by \eqref{eq:ueta-infty} and
\eqref{eq:ueta-one}.

Let \(\Omega=\mathbb{R}^N\setminus\{0\}\) in the first alternative of (H2) and
\(\Omega=\mathbb{R}^N\) in the second one. By the first inequality in
\eqref{eq:lip}, we have \(0\le g_\eta(x,v)\le\abs{x}^{\sigma}v^{p}\) for all
\(x\in\Omega\). Since \(v\) is time independent, \eqref{eq:supersol} gives
\begin{equation}\label{eq:vsuper-eta}
\int_{\Omega}\Bigl(\nabla v^m\cdot\nabla\varphi-g_\eta(x,v)\,\varphi\Bigr)\,dx\ \ge\ 0
\qquad\text{for all }0\le\varphi\in C^\infty_c(\Omega).
\end{equation}
The sign \(\sigma<0\) is used here: a supersolution of the \emph{singular}
stationary equation is automatically a supersolution of \emph{every}
regularized problem \eqref{eq:Reta}.

Assume first that \(v(x)\to+\infty\) as \(x\to0\), and choose
\(\delta\in(0,1)\) such that \(v(x)\ge K\) for all
\(x\in B(0,\delta)\setminus\{0\}\). Then, for every \(t\in[0,T]\),
\begin{equation}\label{eq:awayfromzero}
\bigl(u_\eta(t,\cdot)-v\bigr)_+(x)=0\qquad\text{for a.e. }x\in B(0,\delta).
\end{equation}
In this first alternative we also set
\(V\eqdef\min\{v^m,K^m\}\). Since \(v^m\in H^1_{\mathrm{loc}}(\RN\setminus\{0\})\)
and \(V=K^m\) on \(B(0,\delta)\setminus\{0\}\), we have
\(V\in H^1_{\mathrm{loc}}(\RN)\). In the second alternative of (H2), \(v\) is
bounded on \(B(0,1)\) with \(v^m\in H^1(B(0,1))\), and it is a supersolution on
all of \(\mathbb{R}^N\). We then set \(\delta=0\) and \(V\eqdef v^m\in
H^1_{\mathrm{loc}}(\RN)\). In both cases, since \(u_\eta\le K\) on
\([0,T]\times\RN\),
\begin{equation}\label{eq:Vtrunc}
H\bigl(u_\eta^m-v^m\bigr)=H\bigl(u_\eta^m-V\bigr)\quad\text{a.e. in }(0,T)\times\RN
\end{equation}
for every function \(H\) vanishing on \((-\infty,0]\).

Since \(0\le(u_\eta-v)_+\le u_\eta\) and \(0\le(u_\eta^m-v^m)_+\le u_\eta^m\), we
obtain, for every \(t\in[0,T]\),
\[
\|(u_\eta(t)-v)_+\|_1\le M_T
\qquad\text{and}\qquad
\|(u_\eta^m(t)-v^m)_+\|_1\le\|u_\eta(t)\|_\infty^{m-1}\|u_\eta(t)\|_1
\le K^{m-1}M_T.
\]

Let \(\chi\in C_c^\infty(\mathbb{R}^N)\) satisfy \(0\le\chi\le1\),
\(\chi\equiv1\) on \(B(0,1)\) and \(\supp\chi\subset B(0,2)\), and set
\(\chi_R(x)\eqdef\chi(x/R)\) for \(R\ge1\). Let
\(\{H_\varepsilon\}_{\varepsilon>0}\subset C^{\infty}(\mathbb{R})\) be
nondecreasing approximations of \(\sgn_+\) with bounded derivatives of all
orders and
\[
H_\varepsilon\equiv0\text{ on }(-\infty,0],\qquad 0\le H_\varepsilon\le1,
\qquad H_\varepsilon'\ge0,\qquad
H_\varepsilon(s)\uparrow\mathbf 1_{(0,\infty)}(s)\text{ as }\varepsilon\downarrow0 .
\]
We define the test function
\[
\varphi\eqdef H_\varepsilon\bigl(u_\eta^m-v^m\bigr)\chi_R
=H_\varepsilon\bigl(u_\eta^m-V\bigr)\chi_R ,
\]
the equality being \eqref{eq:Vtrunc}. Then \(\varphi\ge0\) and
\(\|\varphi\|_\infty\le1\). Moreover, \(u_\eta^m-V\in H^1_{\mathrm{loc}}(\RN)\)
for a.e.\ \(t\), since \(u_\eta^m\in L^2(0,T;H^1_{\mathrm{loc}}(\RN))\) and
\(V\in H^1_{\mathrm{loc}}(\RN)\). Hence \(\varphi\in L^2(0,T;H^1(\RN))\),
because \(H_\varepsilon\) is smooth with \(H_\varepsilon'\in L^\infty(\mathbb{R})\)
and \(H_\varepsilon(0)=0\). By \eqref{eq:awayfromzero}, \(\varphi\) vanishes
a.e.\ on \(B(0,\delta)\), and \(\supp\varphi(t)\subset\{\delta\le\abs{x}\le2R\}\).
In particular, \(\varphi(t)\) has compact support in \(\Omega\): in the first
alternative of (H2) it vanishes on \(B(0,\delta)\), and in the second one
\(\Omega=\RN\). Hence \(\varphi(t)\) is an admissible test function in
\eqref{eq:vsuper-eta} for a.e.\ \(t\), after the density argument which
follows Definition~\ref{def:supersol}, and \(\varphi\) is admissible in the
first-order weak formulation of \eqref{eq:Reta}.

Subtracting \eqref{eq:vsuper-eta} from the weak formulation of \eqref{eq:Reta}
tested against \(\varphi\), we obtain, for a.e.\ \(t\in(0,T)\),
\begin{align*}
\frac{d}{dt}\int_{\mathbb{R}^N}\Phi_\varepsilon\bigl(u_\eta-v\bigr)\chi_R\,dx
\le\ &
\underbrace{-\int_{\mathbb{R}^N}H'_\varepsilon\bigl(u_\eta^m-v^m\bigr)
\bigl\lvert\nabla\bigl(u_\eta^m-v^m\bigr)\bigr\rvert^{2}\chi_R\,dx}_{=:J_1}\\
&+\underbrace{\Bigl(-\int_{\mathbb{R}^N}H_\varepsilon\bigl(u_\eta^m-v^m\bigr)\,
\nabla\bigl(u_\eta^m-v^m\bigr)\cdot\nabla\chi_R\,dx\Bigr)}_{=:J_2}\\
&+\underbrace{\int_{\mathbb{R}^N}\bigl[g_\eta(x,u_\eta)-g_\eta(x,v)\bigr]\,
H_\varepsilon\bigl(u_\eta^m-v^m\bigr)\chi_R\,dx}_{=:J_3} .
\end{align*}
Here \(\Phi_\varepsilon\) denotes, for each fixed \(x\), the primitive
\[
\Phi_\varepsilon(\varsigma)\eqdef\int_0^{\varsigma}
H_\varepsilon\bigl((\tau+v(x))^m-v(x)^m\bigr)\,d\tau,\qquad\varsigma\ge-v(x),
\]
which satisfies \(0\le\Phi_\varepsilon(\varsigma)\le\varsigma_+\) and
\(\Phi_\varepsilon(\varsigma)\to\varsigma_+\) as \(\varepsilon\downarrow0\).
This is the classical Kato computation. Since \(v\) does not depend on time,
the only delicate point is the chain rule
\begin{equation}\label{eq:chainrule}
\int_0^T\bigl\langle\partial_tu_\eta,H_\varepsilon(u_\eta^m-V)\chi_R\bigr\rangle\,\theta\,dt
=-\int_0^T\!\!\int_{\RN}\Phi_\varepsilon(u_\eta-v)\,\chi_R\,\theta'\dx\,dt,
\qquad \theta\in C^\infty_c((0,T)),
\end{equation}
where we used \eqref{eq:Vtrunc} once more. We justify it with the chain rule
of Otto \cite[Lemma~1]{Otto96}, which is a refinement of the argument of Alt
and Luckhaus \cite[Lemma~1.5]{AltLuckhaus83}; see also
\cite[Lemma~4]{Carrillo99} for the analogous statement with comparison
functions which do not depend on \(x\). In the notation of
\cite{Otto96}, the unknown is \(z\eqdef u_\eta^m\), the monotone function is
\(b(z)=z^{1/m}\), and the entropy is the convex function \(\Theta_\varepsilon\)
with \(\Theta_\varepsilon'=H_\varepsilon\), whose first and second derivatives are
bounded and continuous. For a time-independent function
\(v_0=V\in H^1_{\mathrm{loc}}(\RN)\), the corresponding entropy flux is
\[
q(z,V)=\int_{V}^{z}H_\varepsilon(\zeta-V)\,d\bigl(\zeta^{1/m}\bigr)
=\int_0^{u_\eta-v}H_\varepsilon\bigl((\tau+v)^m-v^m\bigr)\,d\tau
=\Phi_\varepsilon(u_\eta-v)
\]
on the set \(\{u_\eta>v\}\), after the change of variables
\(\zeta=(\tau+v)^m\), while both sides vanish on \(\{u_\eta\le v\}\). Part (a)
of \cite[Lemma~1]{Otto96} applied to the nondecreasing function \(\Theta_\varepsilon\), and part (b) applied
to the nonincreasing function \(-\Theta_\varepsilon\), give the two inequalities which together are
\eqref{eq:chainrule}, with \(\theta\chi_R\) as nonnegative test function.
Otto's argument requires no time regularity of \(b(z)\) beyond what is
contained in the weak formulation. On a bounded open set
\(\Omega'\supset\supp\chi_R\), it only uses that \(z\in L^2(0,T;H^1(\Omega'))\),
that the energy of \(b(z)\) is bounded in \(L^\infty(0,T;L^1(\Omega'))\), and
the weak inequality \cite[Eq.~(14)]{Otto96} tested with the Steklov averages
\(h^{-1}\int_t^{t+h}\Theta_\varepsilon'(z-V)\,\theta\chi_R\,d\tau\), \(h>0\).
These test functions are admissible because \(\Theta_\varepsilon'\) is bounded
and \(\Theta_\varepsilon''\) is bounded and continuous. Since \(\chi_R\) has
compact support in \(\Omega'\), the boundary condition \cite[Eq.~(19)]{Otto96}
plays no role. These properties hold for \(u_\eta\), because
\(\partial_tu_\eta=\Delta u_\eta^m+g_\eta(x,u_\eta)\) in the weak sense, with
\(\nabla u_\eta^m\in L^2((0,T)\times\RN)\), \(u_\eta\in L^\infty((0,T),L^1\cap
L^\infty(\RN))\) and \(0\le g_\eta(x,u_\eta)\le\eta^{\sigma-1}K\). The reaction
term enters the argument only as a fixed bounded function multiplied by the
test function. In particular, the
\(x\)-dependence of the reaction term plays no role. Moreover, the function
\(t\mapsto\int\Phi_\varepsilon(u_\eta(t)-v)\chi_R\dx\) is continuous on
\([0,T]\), since \(u_\eta\in C([0,T],\Lp1)\) and \(\Phi_\varepsilon\) is
\(1\)-Lipschitz. Hence \eqref{eq:chainrule} yields the differential inequality
above in \(\mathcal D'((0,T))\), together with the value of the left-hand side
at \(t=0\).

We now estimate the three terms \(J_1\), \(J_2\), \(J_3\). The term \(J_1\) is
nonpositive because \(H_\varepsilon'\ge0\), and it can be discarded.

For the term \(J_3\), we use Lemma~\ref{lem:lipschitz}, together with the fact
that \(H_\varepsilon(u_\eta^m-v^m)\ne0\) forces \(u_\eta>v\). By
\eqref{eq:lip}, we get
\[
\bigl[g_\eta(x,u_\eta)-g_\eta(x,v)\bigr]H_\varepsilon(u_\eta^m-v^m)
\le\Lambda_\eta\abs{u_\eta-v}H_\varepsilon(u_\eta^m-v^m)
\le\Lambda_\eta(u_\eta-v)_+,
\]
and therefore
\[
J_3\le\Lambda_\eta\int_{\mathbb{R}^N}(u_\eta-v)_+\chi_R\,dx .
\]
For the term \(J_2\), we recall that \(\nabla\chi_R\) is supported in the
annulus \(\{R\le\abs{x}\le2R\}\), on which \(u_\eta^m-v^m\in H^1\). Letting
\(\varepsilon\downarrow0\) by dominated convergence, and integrating by parts,
we obtain
\begin{align*}
J_2\xrightarrow[\varepsilon\downarrow0]{}
&-\int_{\mathbb{R}^N}\mathbf 1_{\{u_\eta^m>v^m\}}\,
\nabla\bigl(u_\eta^m-v^m\bigr)\cdot\nabla\chi_R\,dx\\
&=-\int_{\mathbb{R}^N}\nabla\bigl(u_\eta^m-v^m\bigr)_+\cdot\nabla\chi_R\,dx
=\int_{\mathbb{R}^N}\bigl(u_\eta^m-v^m\bigr)_+\,\Delta\chi_R\,dx .
\end{align*}
No boundary terms appear in the last integration by parts, because
\(\nabla\chi_R\) has compact support in the annulus, on which
\((u_\eta^m-v^m)_+\in W^{1,1}\). Since \(\Delta\chi_R\) may change sign, we
bound the limit by
\(\int_{\mathbb{R}^N}\bigl(u_\eta^m-v^m\bigr)_+\,\abs{\Delta\chi_R}\,dx\).

We now let \(\varepsilon\downarrow0\) in the differential inequality. We have
\(\Phi_\varepsilon(u_\eta-v)\to(u_\eta-v)_+\) in \(L^1(\mathbb{R}^N)\) by
dominated convergence, since \(0\le\Phi_\varepsilon(u_\eta-v)\le(u_\eta-v)_+\le
u_\eta\) and \(u_\eta\in L^1\). Inserting the above bounds on \(J_1\), \(J_2\),
\(J_3\), we arrive at
\[
\frac{d}{dt}\int_{\mathbb{R}^N}\bigl(u_\eta-v\bigr)_+\chi_R\,dx
\ \le\ \int_{\mathbb{R}^N}\bigl(u_\eta^m-v^m\bigr)_+\,|\Delta\chi_R|\,dx
+\Lambda_\eta\int_{\mathbb{R}^N}\bigl(u_\eta-v\bigr)_+\chi_R\,dx
\qquad\text{in }\mathcal{D}'\bigl((0,T)\bigr).
\]

We integrate this inequality on \((0,t)\). The initial term vanishes, because
\(u_\eta(0)=u_0\mathbf 1_{B(0,1/\eta)}\le u_0\le v\) by (H3). Since
\(\Delta\chi_R(x)=R^{-2}(\Delta\chi)(x/R)\), we have
\(\|\Delta\chi_R\|_\infty\le CR^{-2}\). Using \((u_\eta^m-v^m)_+\le u_\eta^m\),
we deduce that
\[
\int_{\mathbb{R}^N}(u_\eta^m-v^m)_+\,|\Delta\chi_R|\,dx
\le CR^{-2}\int_{\mathbb{R}^N}u_\eta^m\,dx
\le CR^{-2}K^{m-1}M_T,
\]
uniformly in \(t\in[0,T]\). We let \(R\to\infty\), using
\(\chi_R(x)\to1\) for every \(x\in\mathbb{R}^N\) and dominated convergence with
the integrable majorant \(u_\eta\). We obtain
\[
F(t)\ \le\ \Lambda_\eta\int_0^tF(s)\,ds,\qquad t\in[0,T],
\qquad\text{where}\qquad
F(t)\eqdef\bigl\lVert\bigl(u_\eta(t)-v\bigr)_+\bigr\rVert_{1}\ge0 .
\]
Since \(F\in L^\infty(0,T)\), Lemma~\ref{lem:gronwall} applied with
\(\Lambda=\Lambda_\eta\) yields \(F(t)=0\) for a.e.\ \(t\in[0,T]\), and then
for every \(t\in[0,T]\) by the continuity of \(t\mapsto u_\eta(t)\) in
\(L^1\). In other words, \(u_\eta(t)\le v\) a.e.\ on
\(\mathbb{R}^N\setminus\{0\}\). As \(T>0\) and \(\eta\in(0,1)\) are arbitrary,
\eqref{eq:monotone} gives \eqref{eq:comparison}.
\end{proof}
\begin{remark}\rm \label{rem:why-eta}
Two structural facts make Lemma~\ref{lem:comparison} work. First, the
regularized reaction \(g_\eta\) is \emph{globally} Lipschitz in \(s\)
(Lemma~\ref{lem:lipschitz}), so that no local boundedness of the barrier is
needed in the Gronwall step. Second, we have
\((\abs{x}^2+\eta^2)^{\sigma/2}\le\abs{x}^{\sigma}\) because \(\sigma<0\), so
that a supersolution of the singular stationary equation is automatically a
supersolution of every regularized problem \eqref{eq:Reta}. Both facts would
fail for \(\sigma>0\). This is one place where the Hardy case is more
favourable than the H\'enon case.
\end{remark}

The comparison argument also applies to bounded time-dependent
supersolutions, provided their time derivative is a locally integrable
function. In that case the chain rule of Otto is not needed: the time
derivative of \(u_\eta\) is handled by backward Steklov averages and by the
convexity of the Kato entropy, in the spirit of \cite[Lemma~1.5]{AltLuckhaus83}
and \cite[Lemma~4]{Carrillo99}. This extension is only used in
Section~\ref{sec:conjectures}.

\begin{lemma}\label{lem:comparisonU}
Let \(u_0\in\mathcal A_r\), \(r>r_c\), and let \(u\) be the solution of
\eqref{eq:CP}. Let \(T>0\) and let \(U:[0,T)\times\RN\to(0,\infty)\) be a
measurable function such that, for every \(T'\in(0,T)\),
\begin{enumerate}[label=\textup{(U\arabic*)},leftmargin=2.6em]
\item \(U\in L^\infty((0,T')\times\RN)\cap C([0,T'],L^1_{\mathrm{loc}}(\RN))\),
\(U^m\in L^2(0,T';H^1_{\mathrm{loc}}(\RN))\), and the distributional time
derivative satisfies \(\partial_tU\in L^1_{\mathrm{loc}}([0,T)\times\RN)\);
\item for every nonnegative \(\varphi\in C_c^\infty((0,T)\times\RN)\),
\[
\int_0^T\!\!\int_{\RN}\Bigl(\partial_tU\,\varphi+\nabla U^m\cdot\nabla\varphi
-\abs{x}^{\sigma}U^{p}\varphi\Bigr)\dx\,dt\ \ge\ 0 ;
\]
\item \(u_0\le U(0)\) a.e.\ in \(\RN\), where \(U(0)\) is the value at \(t=0\)
of the continuous representative \(t\mapsto U(t)\in L^1_{\mathrm{loc}}(\RN)\).
\end{enumerate}
Then \(u\le U\) a.e.\ in \((0,\min\{T,\Tmax(u_0)\})\times\RN\).
\end{lemma}

\begin{proof}
By \eqref{eq:monotone}, it suffices to show that \(u_\eta\le U\) a.e.\ in
\((0,T')\times\RN\) for every \(\eta\in(0,1)\) and every \(T'\in(0,T)\). We fix
\(\eta\) and \(T'\), and we write \(u\eqdef u_\eta\), \(z\eqdef u_\eta^m\),
\(g(x,s)\eqdef g_\eta(x,s)\), \(K\eqdef K_{\eta,T'}\) and \(M\eqdef M_{\eta,T'}\),
see \eqref{eq:ueta-infty}--\eqref{eq:ueta-one}. All the integrals below are
taken over compact subsets of \((0,T')\times\RN\), on which every function we
meet is integrable.

\smallskip
\noindent\emph{Step 1: preliminaries.} Since \(u\in C([0,T'],\Lp1)\) and
\(\nabla z\in L^2((0,T')\times\RN)\), the first-order weak formulation of
\eqref{eq:Reta}, tested with \(\rho(t)\phi(x)\), where \(\rho\) is a piecewise
linear approximation of \(\mathbf 1_{[s,t]}\) and \(\phi\in C^1_c(\RN)\), gives
\begin{equation}\label{eq:Wint}
\int_{\RN}\bigl(u(t)-u(s)\bigr)\phi\dx
=\int_s^t\!\!\int_{\RN}\Bigl(-\nabla z\cdot\nabla\phi+g(x,u)\,\phi\Bigr)\dx\,d\tau,
\qquad 0\le s<t\le T'.
\end{equation}
Since \(u\) is bounded, \eqref{eq:Wint} extends by density to every
\(\phi\in H^1(\RN)\) with compact support. Next, since \(\partial_tU\in
L^1_{\mathrm{loc}}\), for a.e.\ \(x\) the function \(t\mapsto U(t,x)\) is
absolutely continuous on \([0,T']\), with derivative \(\partial_tU(\cdot,x)\),
and it coincides with the continuous representative of (U1).

We take \(H_\varepsilon(s)\eqdef H(s/\varepsilon)\), where \(H\in C^\infty(\R)\)
is nondecreasing, \(H=0\) on \((-\infty,0]\) and \(H=1\) on \([1,\infty)\).
Thus \(H_\varepsilon\) has the properties used in the proof of
Lemma~\ref{lem:comparison}; moreover \(H_\varepsilon'\ge0\) is supported in
\([0,\varepsilon]\), \(\int_\R H_\varepsilon'=1\), and
\(H_\varepsilon\uparrow\mathbf 1_{(0,\infty)}\) as \(\varepsilon\downarrow0\).
For \(w\ge0\) and \(V>0\), we set
\[
G_\varepsilon(w;V)\eqdef\int_0^{w}H_\varepsilon\bigl(r^m-V^m\bigr)\,dr .
\]
The function \(G_\varepsilon\) has the following properties.
\begin{enumerate}[label=\textup{(G\arabic*)},leftmargin=2.6em]
\item \(w\mapsto G_\varepsilon(w;V)\) is convex on \([0,\infty)\), with
\(\partial_wG_\varepsilon(w;V)=H_\varepsilon(w^m-V^m)\); moreover
\(0\le G_\varepsilon(w;V)\le(w-V)_+\), and
\(G_\varepsilon(w;V)\uparrow(w-V)_+\) as \(\varepsilon\downarrow0\).
\item \(G_\varepsilon\) is of class \(C^1\) on \([0,\infty)\times(0,\infty)\),
and the change of variables \(\zeta=r^m-V^m\) gives
\[
\partial_VG_\varepsilon(w;V)=-mV^{m-1}\int_0^w H_\varepsilon'(r^m-V^m)\,dr
=-\int_0^{(w^m-V^m)_+}H_\varepsilon'(\zeta)\Bigl(\frac{V^m}{V^m+\zeta}\Bigr)^{\frac{m-1}{m}}d\zeta
\in[-1,0].
\]
\item Consequently, since \(H_\varepsilon(w^m-V^m)=\int_0^{(w^m-V^m)_+}H_\varepsilon'(\zeta)\,d\zeta\),
\[
0\le k_\varepsilon(w;V)\eqdef\partial_VG_\varepsilon(w;V)+H_\varepsilon(w^m-V^m)
\le1-\Bigl(\frac{V^m}{V^m+\varepsilon}\Bigr)^{\frac{m-1}{m}}
\xrightarrow[\varepsilon\downarrow0]{}0 .
\]
\end{enumerate}
Indeed, (G1) follows from the monotonicity of \(H_\varepsilon\) and of
\(r\mapsto r^m\), and from \(H_\varepsilon(r^m-V^m)\le\mathbf 1_{\{r>V\}}\);
(G2) and (G3) follow from \(H_\varepsilon'\ge0\),
\(\operatorname{supp}H_\varepsilon'\subset[0,\varepsilon]\) and
\(\int_\R H_\varepsilon'=1\).

\smallskip
\noindent\emph{Step 2: a Kato inequality.} Let \(\chi_R\) be as in the proof
of Lemma~\ref{lem:comparison}, let \(0\le\theta\in C^\infty_c((0,T'))\), and set
\[
w\eqdef z-U^m,\qquad \psi\eqdef H_\varepsilon(w)\,\chi_R\,\theta .
\]
Since \(z\) and \(U^m\) belong to \(L^2(0,T';H^1_{\mathrm{loc}}(\RN))\), and
since \(H_\varepsilon\) has bounded derivatives, the function \(\psi\) belongs
to \(L^2(0,T';H^1(\RN))\). Moreover, \(0\le\psi\le1\), and \(\psi\) is supported
in a fixed compact subset of \((0,T')\times\RN\). Let \(h_0>0\) be such that
\(\supp\theta\subset(2h_0,T'-2h_0)\), and let \(h\in(0,h_0)\). We use the
backward difference quotient and the backward average
\[
\partial^{-h}u(t)\eqdef\frac{u(t)-u(t-h)}{h},\qquad
\mathcal A_hf(t)\eqdef\frac1h\int_{t-h}^{t}f(\tau)\,d\tau .
\]
Applying \eqref{eq:Wint} with \(s=t-h\) and \(\phi=\psi(t)\), and integrating in
\(t\), we obtain
\begin{equation}\label{eq:steklov}
\int_0^{T'}\!\!\int_{\RN}\partial^{-h}u\,\psi\dx\,dt
=\int_0^{T'}\!\!\int_{\RN}\Bigl(-\mathcal A_h(\nabla z)\cdot\nabla\psi
+\mathcal A_h\bigl(g(\cdot,u)\bigr)\,\psi\Bigr)\dx\,dt .
\end{equation}
We now bound the left-hand side from below. Since \(\psi=\partial_wG_\varepsilon(u;U)\,\chi_R\theta\)
by (G1), the convexity in (G1) gives, a.e.\ in \((0,T')\times\RN\),
\[
\bigl(u(t)-u(t-h)\bigr)\,\partial_wG_\varepsilon\bigl(u(t);U(t)\bigr)
\ge G_\varepsilon\bigl(u(t);U(t)\bigr)-G_\varepsilon\bigl(u(t-h);U(t)\bigr).
\]
By (G2) and the absolute continuity of \(U(\cdot,x)\),
\[
G_\varepsilon\bigl(u(t-h);U(t)\bigr)
=G_\varepsilon\bigl(u(t-h);U(t-h)\bigr)
+\int_{t-h}^{t}\partial_VG_\varepsilon\bigl(u(t-h);U(\tau)\bigr)\,\partial_tU(\tau)\,d\tau .
\]
We multiply by \(h^{-1}\chi_R\theta\ge0\) and integrate. With
\(a(t)\eqdef\int_{\RN}G_\varepsilon(u(t);U(t))\chi_R\dx\), which is bounded, we
infer from \eqref{eq:steklov} that
\begin{align*}
-\int_0^{T'}a(t)\,\frac{\theta(t+h)-\theta(t)}{h}\,dt
&-\int_0^{T'}\!\!\int_{\RN}\Bigl(\frac1h\int_{t-h}^{t}\partial_VG_\varepsilon\bigl(u(t-h);U(\tau)\bigr)\,
\partial_tU(\tau)\,d\tau\Bigr)\chi_R\,\theta\dx\,dt\\
&\le\int_0^{T'}\!\!\int_{\RN}\Bigl(-\mathcal A_h(\nabla z)\cdot\nabla\psi
+\mathcal A_h\bigl(g(\cdot,u)\bigr)\,\psi\Bigr)\dx\,dt ,
\end{align*}
where we used the discrete integration by parts
\(h^{-1}\int(a(t)-a(t-h))\theta(t)\,dt=-\int a(t)\,h^{-1}(\theta(t+h)-\theta(t))\,dt\).
We let \(h\downarrow0\). The first term converges to \(-\int a\,\theta'\,dt\).
On the right-hand side, \(\mathcal A_h(\nabla z)\to\nabla z\) in
\(L^2_{\mathrm{loc}}\), and \(\mathcal A_h(g(\cdot,u))\to g(\cdot,u)\) in
\(L^1_{\mathrm{loc}}\) with \(\abs{\mathcal A_h(g(\cdot,u))}\le\eta^{\sigma-1}K\).
In the second term, we write the inner average as
\[
\partial_VG_\varepsilon\bigl(u(t);U(t)\bigr)\,\mathcal A_h(\partial_tU)(t)
+\int_0^1\Bigl[\partial_VG_\varepsilon\bigl(u(t-h);U(t-\lambda h)\bigr)
-\partial_VG_\varepsilon\bigl(u(t);U(t)\bigr)\Bigr]\partial_tU(t-\lambda h)\,d\lambda .
\]
The first part converges to \(\partial_VG_\varepsilon(u;U)\,\partial_tU\) in
\(L^1_{\mathrm{loc}}\), because \(\abs{\partial_VG_\varepsilon}\le1\) and
\(\mathcal A_h(\partial_tU)\to\partial_tU\) in \(L^1_{\mathrm{loc}}\). In the
second part, the bracket is bounded by \(2\) and tends to \(0\) in measure on
compact sets, by the continuity of \(\partial_VG_\varepsilon\) and since
\(u(t-h)\to u(t)\) and \(U(t-\lambda h)\to U(t)\) in \(L^1_{\mathrm{loc}}\),
while \(\partial_tU(t-\lambda h)\to\partial_tU(t)\) in \(L^1_{\mathrm{loc}}\) by
the continuity of translations. Since \(\abs{a_hb_h}\le2\abs{b_h-b}+\abs{a_h}\abs{b}\)
whenever \(\abs{a_h}\le2\), the second part tends to \(0\) in
\(L^1_{\mathrm{loc}}\). We arrive at
\begin{equation}\label{eq:katoU1}
-\int_0^{T'}\!\!\int_{\RN}G_\varepsilon(u;U)\,\chi_R\,\theta'
-\int_0^{T'}\!\!\int_{\RN}\partial_VG_\varepsilon(u;U)\,\partial_tU\,\chi_R\,\theta
\le\int_0^{T'}\!\!\int_{\RN}\Bigl(-\nabla z\cdot\nabla\psi+g(x,u)\,\psi\Bigr).
\end{equation}

\smallskip
\noindent\emph{Step 3: use of the supersolution.} Since \(\psi\ge0\) is bounded,
belongs to \(L^2(0,T';H^1(\RN))\) and has compact support in
\((0,T')\times\RN\), mollification in \((t,x)\) provides nonnegative
\(\varphi_k\in C^\infty_c((0,T)\times\RN)\), with supports in a fixed compact
set and \(0\le\varphi_k\le1\), such that \(\varphi_k\to\psi\) in
\(L^2(0,T';H^1(\RN))\) and a.e. Since \(\partial_tU\) and
\(\abs{x}^{\sigma}U^p\) are locally integrable and \(\nabla U^m\) is locally
square integrable, we may pass to the limit in (U2) and obtain (U2) with
\(\varphi=\psi\). As \(g(x,U)\le\abs{x}^{\sigma}U^p\) by \eqref{eq:lip}, we get
\begin{equation}\label{eq:katoU2}
0\le\int_0^{T'}\!\!\int_{\RN}\Bigl(\partial_tU\,\psi+\nabla U^m\cdot\nabla\psi
-g(x,U)\,\psi\Bigr).
\end{equation}
Adding \eqref{eq:katoU1} and \eqref{eq:katoU2}, and recalling that
\(\psi=H_\varepsilon(w)\chi_R\theta\) and the definition of \(k_\varepsilon\) in
(G3), we obtain
\begin{align*}
-\int_0^{T'}\!\!\int_{\RN}G_\varepsilon(u;U)\,\chi_R\,\theta'
&\le\int_0^{T'}\!\!\int_{\RN}k_\varepsilon(u;U)\,\partial_tU\,\chi_R\theta
-\int_0^{T'}\!\!\int_{\RN}\nabla w\cdot\nabla\psi\\
&\quad+\int_0^{T'}\!\!\int_{\RN}\bigl(g(x,u)-g(x,U)\bigr)\psi .
\end{align*}

\smallskip
\noindent\emph{Step 4: conclusion.} As in the proof of
Lemma~\ref{lem:comparison},
\(-\nabla w\cdot\nabla\psi=-H_\varepsilon'(w)\abs{\nabla w}^2\chi_R\theta
-H_\varepsilon(w)\nabla w\cdot\nabla\chi_R\,\theta
\le-H_\varepsilon(w)\nabla w\cdot\nabla\chi_R\,\theta\), and
\((g(x,u)-g(x,U))H_\varepsilon(w)\le\Lambda_\eta(u-U)_+\) by \eqref{eq:lip},
since \(H_\varepsilon(w)\ne0\) forces \(u>U\). We let
\(\varepsilon\downarrow0\). By (G1) and monotone convergence,
\(G_\varepsilon(u;U)\to(u-U)_+\). By (G3), \(\abs{k_\varepsilon\partial_tU}\le\abs{\partial_tU}\)
and \(k_\varepsilon\to0\) everywhere, because \(U>0\); hence the first term on
the right-hand side tends to \(0\) by dominated convergence. Finally,
\(-\int\!\!\int H_\varepsilon(w)\nabla w\cdot\nabla\chi_R\,\theta\to
-\int\!\!\int\nabla w_+\cdot\nabla\chi_R\,\theta=\int\!\!\int w_+\Delta\chi_R\,\theta\).
Setting \(F_R(t)\eqdef\int_{\RN}(u(t)-U(t))_+\chi_R\dx\), we conclude that
\[
-\int_0^{T'}F_R\,\theta'\,dt\le\int_0^{T'}\Bigl(\int_{\RN}\bigl(z-U^m\bigr)_+\abs{\Delta\chi_R}\dx
+\Lambda_\eta F_R\Bigr)\theta\,dt
\qquad\text{for all }0\le\theta\in C^\infty_c((0,T')).
\]
The function \(F_R\) is continuous on \([0,T']\), because
\(u\in C([0,T'],\Lp1)\), \(U\in C([0,T'],L^1_{\mathrm{loc}}(\RN))\) and
\(\chi_R\) has compact support. Moreover \(F_R(0)=0\), since
\(u(0)=u_0\mathbf 1_{B(0,1/\eta)}\le u_0\le U(0)\) by (U3). A continuous
function whose distributional derivative is bounded above by an integrable
function \(b\) satisfies \(F_R(t)-F_R(0)\le\int_0^tb\). Using
\((z-U^m)_+\le z\le K^{m-1}u\) and \(\norm{\Delta\chi_R}_\infty\le CR^{-2}\), we
deduce that
\[
F_R(t)\le CR^{-2}K^{m-1}MT'+\Lambda_\eta\int_0^tF_R(s)\,ds,\qquad t\in[0,T'].
\]
Let \(F(t)\eqdef\norm{(u(t)-U(t))_+}_1\le M\). Letting \(R\to\infty\) by
monotone convergence, we obtain \(F(t)\le\Lambda_\eta\int_0^tF(s)\,ds\) for
\(t\in[0,T']\).
Lemma~\ref{lem:gronwall} gives \(F=0\) on \([0,T']\), that is, \(u_\eta\le U\)
a.e.\ in \((0,T')\times\RN\). This completes the proof.
\end{proof}

The obstruction described in Remark~\ref{rem:origin} shows that a bounded
barrier with critical tail cannot be obtained by naively smoothing \(S\). We
therefore build it by gluing \(S\) to an explicit inner profile whose
Laplacian reproduces the Hardy weight. The first ingredient is the following
bounded stationary supersolution in the unit ball.

\begin{proposition}\label{prop:inner}
Let \(N\ge3\), \(m>1\), \(\sigma\in(-2,0)\), \(p>m\), and set
\[
q\eqdef\frac{p}{m}>1,\qquad D\eqdef(\sigma+2)(N+\sigma)>0,
\qquad c\eqdef\Bigl(\frac{D}{2}\Bigr)^{\frac{1}{q-1}} .
\]
Then the function
\begin{equation}\label{eq:inner}
w(x)\eqdef c\Bigl(1-\tfrac12\abs{x}^{\sigma+2}\Bigr),\qquad x\in B(0,1),
\end{equation}
satisfies \(\tfrac{c}{2}\le w\le c\) on \(B(0,1)\), \(w\in H^1(B(0,1))\), and
\begin{equation}\label{eq:inner-super}
\Delta w+\abs{x}^{\sigma}w^{\,q}\ \le\ 0\qquad\text{in }\Dprime{B(0,1)} .
\end{equation}
Equivalently, \(\Psi\eqdef w^{1/m}\) is a bounded stationary weak
supersolution of \eqref{eq:CP} in \(B(0,1)\), with
\(\bigl(\tfrac{c}{2}\bigr)^{1/m}\le\Psi\le c^{1/m}\).
\end{proposition}

\begin{proof}
Since \(\sigma+2\in(0,2)\), we have \(\tfrac12\abs{x}^{\sigma+2}\le\tfrac12\) on
\(B(0,1)\), whence \(\tfrac{c}{2}\le w\le c\). In particular, \(w\) is positive
and bounded on \(B(0,1)\). Moreover
\(\abs{\nabla w}=\tfrac{c}{2}(\sigma+2)\abs{x}^{\sigma+1}\), so that
\(\abs{\nabla w}^{2}\) is a multiple of \(\abs{x}^{2\sigma+2}\). This function is
integrable on \(B(0,1)\) because \(2\sigma+2>-2>-N\), and we conclude that
\(w\in H^1(B(0,1))\).

By Lemma~\ref{lem:dirac}(ii) with \(b=\sigma+2>0\), we have
\[
\Delta w=-\frac{c}{2}\,(\sigma+2)(N+\sigma)\,\abs{x}^{\sigma}
=-\frac{c\,D}{2}\,\abs{x}^{\sigma}\qquad\text{in }\Dprime{\RN},
\]
with no Dirac mass at the origin. Note that \(\abs{x}^{\sigma}\in
L^1_{\mathrm{loc}}(\RN)\), because \(\sigma>-N\). Using \(w\le c\) and the
definition of \(c\), we obtain
\[
\Delta w+\abs{x}^{\sigma}w^{q}
\le\abs{x}^{\sigma}\Bigl(c^{q}-\frac{cD}{2}\Bigr)
=c\,\abs{x}^{\sigma}\Bigl(c^{q-1}-\frac{D}{2}\Bigr)=0 .
\]
Both sides are locally integrable functions on \(B(0,1)\), so that this
inequality holds in \(\Dprime{B(0,1)}\), which is \eqref{eq:inner-super}.

Finally, since \(\Psi^m=w\), we have
\(\Delta\Psi^{m}+\abs{x}^{\sigma}\Psi^{p}=\Delta w+\abs{x}^{\sigma}w^{q}\le0\).
Moreover \(\bigl(\tfrac{c}{2}\bigr)^{1/m}\le\Psi\le c^{1/m}\), and
\(\Psi^m=w\in H^1(B(0,1))\). Thus all the requirements of
Definition~\ref{def:supersol} are satisfied in \(B(0,1)\).
\end{proof}
\subsection{{Regular radial steady states}}\label{subsec:regular}

Throughout this subsection we assume \(N\ge3\), \(m>1\), \(\sigma\in(-2,0)\)
and \(p>\pc\), and we set
\begin{equation}\label{eq:EFconst}
q\eqdef\frac pm,\qquad a\eqdef\ac m,\qquad
\kappa_\ast\eqdef a(N-2-a)=A^{p-m},\qquad d\eqdef N-2-2a .
\end{equation}
By Proposition~\ref{prop:steady}, \(\kappa_\ast>0\). Moreover \(d>0\) if and
only if \(p>\pss\) (Remark~\ref{rem:pss}), and \(a(q-1)=\ac(p-m)=\sigma+2\).
In particular, \(\kappa_\ast=(A^m)^{q-1}\).

\begin{lemma}\label{lem:A}
For every \(b>0\), there exist \(R_b\in(0,\infty]\) and a unique function
\(V\in C([0,R_b))\cap C^2((0,R_b))\) such that \(V(0)=b\), \(V>0\) on
\([0,R_b)\), \(V(r)\to0\) as \(r\uparrow R_b\) if \(R_b<\infty\), and
\begin{equation}\label{eq:radialODE}
V''+\frac{N-1}{r}\,V'+r^{\sigma}V^{q}=0\quad\text{on }(0,R_b),
\qquad
r^{N-1}V'(r)=-\int_0^r\varrho^{N-1+\sigma}V(\varrho)^{q}\,d\varrho .
\end{equation}
In particular, \(V\) is decreasing and
\(\abs{V'(r)}\le b^{q}r^{\sigma+1}/(N+\sigma)\) on \((0,R_b)\). Moreover, the
function \(y(s)\eqdef e^{as}V(e^{s})\), \(s<\log R_b\), solves
\begin{equation}\label{eq:EF}
y''+d\,y'-\kappa_\ast\,y+y^{q}=0,
\end{equation}
and, for every \(\vartheta\in(0,1)\), there is \(s_\vartheta<\log R_b\) such
that \(y(s_\vartheta)=\vartheta^mA^m\) and \(y(s)<\vartheta^mA^m\) for all
\(s<s_\vartheta\).
\end{lemma}

\begin{proof}
Since \(\sigma>-2\), we have
\(\int_0^r t^{1-N}\int_0^t\varrho^{N-1+\sigma}\,d\varrho\,dt
=r^{\sigma+2}/\bigl((N+\sigma)(\sigma+2)\bigr)<\infty\). Hence the map
\[
V\longmapsto b-\int_0^r t^{1-N}\int_0^t\varrho^{N-1+\sigma}V(\varrho)^{q}\,d\varrho\,dt
\]
is a contraction on \(\{V\in C([0,r_\ast]):\ b/2\le V\le 3b/2\}\) for
\(r_\ast>0\) small enough, and we obtain a unique local solution of the
integral form of \eqref{eq:radialODE}. For \(r>0\), the equation is a regular
second-order ODE whose nonlinearity is locally Lipschitz on \(\{V>0\}\).
Therefore the solution extends uniquely to a maximal interval \([0,R_b)\) on
which \(V>0\). The identity in \eqref{eq:radialODE} shows that \(V'<0\) on
\((0,R_b)\), so that \(0<V\le b\). It follows that \(V(r)\to0\) as
\(r\uparrow R_b\) whenever \(R_b<\infty\), and that
\(\abs{V'(r)}\le b^{q}r^{1-N}\int_0^r\varrho^{N-1+\sigma}\,d\varrho
=b^{q}r^{\sigma+1}/(N+\sigma)\).

With \(r=e^{s}\), we have
\(e^{2s}\bigl(V''+\frac{N-1}{r}V'\bigr)=\frac{d^2}{ds^2}V(e^s)+(N-2)\frac{d}{ds}V(e^s)\).
Inserting \(V(e^{s})=e^{-as}y(s)\), we get
\(e^{(a+2)s}\bigl(V''+\frac{N-1}rV'\bigr)=y''+d\,y'-\kappa_\ast y\), while
\(e^{(a+2)s}r^{\sigma}V^{q}=e^{(\sigma+2-a(q-1))s}y^{q}=y^{q}\) because
\(a(q-1)=\sigma+2\). This proves \eqref{eq:EF}.

As \(s\to-\infty\), we have \(y(s)\to0\) and
\(y'(s)=e^{as}\bigl(aV(r)+rV'(r)\bigr)_{r=e^{s}}\), where
\(rV'(r)=O(r^{\sigma+2})\) tends to \(0\) and \(V(r)\to b\). Hence \(y'>0\) on
a half-line \((-\infty,s_\ast)\). Let \(s_1\in(-\infty,\log R_b]\) be the
supremum of all such \(s_\ast\), and set \(Y\eqdef\sup_{s<s_1}y(s)\in(0,\infty]\).
On \((-\infty,s_1)\), the function \(y\) is positive and increasing. We claim
that \(Y\ge A^m\), and we distinguish three cases.

The case \(s_1=\log R_b<\infty\) is impossible, since \(y\) would then tend to
\(0\) while being positive and increasing. If \(s_1<\log R_b\), then
\(y'(s_1)=0\) and \(y'>0\) on \((-\infty,s_1)\), so that \(y''(s_1)\le0\). By
\eqref{eq:EF}, we obtain \(\kappa_\ast y(s_1)-y(s_1)^{q}=y''(s_1)\le0\), that
is, \(y(s_1)^{q-1}\ge\kappa_\ast=(A^m)^{q-1}\), whence \(Y=y(s_1)\ge A^m\).
Assume finally that \(s_1=+\infty\) and \(Y<\infty\). Then \(z\eqdef y'>0\) is
integrable on \((0,\infty)\), and \(h\eqdef z'+dz=\kappa_\ast y-y^{q}\) tends to
\(\ell\eqdef\kappa_\ast Y-Y^{q}\) as \(s\to\infty\). Integrating \(z'=h-dz\)
over \([s,s+1]\), we get
\[
z(s+1)-z(s)=\int_s^{s+1}h-d\int_s^{s+1}z\ \longrightarrow\ \ell
\qquad(s\to\infty),
\]
because \(\int_s^{s+1}z\to0\). If \(\ell>0\), then \(z\to\infty\), which
contradicts the integrability of \(z\). If \(\ell<0\), then \(z\) eventually
becomes negative, which contradicts \(z>0\). Hence \(\ell=0\), and since
\(Y>0\), we infer that \(Y^{q-1}=\kappa_\ast\), that is, \(Y=A^m\). In all cases
\(Y\ge A^m\).

Given \(\vartheta\in(0,1)\), the level \(\vartheta^mA^m\) lies in \((0,Y)\).
The continuity and the strict monotonicity of \(y\) on \((-\infty,s_1)\)
provide a unique \(s_\vartheta<s_1\) with \(y(s_\vartheta)=\vartheta^mA^m\) and
\(y<\vartheta^mA^m\) on \((-\infty,s_\vartheta)\).
\end{proof}

\subsection{Finite time blow-up: slow decay and large amplitude}\label{subsec:blowup}
The blow-up half of the picture rests on a single computation, namely the
energy \eqref{eq:energy} of the portion of the datum sitting on the annulus
\(\{\lambda<\abs{x}<2\lambda\}\). The exponent \(\ac\) emerges as the exact
balance between the diffusive and the reactive homogeneities, which explains
why it is the second critical exponent. Throughout this section, we fix once
and for all a function
\begin{equation}\label{eq:zeta}
\zeta\in C^{\infty}_c(\RN),\qquad 0\le\zeta\le1,\qquad
\supp\zeta\subset\{1<\abs{x}<2\},\qquad \zeta\not\equiv0,
\end{equation}
and we set
\begin{equation}\label{eq:I1I2}
I_1\eqdef\int_{\RN}\bigl\lvert\nabla\zeta^{m}\bigr\rvert^{2}\dy\in(0,\infty),
\qquad
I_2\eqdef\int_{\RN}\abs{y}^{\sigma}\zeta^{m+p}\dy\in(0,\infty).
\end{equation}
Both integrals are finite, because \(\zeta^m\in C^1_c(\RN)\) (recall \(m>1\))
and \(\zeta\) is supported away from the origin. They are positive, because
\(\zeta^m\) is a nonzero compactly supported function, hence nonconstant. For
\(\lambda>0\) and \(M>0\) we write
\begin{equation}\label{eq:annulus-data}
\zeta_\lambda(x)\eqdef\zeta(x/\lambda),\qquad
v_{M,\lambda}\eqdef M\,\zeta_\lambda\in \Lpp1\cap\Lp{\infty},
\qquad \supp v_{M,\lambda}\subset\{\lambda<\abs{x}<2\lambda\},
\end{equation}
and we introduce the constant
\begin{equation}\label{eq:M0}
M_0\eqdef\left(\frac{(m+p)\,I_1}{2m\,I_2}\right)^{\frac{1}{p-m}}\in(0,\infty).
\end{equation}

\begin{lemma}\label{lem:energy-annulus}
Let \(m>1\), \(p>m\) and \(\sigma\in(-2,0)\). For all \(M,\lambda>0\), the
datum \(v_{M,\lambda}\) is admissible in Theorem~\ref{thm:IL-super}(b), that is,
\(v_{M,\lambda}\in\Lpp1\cap\Lp\infty\cap L^{m+p}_+(\RN,\abs{x}^{\sigma}\dx)\)
and \(\nabla v_{M,\lambda}^{m}\in\Lp2\). Moreover,
\begin{equation}\label{eq:Escaling}
E\bigl(v_{M,\lambda}\bigr)
=\frac{M^{2m}}{2}\,\lambda^{N-2}\,I_1
-\frac{m}{m+p}\,M^{m+p}\,\lambda^{N+\sigma}\,I_2 .
\end{equation}
In particular, for \(\lambda>0\) fixed,
\begin{equation}\label{eq:Esign}
E\bigl(v_{M,\lambda}\bigr)<0
\iff M^{\,p-m}\lambda^{\sigma+2}>M_0^{\,p-m}.
\end{equation}
Consequently, if \(M=\kappa\lambda^{-\gamma}\) for some \(\kappa>0\) and
\(\gamma>0\), then
\begin{equation}\label{eq:Escaling2}
E\bigl(v_{\kappa\lambda^{-\gamma},\lambda}\bigr)
=\lambda^{N-2-2m\gamma}\left[\frac{\kappa^{2m}}{2}I_1
-\frac{m}{m+p}\,\kappa^{m+p}\,\lambda^{\,\sigma+2-\gamma(p-m)}\,I_2\right],
\end{equation}
and the exponent \(\sigma+2-\gamma(p-m)=(p-m)(\ac-\gamma)\) has the sign of
\(\ac-\gamma\).
\end{lemma}
\begin{proof}
Since \(0\le\zeta\le1\), we have \(v_{M,\lambda}\in L^\infty(\RN)\). Its support
is contained in \(\{\lambda<\abs{x}<2\lambda\}\), which is compact and bounded
away from the origin, so that \(\abs{x}^{\sigma}\) is bounded on
\(\supp v_{M,\lambda}\). Hence
\(v_{M,\lambda}\in L^1(\RN)\cap L^{m+p}(\RN,\abs{x}^{\sigma}\dx)\). Moreover,
\(\nabla v_{M,\lambda}^{m}=M^{m}\lambda^{-1}(\nabla\zeta^{m})(\cdot/\lambda)\)
is continuous with compact support, so that
\(\nabla v_{M,\lambda}^{m}\in L^2(\RN)\). This proves the admissibility, and all
the integrals defining \(E\) are finite.

The substitution \(x=\lambda y\) gives
\[
\int_{\RN}\bigl\lvert\nabla v_{M,\lambda}^{m}\bigr\rvert^{2}\dx
=M^{2m}\lambda^{-2}\lambda^{N}\int_{\RN}\bigl\lvert\nabla\zeta^{m}\bigr\rvert^{2}\dy
=M^{2m}\lambda^{N-2}I_1
\]
and
\[
\int_{\RN}\abs{x}^{\sigma}v_{M,\lambda}^{\,m+p}\dx
=M^{m+p}\lambda^{\sigma}\lambda^{N}\int_{\RN}\abs{y}^{\sigma}\zeta^{m+p}\dy
=M^{m+p}\lambda^{N+\sigma}I_2 .
\]
Inserting these identities into the definition \eqref{eq:energy} of \(E\), we
obtain \eqref{eq:Escaling}. Dividing the inequality \(E(v_{M,\lambda})<0\) by
\(\frac{m}{m+p}M^{2m}\lambda^{N-2}I_2>0\), we get
\[
E\bigl(v_{M,\lambda}\bigr)<0
\iff M^{p-m}\lambda^{\sigma+2}
>\frac{(m+p)I_1}{2mI_2}=M_0^{\,p-m},
\]
which is \eqref{eq:Esign}. For \eqref{eq:Escaling2}, we insert
\(M=\kappa\lambda^{-\gamma}\) into \eqref{eq:Escaling} and factor out
\(\lambda^{N-2-2m\gamma}\). The exponent of \(\lambda\) left in the second term
is
\[
\bigl(N+\sigma-(m+p)\gamma\bigr)-\bigl(N-2-2m\gamma\bigr)
=\sigma+2-\gamma(p-m)=(p-m)(\ac-\gamma),
\]
by \eqref{eq:ac}. Consequently, as \(\lambda\to\infty\), the bracket in
\eqref{eq:Escaling2} tends to \(-\infty\) when \(\gamma<\ac\), it is independent
of \(\lambda\) when \(\gamma=\ac\), and it tends to
\(\tfrac{\kappa^{2m}}{2}I_1>0\) when \(\gamma>\ac\). In the last case the
annulus energy is positive for \(\lambda\) large, and no blow-up can be
extracted from far annuli. This scale-invariant balance identifies \(\ac\) as
the second critical exponent.
\end{proof}

%===================================================
\section{{The singular steady state and global existence at the critical decay}}\label{sec:proofs}

In this section we prove the results of Section~\ref{subsec:mainresults}
which concern the singular steady state and global existence, namely
Propositions~\ref{prop:steady} and \ref{rem:LequalA},
Theorems~\ref{thm:trap}, \ref{thm:barrier}, \ref{thm:global} and
\ref{thm:conj91}, and Corollaries~\ref{cor:sandwich} and \ref{cor:belowS}.

\subsection{{The singular steady state and trapping}}

\begin{proof}[Proof of Proposition~\ref{prop:steady}]
Set \(a\eqdef\ac m>0\). By \eqref{eq:ac}, we have
\(\ac p-\sigma=a+\ac(p-m)-\sigma=a+2\). Hence, on \(\RN\setminus\{0\}\),
\eqref{eq:laplacian} with \(b=a\) gives
\[
\Delta S^m+\abs{x}^{\sigma}S^{p}
=A^m a(a+2-N)\abs{x}^{-a-2}+A^{p}\abs{x}^{-a-2}
=A^m\abs{x}^{-a-2}\Bigl[A^{p-m}-a(N-2-a)\Bigr].
\]
Therefore \eqref{eq:stationary} holds for some \(A>0\) if and only if
\(a(N-2-a)>0\), that is, \(a<N-2\). In that case \(A\) is uniquely given by
\eqref{eq:A}. Since \(a=m(\sigma+2)/(p-m)\), the condition \(a<N-2\) is
equivalent to \(m(\sigma+2)<(N-2)(p-m)\), that is, to
\(p>m+m(\sigma+2)/(N-2)=\pc\). This proves the equivalence of (i) and (ii).

Assume now that \(p>\pc\). Since \(m>1\), we have \(\ac<a<N-2<N\), so that
\(S\) and \(S^m=A^m\abs{x}^{-a}\) belong to \(L^1_{\mathrm{loc}}(\RN)\).
Moreover \(\abs{x}^{\sigma}S^{p}=A^{p}\abs{x}^{-a-2}\) with \(a+2<N\), hence
\(\abs{x}^{\sigma}S^{p}\in L^1_{\mathrm{loc}}(\RN)\). Lemma~\ref{lem:dirac}(i)
with \(b=a\in(0,N-2)\) shows that \(\Delta S^m=A^ma(a+2-N)\abs{x}^{-a-2}\) in
\(\Dprime{\RN}\), and we conclude that \eqref{eq:stationary} holds in
\(\Dprime{\RN}\).

Finally, \(-\sigma(m-1)/2\in(0,m-1)\) gives \(\pG\in(1,m)\), and \(\pF>m\) is
obvious. Next, \(\pc-m=m(\sigma+2)/(N-2)>(\sigma+2)/N=\pF-m\), because
\(m>1>(N-2)/N\). {Lastly},
\(\pss-\pc=m(N+2\sigma+2-\sigma-N)/(N-2)=m(\sigma+2)/(N-2)=\pc-m>0\).
\end{proof}

We now turn to the trapping theorem. It follows directly from the comparison
lemma and from the properties of the singular steady state.

\begin{proof}[Proof of Theorem~\ref{thm:trap}]
By Proposition~\ref{prop:steady}, \(S\) solves \eqref{eq:stationary}
classically on \(\mathbb{R}^N\setminus\{0\}\). We check that \(\vartheta S\)
satisfies (H1) of Lemma~\ref{lem:comparison}, and we begin with the regularity
requirements of Definition~\ref{def:supersol} on \(\Omega=\RN\setminus\{0\}\).
Since \(p>p_c\), we have \(0<\alpha_c m<N-2\), so that \(\alpha_c<N\) and
\(\vartheta S\in L^1_{\mathrm{loc}}(\mathbb{R}^N)\). The function
\((\vartheta S)^m=\vartheta^m A^m\abs{x}^{-\alpha_c m}\) is smooth on
\(\mathbb{R}^N\setminus\{0\}\), so that
\((\vartheta S)^m\in H^1_{\mathrm{loc}}(\mathbb{R}^N\setminus\{0\})\). By
\eqref{eq:ac}, we have
\(\abs{x}^{\sigma}(\vartheta S)^{p}=\vartheta^p A^p\abs{x}^{-\alpha_c m-2}\),
which is locally integrable on \(\RN\) since \(\alpha_c m+2<N\). Moreover,
\(\vartheta S(x)\le\vartheta A\delta^{-\alpha_c}\) for \(\abs{x}\ge\delta\), so
that \(\vartheta S\in L^\infty(\{\abs{x}\ge\delta\})\) for every \(\delta>0\).

Next, for every \(0\le\varphi\in C_c^\infty(\mathbb{R}^N\setminus\{0\})\), we
have
\begin{align*}
&\int_{\mathbb{R}^N\setminus\{0\}}
\Bigl(\nabla(\vartheta S)^m\cdot\nabla\varphi
-\abs{x}^{\sigma}(\vartheta S)^{p}\varphi\Bigr)\dx\\
&\qquad=\vartheta^m
\underbrace{\int_{\mathbb{R}^N\setminus\{0\}}
\Bigl(\nabla S^m\cdot\nabla\varphi-\abs{x}^{\sigma}S^{p}\varphi\Bigr)\dx}_{=\,0}
+(\vartheta^m-\vartheta^p)
\int_{\mathbb{R}^N\setminus\{0\}}\abs{x}^{\sigma}S^{p}\varphi\dx\ \ge\ 0 .
\end{align*}
The first term vanishes because \(S\) is a solution on
\(\RN\setminus\{0\}\), and the second one is nonnegative because
\(\vartheta^p\le\vartheta^m\) and \(\abs{x}^{\sigma}S^{p}\varphi\ge0\). This is
Lemma~\ref{lem:submultiple} in the present situation, and it shows that
\(\vartheta S\) is a stationary weak supersolution of \eqref{eq:CP} in
\(\RN\setminus\{0\}\). Hence (H1) holds.

Since \(\alpha_c>0\), we have \(\vartheta S(x)\to+\infty\) as \(x\to0\), which
is the first alternative of (H2). Finally, (H3) is exactly \eqref{eq:below}.
Lemma~\ref{lem:comparison} yields \eqref{eq:trapped}, and the uniform bound
away from the origin follows at once.
\end{proof}

\subsection{{The explicit barrier and global existence}}

We now prove Theorem~\ref{thm:barrier}. The proof is divided into three steps.
We first locate the matching radius $\rho$, and we establish the crossing
property \eqref{eq:crossing} between the inner profile $w$ and the outer
profile $\Sigma=\mu^m S^m$. We then check the continuity and the regularity of
the resulting barrier $\overline{u}$, as well as its membership in
$L^{r+1}(\mathbb{R}^N)$ for $r>r_0$. Finally, we verify the supersolution
property. We treat separately the absolutely continuous part of the
distributional Laplacian and the singular measure supported on the sphere
$\{|x|=\rho\}$. The key point is that the jump of the normal derivative is
nonpositive, so that the singular part reinforces the supersolution
inequality.

\begin{proof}[Proof of Theorem~\ref{thm:barrier}]
Write $\Sigma(\varrho)\eqdef\mu^{m}A^{m}\varrho^{-\alpha_c m}$ for the radial
profile of $(\mu S)^{m}$, and recall that
$w(\varrho)=c(1-\tfrac12\varrho^{\sigma+2})\in[c/2,c]$ for $\varrho\in(0,1]$.

\smallskip
\noindent\emph{Step 1: the matching radius.} By the definition of $\mu$ in
\eqref{eq:mu}, for $\varrho\in[\tfrac12,1]$ we have
\[
\Sigma(\varrho)\le \mu^{m}A^{m}2^{\alpha_c m}=\bigl(\mu A2^{\alpha_c}\bigr)^{m}\le\frac{c}{2}
\le w(\varrho).
\]
Hence the set which defines $\rho$ contains $[\tfrac12,1]$, and therefore
$\rho\le\tfrac12$. On the other hand, $\Sigma(\varrho)\to+\infty$ as
$\varrho\downarrow0$, while $w\le c$, so that $\rho>0$. The function
$\Sigma-w$ is continuous on $(0,1]$. It is positive on $(0,\rho)$ by the
minimality of $\rho$, and it is nonpositive at points of the defining set
arbitrarily close to $\rho$. We deduce that
\begin{equation}\label{eq:crossing}
\Sigma(\varrho)>w(\varrho)\quad\text{for }0<\varrho<\rho,
\qquad \Sigma(\rho)=w(\rho).
\end{equation}

\smallskip
\noindent\emph{Step 2: continuity, boundedness and integrability.} By
\eqref{eq:crossing}, both branches of \(\overline{u}\) coincide on
\(\{\abs{x}=\rho\}\), so that \(\overline{u}\) is continuous on \(\RN\). Since
\(\Psi\) and \(\mu S\) are both radially decreasing, \(\overline{u}\) is radially
nonincreasing, and
\(\|\overline{u}\|_\infty=\Psi(0)=c^{1/m}\). Because \(\rho\le\tfrac12\), the
identity \(\overline{u}=\mu A\abs{x}^{-\alpha_c}\) holds on
\(\{\abs{x}\ge\tfrac12\}\), which is \eqref{eq:barrier}. Moreover,
\(\overline{u}\in L^{r+1}(\mathbb{R}^N)\) for \(r>r_0\), because \(\overline{u}\)
is bounded and \(\alpha_c(r+1)>\alpha_c(r_0+1)=N\) by
Remark~\ref{rem:duality}. As for the local regularity, set
\[
W\eqdef\overline{u}^{\,m},
\qquad\text{so that}\qquad
W=w\text{ on }B(0,\rho),
\quad W=\mu^{m}A^{m}\abs{x}^{-\alpha_c m}\text{ outside}.
\]
Both branches are smooth away from the origin and from
\(\{\abs{x}=\rho\}\), and \(W\) is continuous across \(\{\abs{x}=\rho\}\). By
Proposition~\ref{prop:inner}, \(w\in H^1(B(0,1))\), and the outer branch is
smooth on \(\{\abs{x}\ge\rho\}\). Hence \(W\in H^1_{\mathrm{loc}}(\mathbb{R}^N)\).
Moreover \(\abs{x}^{\sigma}\overline{u}^{\,p}\in L^1_{\mathrm{loc}}(\mathbb{R}^N)\),
because \(\abs{x}^{\sigma}\overline{u}^{\,p}\le c^{p/m}\abs{x}^{\sigma}\) near the
origin, and \(\sigma>-N\).

\smallskip
\noindent\emph{Step 3: the supersolution property.} Let
\(\mathbb S_\rho\eqdef\{\abs{x}=\rho\}\). Since \(W\) is Lipschitz across
\(\mathbb S_\rho\) and smooth on either side, its distributional Laplacian is
\begin{equation}\label{eq:jump}
\Delta W=\bigl(\Delta W\bigr)_{\mathrm{ac}}
+\Bigl[\partial_\varrho \Sigma(\rho^{+})-\partial_\varrho w(\rho^{-})\Bigr]\,
\mathcal{H}^{N-1}\lfloor_{\mathbb S_\rho}.
\end{equation}
Here the absolutely continuous part is computed branchwise, and the coefficient
of the singular part is the jump of the outward normal derivative. By
\eqref{eq:crossing}, the function \(\Sigma-w\) is positive on \((0,\rho)\) and
vanishes at \(\rho\). Both branches being smooth at \(\varrho=\rho>0\), this
forces
\[
\partial_\varrho\bigl(\Sigma-w\bigr)(\rho)\le0,
\qquad\text{that is,}\qquad
\partial_\varrho \Sigma(\rho)-\partial_\varrho w(\rho)\le0,
\]
so that the singular part in \eqref{eq:jump} is a nonpositive measure.

As for the absolutely continuous part, on \(\{\abs{x}<\rho\}\) we have, by
Proposition~\ref{prop:inner} and \(\rho\le1\),
\[
\Delta W+\abs{x}^{\sigma}\overline{u}^{\,p}
=\Delta w+\abs{x}^{\sigma}w^{q}\le0,
\]
and Lemma~\ref{lem:dirac}(ii) guarantees that no Dirac mass is created at the
origin. On \(\{\abs{x}>\rho\}\), using \(\mu\le1\), \(p>m\) and
Proposition~\ref{prop:steady}, we get
\[
\Delta(\mu S)^{m}+\abs{x}^{\sigma}(\mu S)^{p}
=\mu^{m}\Delta S^{m}+\mu^{p}\abs{x}^{\sigma}S^{p}
=\bigl(\mu^{p}-\mu^{m}\bigr)\abs{x}^{\sigma}S^{p}\le0.
\]
Adding up the two contributions, we conclude that
\[
\Delta \overline{u}^{\,m}+\abs{x}^{\sigma}\overline{u}^{\,p}\le0
\qquad\text{in }\mathcal{D}'(\mathbb{R}^N),
\]
which, together with Step~2, shows that \(\ol u\) is a stationary weak
supersolution in \(\RN\).
\end{proof}
\begin{remark}\rm \label{rem:explicit}
All the constants produced by the construction are explicit. Indeed,
$c^{1/m}=(D/2)^{1/(p-m)}$ since $q-1=(p-m)/m$, so that
\begin{equation}\label{eq:explicit}
\|\overline{u}\|_{\infty}=\Bigl(\frac{(\sigma+2)(N+\sigma)}{2}\Bigr)^{\frac{1}{p-m}},
\qquad
A=\bigl[\alpha_c m(N-2-\alpha_c m)\bigr]^{\frac{1}{p-m}},
\end{equation}
and consequently
\[
\mu=\min\left\{1,\ 2^{-\alpha_c-\frac1m}
\left(\frac{(\sigma+2)(N+\sigma)}{2\,\alpha_c m\,(N-2-\alpha_c m)}\right)^{\frac{1}{p-m}}\right\}.
\]
Both quantities in \eqref{eq:explicit} are of the form
$(\,\cdot\,)^{1/(p-m)}$. This is the natural homogeneity of amplitudes of
stationary profiles. Indeed, in the variable $v=u^m$, the amplitudes $c$ and
$A^m$ are determined by algebraic relations of the form $c^{q-1}=\mathrm{const}$
(see Proposition~\ref{prop:inner} and \eqref{eq:A}), and taking the $m$-th root
gives the exponent $1/(m(q-1))=1/(p-m)$. For instance, for $N=3$, $m=2$,
$\sigma=-1$ and $p=5$, one finds $\ac=1/3$, $A\approx0.606$ and
$\mu\approx0.927$.
\end{remark}
We now proceed to the proof of Theorem~\ref{thm:global}. The argument combines
the comparison lemma with the bounded barrier constructed in
Theorem~\ref{thm:barrier}. We first show that every sub-multiple
$\vartheta\overline{u}$ is an admissible supersolution for
Lemma~\ref{lem:comparison}. We then deduce the uniform bound
\eqref{eq:global-bound} and the global existence of the solution. Finally, we
verify that the data $L\phi$, with $0<L\le\mu A$, are dominated by the
barrier, which yields \eqref{eq:global-bound-L}.

\begin{proof}[Proof of Theorem~\ref{thm:global}]
By Theorem~\ref{thm:barrier} and Lemma~\ref{lem:submultiple},
$\vartheta\overline{u}$ is a stationary weak supersolution of \eqref{eq:CP}
in $\RN$. Explicitly, for every $0\le\varphi\in C_c^\infty(\mathbb{R}^N)$, since
$p>m$ and $\vartheta\in(0,1]$ imply $\vartheta^{p}\le\vartheta^{m}$, we have
\begin{align*}
\int_{\mathbb{R}^N}\Bigl(\nabla(\vartheta\overline{u})^{m}\cdot\nabla\varphi
-|x|^{\sigma}(\vartheta\overline{u})^{p}\varphi\Bigr)\,dx
&=\vartheta^{m}\underbrace{\int_{\mathbb{R}^N}\Bigl(\nabla\overline{u}^{m}\cdot\nabla\varphi
-|x|^{\sigma}\overline{u}^{p}\varphi\Bigr)\,dx}_{\ge0}\\
&\quad+(\vartheta^{m}-\vartheta^{p})\int_{\mathbb{R}^N}|x|^{\sigma}\overline{u}^{p}\varphi\,dx\ \ge\ 0 .
\end{align*}
Moreover, $\vartheta\overline{u}$ is bounded and positive, and
$(\vartheta\overline{u})^{m}\in H^1(B(0,1))$, since
$\overline{u}^{m}\in H^1_{\mathrm{loc}}(\mathbb{R}^N)$ by
Theorem~\ref{thm:barrier}. Hence $\vartheta\overline{u}$ satisfies (H1) and the
second alternative of (H2) of Lemma~\ref{lem:comparison}, while (H3) is
\eqref{eq:below-barrier}.

Lemma~\ref{lem:comparison} gives $u(t,\cdot)\le\vartheta\overline{u}$ on
$[0,T_{\max}(u_0))$. We infer that
$\|u(t)\|_{\infty}\le\vartheta\|\overline{u}\|_{\infty}=\vartheta c^{1/m}$
for every $t<T_{\max}(u_0)$. If $T_{\max}(u_0)$ were finite,
Assumption~\ref{ass:IL}(A) would force
$\limsup_{t\uparrow T_{\max}}\|u(t)\|_{\infty}=\infty$, which is a
contradiction. Therefore $T_{\max}(u_0)=\infty$, and \eqref{eq:global-bound}
holds.

It remains to check that $L\phi\le\vartheta\overline{u}$ with
$\vartheta\eqdef L/(\mu A)\in(0,1]$ when $0<L\le\mu A$, that is, that
$\mu A\,\phi\le\overline{u}$. We distinguish three regions. For $|x|\ge1$, we
have $\phi(x)=|x|^{-\alpha_c}$ and, by \eqref{eq:barrier},
$\overline{u}(x)=\mu A|x|^{-\alpha_c}=\mu A\phi(x)$. For
$\tfrac12\le|x|\le1$, we have $\phi(x)=1$ and
$\overline{u}(x)=\mu A|x|^{-\alpha_c}\ge\mu A$. For $|x|\le\tfrac12$, we have
$\phi(x)=1$, and, $\overline{u}$ being radially nonincreasing,
$\overline{u}(x)\ge\overline{u}\bigl(\tfrac12 e\bigr)=\mu A\,2^{\alpha_c}\ge\mu A$
for any unit vector $e$. Hence $\mu A\phi\le\overline{u}$, and
\eqref{eq:global-bound-L} follows from \eqref{eq:global-bound}. Note that
every $u_0$ with $0\le u_0\le L\phi$ belongs to $\mathcal{A}_r$ for every
$r>r_0$, since $\phi$ is bounded and $\alpha_c(r+1)>N$
(Remark~\ref{rem:duality}). Finally, for $\gamma\ge\alpha_c$ we have
$\phi_{\gamma}\le\phi$, and the last assertion follows from the previous
one.
\end{proof}

\subsection{{\texorpdfstring{{The improved lower bound}}{The improved lower bound}}}

\begin{remark}\rm \label{rem:strategy}
A natural way to improve the lower bound \(L_\ast\ge\mu A\) would be to exploit
the strict inequality \(u<S\) near the origin in order to apply an
\(L^\infty\) smoothing estimate. However, by Remark~\ref{rem:origin}, the
source term \(\abs{x}^{\sigma}u^{p}\) belongs to \(L^q_{\mathrm{loc}}\) only for
\(q<N/2\), below the threshold \(q>N/2\), and the gap \(u<S\) does not improve
this integrability exponent. We therefore follow a different route. We replace
the explicit inner profile \(w\) of Proposition~\ref{prop:inner} by the
\emph{regular radial steady states} of \eqref{eq:CP}. In Emden--Fowler
variables, these steady states reach every level below the singular amplitude
(Lemma~\ref{lem:A}). Gluing them to \(\vartheta S\), \(\vartheta<1\), at their
first crossing produces bounded barriers with tail \(\vartheta A\)
(Theorem~\ref{thm:conj91}).
\end{remark}

\begin{proof}[Proof of Theorem~\ref{thm:conj91}]
Fix \(\vartheta\in(0,1)\), let \(V\) be given by Lemma~\ref{lem:A} with \(b=1\),
and set \(\rho_1\eqdef e^{s_\vartheta}<R_1\). We define
\[
W(r)\eqdef
\begin{cases}
V(r), & 0\le r<\rho_1,\\
\Sigma(r)\eqdef\vartheta^mA^m r^{-a}, & r\ge\rho_1,
\end{cases}
\qquad\qquad
\mathcal U(x)\eqdef W(\abs{x})^{1/m}.
\]
Since \(V(r)=r^{-a}y(\log r)\), Lemma~\ref{lem:A} gives
\(\Sigma(r)-V(r)=r^{-a}\bigl(\vartheta^mA^m-y(\log r)\bigr)>0\) for
\(0<r<\rho_1\), and \(\Sigma(\rho_1)=V(\rho_1)\). Hence \(W\) is continuous,
positive, nonincreasing and bounded by \(V(0)=1\). Moreover, the function
\(\Sigma-V\) is smooth near \(\rho_1>0\), positive on \((0,\rho_1)\), and it
vanishes at \(\rho_1\). We deduce that \((\Sigma-V)'(\rho_1)\le0\).

On \(B(0,\rho_1)\), the function \(W=V\) solves
\(\Delta W+\abs{x}^{\sigma}W^{q}=0\) in \(\Dprime{B(0,\rho_1)}\), without a
Dirac mass at the origin. Indeed, \(V\) is continuous at \(0\) and
\(r^{N-1}V'(r)=O(r^{N+\sigma})\to0\) as \(r\to0\), so that the boundary terms in
Green's formula \eqref{eq:green} on \(\{\abs{x}>\eps\}\) vanish as
\(\eps\downarrow0\), exactly as in Lemma~\ref{lem:dirac}. On
\(\{\abs{x}>\rho_1\}\), the function \(W=(\vartheta S)^m\) satisfies
\(\Delta W+\abs{x}^{\sigma}W^{q}=(\vartheta^p-\vartheta^m)\abs{x}^{\sigma}S^{p}\le0\)
by Proposition~\ref{prop:steady}. Since \(W\) is Lipschitz across
\(\{\abs{x}=\rho_1\}\) and smooth on either side, its distributional Laplacian
is given by the analogue of \eqref{eq:jump}. Its singular part
\(\bigl(\Sigma'(\rho_1)-V'(\rho_1)\bigr)\Hmeas\lfloor_{\{\abs{x}=\rho_1\}}\) is a
nonpositive measure. Therefore
\(\Delta \mathcal U^m+\abs{x}^{\sigma}\mathcal U^{p}\le0\) in \(\Dprime{\RN}\).
In addition, \(\mathcal U^m=W\in H^1_{\mathrm{loc}}(\RN)\), because
\(\abs{V'(r)}\le r^{\sigma+1}/(N+\sigma)\) and \(2\sigma+2>-N\). Finally,
\(\abs{x}^{\sigma}\mathcal U^{p}\in L^1_{\mathrm{loc}}(\RN)\), because
\(\mathcal U\) is bounded and \(\sigma>-N\). Thus \(\mathcal U\) is a bounded
stationary weak supersolution of \eqref{eq:CP} in \(\RN\).

By the scaling of Remark~\ref{rem:scaling}, the function
\(\ol u_\vartheta(x)\eqdef\rho_1^{\ac}\,\mathcal U(\rho_1x)\) is again a
stationary weak supersolution in \(\RN\). It equals
\(\vartheta A\abs{x}^{-\ac}\) for \(\abs{x}\ge1\). Being radially nonincreasing,
it satisfies \(\ol u_\vartheta(x)\ge\ol u_\vartheta(e_1)=\vartheta A\) for
\(\abs{x}\le1\), where \(e_1\) is any unit vector. This is
\eqref{eq:barrier-theta}. The function \(\ol u_\vartheta\) is bounded,
\(\ol u_\vartheta^{\,m}\in H^1(B(0,1))\), and it is a stationary weak
supersolution in \(\RN\), hence also in \(\RN\setminus\{0\}\). Consequently,
(H1), the second alternative of (H2), and (H3) of Lemma~\ref{lem:comparison}
hold whenever \(u_0\le\ol u_\vartheta\), and Lemma~\ref{lem:comparison} gives
\(u(t)\le\ol u_\vartheta\) on \([0,\Tmax(u_0))\). The blow-up alternative of
Assumption~\ref{ass:IL}(A) then yields \(\Tmax(u_0)=\infty\), and
\eqref{B-Infty} follows.

Finally, let \(0<L<A\) and choose \(\vartheta\in(L/A,1)\). For \(\abs{x}\ge1\)
we have \(L\phi_{\ac}(x)=L\abs{x}^{-\ac}\le\vartheta A\abs{x}^{-\ac}
=\ol u_\vartheta(x)\), and for \(\abs{x}\le1\) we have
\(L\phi_{\ac}(x)=L\le\vartheta A\le\ol u_\vartheta(x)\). Since
\(L\phi_{\ac}\in\mathcal A_r\) for every \(r>r_0\), the first part of the
theorem applies.
\end{proof}

\begin{proof}[Proof of Corollary~\ref{cor:sandwich}]
By Theorem~\ref{thm:conj91}, we have \((0,A)\subset\mathcal G_{\ac}\), hence
\(L_\ast(\ac)\ge A\). The upper bound is Theorem~\ref{thm:threshold}(v), which
holds for every admissible \(\zeta\). Finally, if \(\gamma\ge\ac\), then
\(L_\ast(\gamma)\ge L_\ast(\ac)\ge A\) by the monotonicity of \(L_\ast\) in
Theorem~\ref{thm:threshold}(ii).
\end{proof}

\begin{remark}\rm\label{rem:supersede}
Since \(\mu\le1\), Corollary~\ref{cor:sandwich} supersedes the lower bound of
Theorem~\ref{thm:threshold}(v). Theorems~\ref{thm:barrier} and
\ref{thm:global} remain useful, because their barrier is completely explicit.
\end{remark}

\begin{figure}[ht]
\centering
\begin{tikzpicture}
\begin{axis}[
  width=12.5cm, height=7.2cm,
  xmin=-5, xmax=10, ymin=0, ymax=2.15,
  xlabel={$s=\log r$}, ylabel={$y(s)/A^m$}]
\addplot[very thick,red!75!black] coordinates {(-5.0000,0.0969) (-4.9445,0.1005) (-4.8891,0.1043) (-4.8336,0.1082) (-4.7782,0.1123) (-4.7227,0.1165) (-4.6673,0.1208) (-4.6118,0.1253) (-4.5564,0.1300) (-4.5009,0.1349) (-4.4455,0.1399) (-4.3900,0.1451) (-4.3346,0.1505) (-4.2791,0.1561) (-4.2237,0.1620) (-4.1682,0.1680) (-4.1128,0.1742) (-4.0573,0.1807) (-4.0019,0.1874) (-3.9464,0.1944) (-3.8910,0.2016) (-3.8355,0.2091) (-3.7800,0.2168) (-3.7246,0.2248) (-3.6691,0.2332) (-3.6137,0.2418) (-3.5582,0.2507) (-3.5028,0.2599) (-3.4473,0.2695) (-3.3919,0.2794) (-3.3364,0.2896) (-3.2810,0.3002) (-3.2255,0.3112) (-3.1701,0.3225) (-3.1146,0.3343) (-3.0592,0.3465) (-3.0037,0.3590) (-2.9483,0.3720) (-2.8928,0.3855) (-2.8374,0.3994) (-2.7819,0.4137) (-2.7264,0.4286) (-2.6710,0.4439) (-2.6155,0.4597) (-2.5601,0.4761) (-2.5046,0.4929) (-2.4492,0.5104) (-2.3937,0.5283) (-2.3383,0.5468) (-2.2828,0.5659) (-2.2274,0.5856) (-2.1719,0.6058) (-2.1165,0.6267) (-2.0610,0.6481) (-2.0056,0.6702) (-1.9501,0.6929) (-1.8947,0.7162) (-1.8392,0.7402) (-1.7838,0.7648) (-1.7283,0.7900) (-1.6729,0.8158) (-1.6174,0.8423) (-1.5619,0.8694) (-1.5065,0.8971) (-1.4510,0.9255) (-1.3956,0.9544) (-1.3401,0.9839) (-1.2847,1.0140) (-1.2292,1.0446) (-1.1738,1.0758) (-1.1183,1.1074) (-1.0629,1.1395) (-1.0074,1.1720) (-0.9520,1.2050) (-0.8965,1.2382) (-0.8411,1.2718) (-0.7856,1.3056) (-0.7302,1.3396) (-0.6747,1.3738) (-0.6193,1.4079) (-0.5638,1.4421) (-0.5083,1.4762) (-0.4529,1.5101) (-0.3974,1.5438) (-0.3420,1.5771) (-0.2865,1.6099) (-0.2311,1.6422) (-0.1756,1.6739) (-0.1202,1.7048) (-0.0647,1.7349) (-0.0093,1.7639) (0.0462,1.7919) (0.1016,1.8187) (0.1571,1.8441) (0.2125,1.8681) (0.2680,1.8906) (0.3234,1.9114) (0.3789,1.9304) (0.4343,1.9475) (0.4898,1.9627) (0.5452,1.9757) (0.6007,1.9866) (0.6562,1.9952) (0.7116,2.0015) (0.7671,2.0054) (0.8225,2.0068) (0.8780,2.0056) (0.9334,2.0019) (0.9889,1.9956) (1.0443,1.9867) (1.0998,1.9752) (1.1552,1.9610) (1.2107,1.9442) (1.2661,1.9248) (1.3216,1.9028) (1.3770,1.8783) (1.4325,1.8512) (1.4879,1.8218) (1.5434,1.7899) (1.5988,1.7558) (1.6543,1.7193) (1.7097,1.6808) (1.7652,1.6401) (1.8207,1.5974) (1.8761,1.5528) (1.9316,1.5064) (1.9870,1.4582) (2.0425,1.4084) (2.0979,1.3570) (2.1534,1.3041) (2.2088,1.2497) (2.2643,1.1940) (2.3197,1.1370) (2.3752,1.0788) (2.4306,1.0194) (2.4861,0.9588) (2.5415,0.8972) (2.5970,0.8345) (2.6524,0.7708) (2.7079,0.7061) (2.7633,0.6403) (2.8188,0.5736) (2.8743,0.5058) (2.9297,0.4370) (2.9852,0.3671) (3.0406,0.2961) (3.0961,0.2240) (3.1515,0.1506) (3.2070,0.0760) (3.2624,0.0000)};
\addplot[very thick,blue!70!black] coordinates {(-5.0000,0.3108) (-4.8993,0.3213) (-4.7987,0.3321) (-4.6980,0.3433) (-4.5973,0.3548) (-4.4966,0.3668) (-4.3960,0.3790) (-4.2953,0.3917) (-4.1946,0.4048) (-4.0940,0.4183) (-3.9933,0.4322) (-3.8926,0.4465) (-3.7919,0.4613) (-3.6913,0.4764) (-3.5906,0.4921) (-3.4899,0.5081) (-3.3893,0.5246) (-3.2886,0.5416) (-3.1879,0.5590) (-3.0872,0.5769) (-2.9866,0.5952) (-2.8859,0.6139) (-2.7852,0.6331) (-2.6846,0.6527) (-2.5839,0.6726) (-2.4832,0.6930) (-2.3826,0.7137) (-2.2819,0.7347) (-2.1812,0.7560) (-2.0805,0.7776) (-1.9799,0.7994) (-1.8792,0.8213) (-1.7785,0.8433) (-1.6779,0.8653) (-1.5772,0.8872) (-1.4765,0.9091) (-1.3758,0.9307) (-1.2752,0.9520) (-1.1745,0.9729) (-1.0738,0.9932) (-0.9732,1.0130) (-0.8725,1.0320) (-0.7718,1.0501) (-0.6711,1.0673) (-0.5705,1.0834) (-0.4698,1.0983) (-0.3691,1.1119) (-0.2685,1.1242) (-0.1678,1.1350) (-0.0671,1.1443) (0.0336,1.1520) (0.1342,1.1581) (0.2349,1.1626) (0.3356,1.1655) (0.4362,1.1668) (0.5369,1.1665) (0.6376,1.1647) (0.7383,1.1615) (0.8389,1.1569) (0.9396,1.1511) (1.0403,1.1441) (1.1409,1.1361) (1.2416,1.1272) (1.3423,1.1175) (1.4430,1.1072) (1.5436,1.0963) (1.6443,1.0850) (1.7450,1.0734) (1.8456,1.0616) (1.9463,1.0498) (2.0470,1.0380) (2.1477,1.0263) (2.2483,1.0148) (2.3490,1.0035) (2.4497,0.9927) (2.5503,0.9822) (2.6510,0.9722) (2.7517,0.9627) (2.8523,0.9537) (2.9530,0.9453) (3.0537,0.9376) (3.1544,0.9304) (3.2550,0.9239) (3.3557,0.9180) (3.4564,0.9128) (3.5570,0.9082) (3.6577,0.9043) (3.7584,0.9010) (3.8591,0.8984) (3.9597,0.8964) (4.0604,0.8951) (4.1611,0.8943) (4.2617,0.8941) (4.3624,0.8945) (4.4631,0.8955) (4.5638,0.8970) (4.6644,0.8990) (4.7651,0.9014) (4.8658,0.9043) (4.9664,0.9077) (5.0671,0.9114) (5.1678,0.9155) (5.2685,0.9199) (5.3691,0.9247) (5.4698,0.9297) (5.5705,0.9350) (5.6711,0.9404) (5.7718,0.9460) (5.8725,0.9518) (5.9732,0.9577) (6.0738,0.9636) (6.1745,0.9696) (6.2752,0.9755) (6.3758,0.9814) (6.4765,0.9872) (6.5772,0.9930) (6.6779,0.9985) (6.7785,1.0040) (6.8792,1.0092) (6.9799,1.0141) (7.0805,1.0188) (7.1812,1.0232) (7.2819,1.0274) (7.3826,1.0311) (7.4832,1.0346) (7.5839,1.0377) (7.6846,1.0404) (7.7852,1.0427) (7.8859,1.0447) (7.9866,1.0462) (8.0872,1.0474) (8.1879,1.0482) (8.2886,1.0486) (8.3893,1.0486) (8.4899,1.0483) (8.5906,1.0476) (8.6913,1.0467) (8.7919,1.0454) (8.8926,1.0438) (8.9933,1.0419) (9.0940,1.0399) (9.1946,1.0376) (9.2953,1.0351) (9.3960,1.0324) (9.4966,1.0296) (9.5973,1.0267) (9.6980,1.0237) (9.7987,1.0206) (9.8993,1.0175) (10.0000,1.0143)};
\addplot[very thick,green!45!black,dashed] coordinates {(-5.0000,0.5674) (-4.8993,0.5738) (-4.7987,0.5802) (-4.6980,0.5867) (-4.5973,0.5932) (-4.4966,0.5998) (-4.3960,0.6065) (-4.2953,0.6133) (-4.1946,0.6201) (-4.0940,0.6270) (-3.9933,0.6339) (-3.8926,0.6409) (-3.7919,0.6480) (-3.6913,0.6552) (-3.5906,0.6624) (-3.4899,0.6697) (-3.3893,0.6770) (-3.2886,0.6844) (-3.1879,0.6919) (-3.0872,0.6994) (-2.9866,0.7069) (-2.8859,0.7146) (-2.7852,0.7222) (-2.6846,0.7299) (-2.5839,0.7377) (-2.4832,0.7455) (-2.3826,0.7533) (-2.2819,0.7611) (-2.1812,0.7690) (-2.0805,0.7769) (-1.9799,0.7847) (-1.8792,0.7926) (-1.7785,0.8005) (-1.6779,0.8083) (-1.5772,0.8161) (-1.4765,0.8239) (-1.3758,0.8316) (-1.2752,0.8393) (-1.1745,0.8468) (-1.0738,0.8543) (-0.9732,0.8617) (-0.8725,0.8690) (-0.7718,0.8761) (-0.6711,0.8830) (-0.5705,0.8899) (-0.4698,0.8965) (-0.3691,0.9029) (-0.2685,0.9092) (-0.1678,0.9152) (-0.0671,0.9210) (0.0336,0.9266) (0.1342,0.9319) (0.2349,0.9370) (0.3356,0.9418) (0.4362,0.9464) (0.5369,0.9507) (0.6376,0.9548) (0.7383,0.9586) (0.8389,0.9621) (0.9396,0.9655) (1.0403,0.9685) (1.1409,0.9714) (1.2416,0.9740) (1.3423,0.9765) (1.4430,0.9787) (1.5436,0.9808) (1.6443,0.9826) (1.7450,0.9843) (1.8456,0.9859) (1.9463,0.9873) (2.0470,0.9886) (2.1477,0.9898) (2.2483,0.9908) (2.3490,0.9918) (2.4497,0.9926) (2.5503,0.9934) (2.6510,0.9941) (2.7517,0.9947) (2.8523,0.9953) (2.9530,0.9958) (3.0537,0.9962) (3.1544,0.9966) (3.2550,0.9970) (3.3557,0.9973) (3.4564,0.9976) (3.5570,0.9978) (3.6577,0.9981) (3.7584,0.9983) (3.8591,0.9985) (3.9597,0.9986) (4.0604,0.9988) (4.1611,0.9989) (4.2617,0.9990) (4.3624,0.9991) (4.4631,0.9992) (4.5638,0.9993) (4.6644,0.9994) (4.7651,0.9995) (4.8658,0.9995) (4.9664,0.9996) (5.0671,0.9996) (5.1678,0.9997) (5.2685,0.9997) (5.3691,0.9997) (5.4698,0.9998) (5.5705,0.9998) (5.6711,0.9998) (5.7718,0.9998) (5.8725,0.9998) (5.9732,0.9999) (6.0738,0.9999) (6.1745,0.9999) (6.2752,0.9999) (6.3758,0.9999) (6.4765,0.9999) (6.5772,0.9999) (6.6779,0.9999) (6.7785,0.9999) (6.8792,1.0000) (6.9799,1.0000) (7.0805,1.0000) (7.1812,1.0000) (7.2819,1.0000) (7.3826,1.0000) (7.4832,1.0000) (7.5839,1.0000) (7.6846,1.0000) (7.7852,1.0000) (7.8859,1.0000) (7.9866,1.0000) (8.0872,1.0000) (8.1879,1.0000) (8.2886,1.0000) (8.3893,1.0000) (8.4899,1.0000) (8.5906,1.0000) (8.6913,1.0000) (8.7919,1.0000) (8.8926,1.0000) (8.9933,1.0000) (9.0940,1.0000) (9.1946,1.0000) (9.2953,1.0000) (9.3960,1.0000) (9.4966,1.0000) (9.5973,1.0000) (9.6980,1.0000) (9.7987,1.0000) (9.8993,1.0000) (10.0000,1.0000)};
\node[red!75!black,anchor=west] at (axis cs:1.05,1.98) {(a)};
\node[blue!70!black,anchor=south] at (axis cs:0.5,1.19) {(b)};
\node[green!45!black,anchor=south] at (axis cs:5.0,1.03) {(c)};
\addplot[thin,gray,domain=-5:10] {1};
\addplot[thin,gray,dotted,domain=-5:10] {0.81};
\node[anchor=south west,font=\small,gray!70!black] at (axis cs:-4.9,1.0) {level $A^m$};
\node[anchor=north east,font=\small,gray!70!black] at (axis cs:9.9,0.80) {level $\vartheta^mA^m$, $\vartheta=0.9$};
\end{axis}
\end{tikzpicture}
\caption{The Emden--Fowler function $y(s)=e^{as}V(e^s)$ of Lemma~\ref{lem:A}
(with $b=1$, $m=2$, $\sigma=-1$), normalized by the equilibrium value $A^m$:
(a) $N=3$, $p=5$, so that $p_c<p<p^{**}$ (red); (b) $N=3$, $p=8$, so that
$p>p^{**}$ and \eqref{eq:JL} fails (blue); (c) $N=12$, $p=20$, for which
\eqref{eq:JL} holds (dashed green). In
all cases $y$ reaches every level $\vartheta^mA^m$, $\vartheta<1$
(Lemma~\ref{lem:A}). In cases (a) and (b), it also reaches the level $A^m$,
which is the mechanism of Proposition~\ref{rem:LequalA}: in case (a) it overshoots
and the steady state vanishes at a finite radius, and in case (b) it
oscillates around the stable focus $A^m$. In case (c), the Joseph--Lundgren
condition holds, $A^m$ is a stable node, and $y$ increases to $A^m$ without
crossing it.}
\label{fig:EF}
\end{figure}
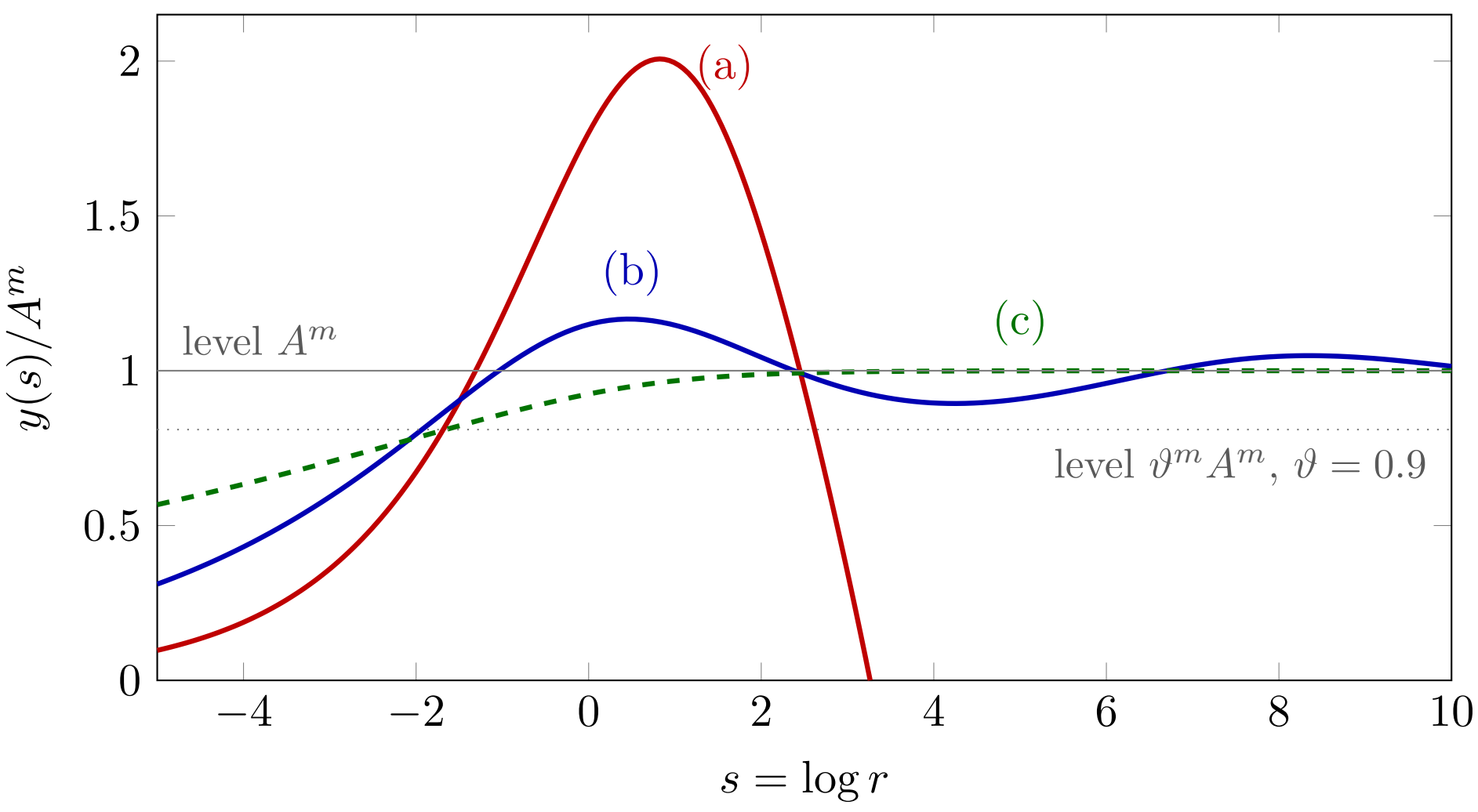

\begin{proof}[Proof of Proposition~\ref{rem:LequalA}]
The proof of Theorem~\ref{thm:conj91} also applies with \(\vartheta=1\), as
soon as the function \(y\) of Lemma~\ref{lem:A} (with \(b=1\)) reaches the
level \(A^m\). More precisely, it suffices that there is \(\hat s<\log R_1\)
with \(y(\hat s)=A^m\) and \(y<A^m\) on \((-\infty,\hat s)\), since the gluing
argument only uses these two properties. In what follows, \(s_1\) and \(Y\)
are as in the proof of Lemma~\ref{lem:A}, and we use the notation
\eqref{eq:EFconst}, so that \eqref{eq:JLfails} reads \(d^2<4(q-1)\kappa_\ast\).
We show that \(\hat s\) exists in the two cases of the proposition, see also
Figure~\ref{fig:EF}. Set
\[
F(y)\eqdef-\frac{\kappa_\ast}{2}y^2+\frac{y^{q+1}}{q+1},\qquad
E_y\eqdef\frac12y'^2+F(y).
\]
By \eqref{eq:EF}, we have \(E_y'=-d\,y'^2\). As \(s\to-\infty\), we have
\(y\to0\) and \(y'\to0\) (see the proof of Lemma~\ref{lem:A}), so that
\(E_y(s)\to0\). Moreover, since \((A^m)^{q-1}=\kappa_\ast\), we have
\(F(A^m)=\kappa_\ast A^{2m}\bigl(\frac1{q+1}-\frac12\bigr)<0\).

\smallskip
\noindent\emph{(a) The case \(\pc<p\le\pss\), that is, \(d\le0\).} Then \(E_y\)
is nondecreasing, and therefore \(E_y\ge0\). We go through the three cases of
the proof of Lemma~\ref{lem:A}. Assume first that \(s_1=+\infty\) and \(Y=A^m\).
Then \(F(y(s))\to F(A^m)<0\), and
\(\frac12 y'^2=E_y-F(y)\ge -F(y(s))\ge-\frac12F(A^m)>0\) for \(s\) large. Thus
\(y'\) is bounded from below by a positive constant for \(s\) large, which
forces \(y\to\infty\) and contradicts \(Y=A^m\). If \(s_1<\log R_1\), then
\(y'(s_1)=0\) and \(E_y(s_1)=F(Y)\ge0\). This means that
\(Y^{q-1}\ge(q+1)\kappa_\ast/2>\kappa_\ast\), whence \(Y>A^m\). In the
remaining case \(s_1=+\infty\) and \(Y=\infty\). In all cases \(Y>A^m\), and
\(y\) reaches the level \(A^m\) on its initial increasing branch.

\smallskip
\noindent\emph{(b) The case \(p>\pss\) when \eqref{eq:JLfails} holds, that is,
\(d>0\) and \(d^2<4(q-1)\kappa_\ast\).} Then \(E_y\) is nonincreasing, and it is
strictly decreasing on \((-\infty,s_1)\), where \(y'>0\). Fix \(s_0<s_1\). We
infer that \(E_y(s)\le E_y(s_0)=:-e_0<0\) for all \(s\ge s_0\). Let
\(Y_1\eqdef\bigl((q+1)\kappa_\ast/2\bigr)^{1/(q-1)}\) be the positive zero of
\(F\). Since \(F(y(s))\le E_y(s)\le-e_0\) for \(s\ge s_0\), and since \(F\) is
negative exactly on \((0,Y_1)\), there are \(0<y_-<y_+<Y_1\) such that
\(y(s)\in[y_-,y_+]\) for all \(s\ge s_0\). In particular, \(y\) does not tend
to \(0\), so that \(R_1=\infty\). Moreover,
\(\frac12y'^2\le-F(y)\le\max_{[0,Y_1]}\abs{F}\). Hence the orbit
\(s\mapsto(y,y')(s)\), \(s\ge s_0\), of the planar system
\(y'=z\), \(z'=-dz+\kappa_\ast y-y^q\) is bounded, and \(E_y\) is a Lyapunov
function for this system. By LaSalle's invariance principle \cite{LaSalle76}, the
\(\omega\)-limit set of the orbit is a nonempty compact invariant set on which
\(E_y'=-dz^2\) vanishes identically. It therefore consists of equilibria lying
in \([y_-,y_+]\times\{0\}\), and it reduces to \(\{(A^m,0)\}\). We conclude that
\((y,y')(s)\to(A^m,0)\) as \(s\to\infty\).

We now use the focus structure of this equilibrium. Set
\(\beta\eqdef\bigl((q-1)\kappa_\ast-d^2/4\bigr)^{1/2}>0\),
\[
X\eqdef y-A^m,\qquad Z\eqdef\frac{1}{\beta}\Bigl(y'+\frac d2X\Bigr),\qquad
g(X)\eqdef\kappa_\ast(A^m+X)-(A^m+X)^q+(q-1)\kappa_\ast X .
\]
Since \(\kappa_\ast A^m-(A^m)^q=0\) and \(q(A^m)^{q-1}=q\kappa_\ast\), we have
\(\abs{g(X)}\le CX^2\) for \(\abs{X}\le A^m/2\). A direct computation based on
\eqref{eq:EF}, which reads \(X''=-dX'-(q-1)\kappa_\ast X+g(X)\), gives
\[
X'=\beta Z-\frac d2X,\qquad Z'=-\beta X-\frac d2Z+\frac{g(X)}{\beta}.
\]
By uniqueness for the planar system, \((X,Z)(s)\ne(0,0)\) for every \(s\).
Hence the polar angle \(\psi\) of \((X,Z)\) is well defined and smooth, and
\[
\psi'=\frac{XZ'-ZX'}{X^2+Z^2}=-\beta+\frac{X\,g(X)}{\beta\,(X^2+Z^2)} .
\]
Since \(\abs{X\,g(X)}\le C\abs{X}^3\le C(X^2+Z^2)^{3/2}\) for \(s\) large, and
since \((X,Z)(s)\to(0,0)\), we obtain \(\psi'(s)\to-\beta<0\). We deduce that
\(\psi(s)\to-\infty\), so that \(X=(X^2+Z^2)^{1/2}\cos\psi\) takes positive
values for arbitrarily large \(s\). Consequently \(y>A^m\) at some point. Since
\(y<A^m\) near \(-\infty\), the first hitting time
\(\hat s\eqdef\inf\{s:\ y(s)=A^m\}\) is finite, and \(y<A^m\) on
\((-\infty,\hat s)\).

\smallskip
In both cases, the proof of Theorem~\ref{thm:conj91} with \(\vartheta=1\) and
\(\rho_1\eqdef e^{\hat s}\) provides a bounded, radially nonincreasing
stationary weak supersolution \(\ol u_1\) of \eqref{eq:CP} in \(\RN\) such
that \(\ol u_1(x)=A\abs{x}^{-\ac}\) for \(\abs{x}\ge1\) and \(\ol u_1\ge A\) for
\(\abs{x}\le1\). Lemma~\ref{lem:comparison}, in its second alternative, and the
blow-up alternative of Assumption~\ref{ass:IL}(A) show that every
\(u_0\in\mathcal A_r\) with \(u_0\le\ol u_1\) generates a global solution
bounded by \(\norm{\ol u_1}_\infty\). Since \(A\phi_{\ac}\le\ol u_1\), we
conclude that \(A\in\mathcal G_{\ac}\).
\end{proof}

\begin{remark}\rm\label{rem:JLnumerics}
When \eqref{eq:JL} holds, the equilibrium \(A^m\) of \eqref{eq:EF} is a stable
node. The numerical computations of Figure~\ref{fig:EF} then suggest that
\(y\) increases to \(A^m\) without crossing it, so that the above gluing is
not available with \(\vartheta=1\). This is consistent with
Conjecture~\ref{conj:threshold}.
\end{remark}

\begin{proof}[Proof of Corollary~\ref{cor:belowS}]
Let \(\vartheta\in(0,1)\), let \(\ol u_\vartheta\) be the barrier of
Theorem~\ref{thm:conj91}, and set
\(v(x)\eqdef\lambda^{\ac}\ol u_\vartheta(\lambda x)\), with \(\lambda\) as in
the statement. By Remark~\ref{rem:scaling}, \(v\) is a bounded, radially
nonincreasing stationary weak supersolution of \eqref{eq:CP} in \(\RN\), and
\(v^m\in H^1(B(0,1))\). By \eqref{eq:barrier-theta}, we have
\(v(x)=\vartheta A\abs{x}^{-\ac}=\vartheta S(x)\) for \(\abs{x}\ge1/\lambda\),
and \(v(x)\ge\lambda^{\ac}\vartheta A=\norm{u_0}_\infty\) for
\(\abs{x}\le1/\lambda\). Hence \(u_0\le\vartheta S=v\) a.e.\ on
\(\{\abs{x}\ge1/\lambda\}\), and \(u_0\le\norm{u_0}_\infty\le v\) a.e.\ on
\(\{\abs{x}<1/\lambda\}\). Lemma~\ref{lem:comparison}, applied in its second
alternative, gives \(u(t)\le v\) on \([0,\Tmax(u_0))\). In particular,
\(\norm{u(t)}_\infty\le\lambda^{\ac}\norm{\ol u_\vartheta}_\infty\), and the
blow-up alternative of Assumption~\ref{ass:IL}(A) yields \(\Tmax(u_0)=\infty\).
This proves \eqref{eq:belowS}. Under the assumptions of
Proposition~\ref{rem:LequalA}, the same argument applies with \(\vartheta=1\)
and \(\ol u_1\) in place of \(\ol u_\vartheta\), since \(\ol u_1\) has the same
two properties \eqref{eq:barrier-theta} with \(\vartheta=1\).
\end{proof}

%===================================================
\section{{Blow-up and the second critical exponent}}\label{sec:proofs-blowup}

In this section we prove the blow-up results of
Section~\ref{subsec:mainresults} and the characterization of the second
critical exponent.

\subsection{{Blow-up for slow decay and for large amplitude threshold}}
We begin with Theorem~\ref{thm:slow}. The strategy is the following. The slow
decay assumption \eqref{eq:slowtail} allows us to bound $u_0$ from below on a
far annulus $\{\lambda<|x|<2\lambda\}$ by a compactly supported datum of the
form $v_{\kappa\lambda^{-\gamma},\lambda}$. By the scaling computation of
Lemma~\ref{lem:energy-annulus}, the energy of this datum becomes negative for
$\lambda$ large as soon as $\gamma<\alpha_c$, whatever the amplitude. The
negative energy criterion of \cite[Theorem~2.9]{IL1} then gives finite-time
blow-up for $v_{\kappa\lambda^{-\gamma},\lambda}$, and the comparison principle
transfers it to $u_0$.

\begin{proof}[Proof of Theorem~\ref{thm:slow}]
Set \(\kappa_1\eqdef\min\{1,\kappa_0/2\}\in(0,\kappa_0)\), with the convention
\(\kappa_0/2=\infty\) if \(\kappa_0=\infty\). By \eqref{eq:slowtail}, there is
\(R_0\ge1\) such that \(\abs{x}^{\gamma}u_0(x)\ge\kappa_1\) for a.e.\
\(\abs{x}\ge R_0\). Indeed, \(\kappa_1<\kappa_0\), and the definition of the
lower limit provides such an \(R_0\).

Let \(\lambda\ge R_0\). On the annulus \(\{\lambda<|x|<2\lambda\}\), we then have
\[
u_0(x)\ \ge\ \kappa_1\,(2\lambda)^{-\gamma}
=\kappa\,\lambda^{-\gamma},\qquad \kappa\eqdef 2^{-\gamma}\kappa_1>0 .
\]
With the notation \eqref{eq:annulus-data}, and since \(0\le\zeta\le1\) is
supported in \(\{1<\abs{y}<2\}\), we infer that
\[
u_0\ \ge\ v_{\kappa\lambda^{-\gamma},\lambda}\qquad\text{a.e. on }\mathbb{R}^N .
\]
By Lemma~\ref{lem:energy-annulus} and \(\gamma<\alpha_c\), the exponent
\(\sigma+2-\gamma(p-m)\) is positive, so that the bracket in
\eqref{eq:Escaling2} tends to \(-\infty\) as \(\lambda\to\infty\). Hence there
is \(\lambda_0=\lambda_0(N,m,p,\sigma,\gamma,\kappa,\zeta)\ge R_0\) with
\[
E\bigl(v_{\kappa\lambda^{-\gamma},\lambda}\bigr)<0
\qquad\text{for all }\lambda\ge\lambda_0 .
\]
Fix such a \(\lambda\). Since \(v_{\kappa\lambda^{-\gamma},\lambda}\) is bounded
with compact support, it belongs to \(\mathcal{A}_r\). By
Lemma~\ref{lem:energy-annulus}, it is admissible in
Theorem~\ref{thm:IL-super}(b), which gives
\(\Tmax(v_{\kappa\lambda^{-\gamma},\lambda})<\infty\). Moreover,
\(u_0\in\mathcal{A}_r\) by hypothesis and
\(v_{\kappa\lambda^{-\gamma},\lambda}\le u_0\). The comparison principle of
Assumption~\ref{ass:IL}(B) therefore yields
\[
T_{\max}(u_0)\ \le\ T_{\max}\bigl(v_{\kappa\lambda^{-\gamma},\lambda}\bigr)\ <\ \infty . \qedhere
\]
\end{proof}
We now come to the proof of Theorem~\ref{thm:blowup-crit}. The argument
parallels that of Theorem~\ref{thm:slow}, with one essential difference. At
$\gamma=\alpha_c$, the exponent $\sigma+2-\alpha_c(p-m)$ vanishes, so that
the bracket in the energy of the compactly supported datum
$v_{\kappa\lambda^{-\alpha_c},\lambda}$ is independent of $\lambda$. We show
that this energy is negative as soon as the tail amplitude exceeds $L_0$, and
we conclude by the negative energy criterion and the comparison principle.

\begin{proof}[Proof of Theorem~\ref{thm:blowup-crit}]
Since \(\liminf_{|x|\to\infty}|x|^{\alpha_c}u_0(x)>L_0\), there exist
\(\delta>0\) and \(R_0\ge1\) such that
\[
u_0(x)\ge\kappa_1|x|^{-\alpha_c}\qquad\text{for a.e. }|x|\ge R_0,
\qquad\text{with }\kappa_1\eqdef L_0+\delta>L_0 .
\]
We set \(\kappa\eqdef 2^{-\alpha_c}\kappa_1\). As in the proof of
Theorem~\ref{thm:slow}, we have
\(u_0\ge v_{\kappa\lambda^{-\alpha_c},\lambda}\) for every \(\lambda\ge R_0\).
Since \(\sigma+2-\alpha_c(p-m)=0\) by \eqref{eq:ac}, the factor
\(\lambda^{\sigma+2-\alpha_c(p-m)}\) in \eqref{eq:Escaling2} equals \(1\), and
we obtain
\[
E\bigl(v_{\kappa\lambda^{-\alpha_c},\lambda}\bigr)
=\lambda^{N-2-2m\alpha_c}\left[\frac{\kappa^{2m}}{2}\,I_1
-\frac{m}{m+p}\,\kappa^{m+p}\,I_2\right].
\]
The prefactor \(\lambda^{N-2-2m\alpha_c}\) is positive, so that the sign of
\(E\) is that of the bracket. Dividing the bracket by
\(\frac{m}{m+p}\kappa^{2m}I_2>0\), we see that
\(E(v_{\kappa\lambda^{-\alpha_c},\lambda})<0\) if and only if
\[
\kappa^{\,p-m}>\frac{(m+p)I_1}{2m\,I_2}=M_0^{\,p-m},
\qquad\text{that is,}\qquad
\kappa_1=2^{\alpha_c}\kappa>2^{\alpha_c}M_0=L_0,
\]
by \eqref{eq:M0} and \eqref{eq:L0}. This holds by the choice of \(\kappa_1\).

Fix any \(\lambda\ge R_0\). By Lemma~\ref{lem:energy-annulus}, the datum
\(v_{\kappa\lambda^{-\alpha_c},\lambda}\) is admissible in
Theorem~\ref{thm:IL-super}(b), which gives finite-time blow-up for it.
Assumption~\ref{ass:IL}(B) then transfers the blow-up to \(u_0\). Finally, for
\(u_0=L\phi\), we have \(\lim_{|x|\to\infty}|x|^{\alpha_c}L\phi(x)=L\), so that
\eqref{eq:crittail} holds whenever \(L>L_0\).
\end{proof}

\subsection{{The second critical exponent and the threshold function}}

We now prove Theorem~\ref{thm:sce}. Part (i) is a direct consequence of the
blow-up theorem for slow decay. Part (ii) relies on the small-data result of
\cite{IL1} in the critical Lebesgue space, together with the bounded barrier at
criticality.
\begin{proof}[Proof of Theorem~\ref{thm:sce}]
We begin with (i). Let \(\gamma\in(0,\alpha_c)\), and assume that
\(u_0\in\mathcal{A}_r\) satisfies
\(\liminf_{|x|\to\infty}|x|^{\gamma}u_0(x)>0\). This lower limit is
\(\kappa_0\in(0,\infty]\), so that the assumptions of Theorem~\ref{thm:slow}
are satisfied, and Theorem~\ref{thm:slow} yields \(T_{\max}(u_0)<\infty\). Data
satisfying \eqref{eq:tail} are a particular case, with \(\kappa_0=L\).

For (ii), we first note that \eqref{eq:smallweighted} is equivalent to
\(u_0(x)\le\varepsilon_\gamma\min\{1,|x|^{-\gamma}\}=\varepsilon_\gamma\phi_\gamma(x)\)
for all \(x\), by \eqref{eq:phigamma}. Since \(\gamma>\alpha_c=N/(r_0+1)\), we
have \(\gamma(r_0+1)>N\), and therefore
\begin{align*}
\|\phi_\gamma\|_{r_0+1}^{r_0+1}
&=\int_{|x|\le1}dx+\int_{|x|>1}|x|^{-\gamma(r_0+1)}\,dx\\
&=\omega_N+N\omega_N\int_1^\infty\varrho^{N-1-\gamma(r_0+1)}\,d\varrho
=\omega_N+\frac{N\omega_N}{\gamma(r_0+1)-N},
\end{align*}
which is \eqref{eq:Cgamma}. We infer that
\(u_0\in L^{r_0+1}_+(\mathbb{R}^N)\cap L^\infty(\mathbb{R}^N)\) and
\[
\|u_0\|_{r_0+1}\le\varepsilon_\gamma\|\phi_\gamma\|_{r_0+1}
=\varepsilon_\gamma C_\gamma=C_0 .
\]
Theorem~\ref{thm:IL-super}(a) then gives \(T_{\max}(u_0)=\infty\) and
\(u\in L^{\infty}\bigl((0,\infty),L^{r}(\mathbb{R}^N)\bigr)\) for all
\(r\in[r_0+1,\infty)\). Finally, if \(p>p_c\) and
\(u_0\le\vartheta\overline{u}\) for some \(\vartheta\in(0,1]\), then
Theorem~\ref{thm:global} gives \eqref{eq:global-bound}.
\end{proof}

We next prove Proposition~\ref{prop:nodecayonly}. It shows that the smallness
condition \eqref{eq:smallweighted} cannot be replaced by a condition on \(L\)
and \(\norm{u_0}_\infty\) alone.

\begin{proof}[Proof of Proposition~\ref{prop:nodecayonly}]
Fix \(k>\gamma\) and, for \(\lambda\ge1\), set
\[
u_0^{(\lambda)}(x)\eqdef\max\Bigl\{\eps\min\bigl\{1,(\abs{x}/\lambda)^{-k}\bigr\},\
L\,\phi_\gamma(x)\Bigr\}.
\]
Then \(u_0^{(\lambda)}\) is bounded, with \(\norm{u_0^{(\lambda)}}_{\infty}=\eps\)
because \(L\le\eps\), the value \(\eps\) being attained at \(x=0\). For
\(\abs{x}\ge\lambda\ (\ge1)\), we have
\[
\abs{x}^{\gamma}u_0^{(\lambda)}(x)
=\max\bigl\{\eps\lambda^k\abs{x}^{\gamma-k},L\bigr\}
\xrightarrow[\abs{x}\to\infty]{}L,
\]
since \(k>\gamma\). In particular, \(u_0^{(\lambda)}\lesssim\abs{x}^{-\gamma}\) at
infinity, so that \(u_0^{(\lambda)}\in\mathcal A_r\) for every
\(r>\max\{r_0,N/\gamma-1\}\). On the other hand, \(u_0^{(\lambda)}\ge\eps\) on
\(B(0,\lambda)\). With \(\zeta\) as in \eqref{eq:zeta} and
\(\lambda'\eqdef\lambda/2\), we infer that
\[
u_0^{(\lambda)}\ \ge\ \eps\,\zeta_{\lambda'}=v_{\eps,\lambda'}
\qquad\text{on }\RN,
\]
because \(\supp\zeta_{\lambda'}\subset\{\lambda'<\abs{x}<2\lambda'\}
\subset B(0,\lambda)\) and \(0\le\zeta\le1\). By \eqref{eq:Esign}, we have
\(E(v_{\eps,\lambda'})<0\) as soon as
\(\lambda'^{\,\sigma+2}>(M_0/\eps)^{p-m}\), which holds for \(\lambda\) large
since \(\sigma+2>0\). For such a \(\lambda\), Theorem~\ref{thm:IL-super}(b)
gives \(\Tmax(v_{\eps,\lambda'})<\infty\), and Assumption~\ref{ass:IL}(B)
yields \(\Tmax(u_0^{(\lambda)})\le \Tmax(v_{\eps,\lambda'})<\infty\). Thus
\(u_0=u_0^{(\lambda)}\) has the required properties.
\end{proof}

We turn to the proof of Theorem~\ref{thm:threshold}, which is divided into
five parts. Part (i) establishes the dichotomy along each ray by comparison.
Part (ii) shows the finiteness and the monotonicity of $L_\ast$, using the
energy criterion and the comparison principle. Part (iii) follows from the
blow-up theorem for slow decay, part (iv) from the small-data theorem for fast
decay, and part (v) combines the bounded barrier at criticality with the
blow-up result for large amplitude.
\begin{proof}[Proof of Theorem~\ref{thm:threshold}]
Throughout the proof, \(L\phi_\gamma\in\mathcal{A}_r\) for every
\(r>\max\{r_0,N/\gamma-1\}\). Indeed, \(L\phi_\gamma\) is bounded, and
\((L\phi_\gamma)^{r+1}\le L^{r+1}|x|^{-\gamma(r+1)}\) for \(|x|\ge1\), which is
integrable at infinity when \(\gamma(r+1)>N\). Hence Assumption~\ref{ass:IL}
applies to all the data considered, and, since the classes \(\mathcal A_r\) are
increasing in \(r\), any two of them belong to a common class.

For (i), let \(0<L<L_\ast(\gamma)\). Since
\(L_\ast(\gamma)=\sup\mathcal{G}_\gamma>L\), there is
\(L'\in\mathcal{G}_\gamma\) with \(L'>L\). Then \(L\phi_\gamma\le L'\phi_\gamma\),
and Assumption~\ref{ass:IL}(B) gives
\(T_{\max}(L\phi_\gamma)\ge T_{\max}(L'\phi_\gamma)=\infty\). If
\(L>L_\ast(\gamma)\), then \(L\notin\mathcal{G}_\gamma\), that is,
\(T_{\max}(L\phi_\gamma)<\infty\).

For (ii), let \(M\eqdef 2^{-\gamma}L\). On \(\supp\zeta\subset\{1<|x|<2\}\), we
have \(\phi_\gamma(x)=|x|^{-\gamma}\ge2^{-\gamma}\) and \(\zeta\le1\), so that
\(M\zeta(x)\le2^{-\gamma}L\le L\phi_\gamma(x)\). We deduce that
\(v_{M,1}\le L\phi_\gamma\) in \(\RN\). By \eqref{eq:Esign} with
\(\lambda=1\), we have \(E(v_{M,1})<0\) as soon as \(M>M_0\). For such \(M\),
Theorem~\ref{thm:IL-super}(b) gives \(T_{\max}(v_{M,1})<\infty\), and
Assumption~\ref{ass:IL}(B) yields
\(T_{\max}(L\phi_\gamma)\le T_{\max}(v_{M,1})<\infty\). Hence
\(L\notin\mathcal{G}_\gamma\) whenever \(L>2^{\gamma}M_0\), which proves
\(L_\ast(\gamma)\le2^{\gamma}M_0\). Monotonicity is immediate. If
\(\gamma_1\le\gamma_2\), then \(\phi_{\gamma_2}\le\phi_{\gamma_1}\), so that
\(\mathcal{G}_{\gamma_1}\subseteq\mathcal{G}_{\gamma_2}\) by
Assumption~\ref{ass:IL}(B), and therefore \(L_\ast(\gamma_1)\le L_\ast(\gamma_2)\).

For (iii), let \(\gamma\in(0,\alpha_c)\) and \(L>0\). Since
\(|x|^{\gamma}L\phi_\gamma(x)=L\) for \(|x|\ge1\), Theorem~\ref{thm:slow}
gives \(T_{\max}(L\phi_\gamma)<\infty\). Hence \(\mathcal{G}_\gamma=\emptyset\)
and \(L_\ast(\gamma)=0\).

For (iv), \(L\phi_\gamma\) satisfies \eqref{eq:smallweighted} whenever
\(L\le\varepsilon_\gamma\), so that Theorem~\ref{thm:sce}(ii) gives
\((0,\varepsilon_\gamma]\subset\mathcal{G}_\gamma\).

Finally, for (v), let \(0<L\le\mu A\). By the proof of
Theorem~\ref{thm:global}, we have \(L\phi\le\mu A\phi\le\overline{u}\), and
Theorem~\ref{thm:global} gives \(T_{\max}(L\phi)=\infty\) together with
\eqref{eq:global-bound-L}. Hence \((0,\mu A]\subset\mathcal{G}_{\alpha_c}\) and
\(L_\ast(\alpha_c)\ge\mu A\). If \(L>L_0\), then
\(\lim_{|x|\to\infty}|x|^{\alpha_c}L\phi(x)=L>L_0\), and
Theorem~\ref{thm:blowup-crit} yields \(T_{\max}(L\phi)<\infty\). Hence
\(\mathcal{G}_{\alpha_c}\cap(L_0,\infty)=\emptyset\), and
\(L_\ast(\alpha_c)\le L_0\).
\end{proof}

We conclude this section with the proof of Corollary~\ref{cor:char}.

\begin{proof}[Proof of Corollary~\ref{cor:char}]
By Theorem~\ref{thm:threshold}(iii), the set \(\{\gamma>0:\ L_\ast(\gamma)=0\}\)
contains \((0,\ac)\), and by Theorem~\ref{thm:threshold}(iv) it is contained in
\((0,\ac]\). Both displayed quantities therefore equal \(\ac\), whether or not
\(L_\ast(\ac)=0\). The last assertion is Theorem~\ref{thm:threshold}(v).
\end{proof}

%===================================================
\section{{Refinements, conjectures and open problems}}\label{sec:conjectures}
%===================================================

Theorems~\ref{thm:sce}, \ref{thm:threshold} and \ref{thm:conj91} determine the
second critical exponent, and they show that \(A\le L_\ast(\ac)\le L_0\) for
\(p>\pc\). In this section we discuss the identification of \(L_\ast(\ac)\),
the influence of slowly varying corrections to the critical tail, and the range
\(\pF<p\le\pc\), in which no singular steady state is available. Throughout,
proved statements and conjectures are kept strictly separate.

\subsection{Towards the identification of \texorpdfstring{$L_\ast(\ac)$}{L*(alpha\_c)}}

Linearizing \eqref{eq:EF} at the equilibrium \(y=A^m\), we obtain
\(z''+dz'+(q-1)\kappa_\ast z=0\), where \((q-1)\kappa_\ast=(\sigma+2)(N-2-a)\).
For \(p>\pss\), that is, \(d>0\), the equilibrium is therefore a stable node if
and only if
\begin{equation}\label{eq:JL}
\bigl(N-2-2a\bigr)^2\ \ge\ 4(\sigma+2)(N-2-a),\qquad a=\ac m=\frac{m(\sigma+2)}{p-m},
\end{equation}
and it is a stable focus otherwise. We say that \(p\) satisfies the
\emph{Joseph--Lundgren condition} if \(p>\pss\) and \eqref{eq:JL} holds. For
\(m=1\) and \(\sigma=0\), this is the classical condition \(p\ge p_{JL}\) of
Joseph and Lundgren \cite{JosephLundgren73}, which governs the ordering of the
regular steady states below the singular one and the stability of the latter
\cite{GuiNiWang92,Wang93,PolacikYanagida03}. As \(p\to\infty\), we have
\(a\to0\), and \eqref{eq:JL} reduces to \(N-2\ge4(\sigma+2)\). Thus the
Joseph--Lundgren condition holds for all large \(p\) if \(N>10+4\sigma\).
{This condition coincides with the Joseph--Lundgren exponent
\[
p_{JL}(\sigma)\eqdef m\,\frac{(N-2)^2-2(N+\sigma)(\sigma+2)
+2(\sigma+2)\sqrt{(2N+\sigma-2)(\sigma+2)}}{(N-2)(N-4\sigma-10)}
\quad\text{if }N>10+4\sigma,\qquad p_{JL}(\sigma)\eqdef\infty\quad\text{otherwise},
\]
introduced in \cite[Eq.~(2.4)]{ISSpeak}: \(p\) satisfies the
Joseph--Lundgren condition if and only if \(p\ge p_{JL}(\sigma)\). Indeed,
set \(t\eqdef p-m\) and \(t^{\ast\ast}\eqdef\pss-m\). Since
\(\sigma+2=a(q-1)\), the difference of the two sides of \eqref{eq:JL} equals
\((N-2)^2-4qa(N-2-a)\), and multiplying it by \(t^2\) we find that
\eqref{eq:JL} is equivalent to
\[
\mathcal Q(t)\eqdef(N-2)(N-4\sigma-10)\,t^{2}-4m(\sigma+2)(N-\sigma-4)\,t
+4m^{2}(\sigma+2)^{2}\ \ge\ 0 .
\]
At \(t=t^{\ast\ast}\) we have \(a=(N-2)/2\), so that the left-hand side of
\eqref{eq:JL} vanishes while its right-hand side is positive, whence
\(\mathcal Q(t^{\ast\ast})<0\). Moreover \(\mathcal Q(0)>0\). If
\(N\le10+4\sigma\), the function \(\mathcal Q\) is concave or affine, and
therefore \(\mathcal Q<0\) on \([t^{\ast\ast},\infty)\). If
\(N>10+4\sigma\), the function \(\mathcal Q\) is convex and \(t^{\ast\ast}\)
lies strictly between its two roots, so that, for \(t>t^{\ast\ast}\),
\(\mathcal Q(t)\ge0\) if and only if \(t\) is at least the larger root.
The discriminant of \(\mathcal Q\) equals
\(16m^2(\sigma+2)^3(2N+\sigma-2)\), and a direct computation shows that the
larger root is \(p_{JL}(\sigma)-m\). For \(\sigma\in(-2,0)\), the exponent
\(p_{JL}(\sigma)\) may be finite in dimensions \(N\le10\), as already
observed in \cite{ISSpeak}.}
Since \(a(N-2-a)=\kappa_\ast\), a direct expansion shows that \eqref{eq:JL}
is equivalent to
\[
q\,\kappa_\ast=\frac pm\,A^{p-m}\le\frac{(N-2)^2}{4}.
\]
This is the condition under which the Hardy inequality controls the
linearization of \(v\mapsto\Delta v+\abs{x}^{\sigma}v^{q}\) at \(v=S^m\),
namely
\[
\varphi\longmapsto\Delta\varphi+q\abs{x}^{\sigma}(S^m)^{q-1}\varphi
=\Delta\varphi+q\kappa_\ast\abs{x}^{-2}\varphi .
\] By Proposition~\ref{rem:LequalA}, the
datum \(A\phi_{\ac}\) generates a global and bounded solution whenever
\(p>\pc\) and the Joseph--Lundgren condition fails.

\begin{remark}\rm \label{rem:noglue}
For \(\vartheta>1\), the function \(\vartheta S\) is a strict subsolution of the
stationary equation in \(\RN\setminus\{0\}\), since
\(\vartheta^{p}>\vartheta^{m}\) (compare with Lemma~\ref{lem:submultiple}).
However, this does not directly produce bounded subsolutions with tail
amplitude \(\vartheta A\). Indeed, a function glued from two radial branches
is a subsolution only if the jump of the normal derivative of its \(m\)-th
power across the matching sphere is nonnegative. Since \(\vartheta S\)
dominates every bounded function near the origin, a bounded cap lies below
\(\vartheta S\) inside the matching sphere. As in the proof of
Theorem~\ref{thm:barrier}, the jump is then nonpositive, that is, it has the
supersolution sign \eqref{eq:jump}, and it vanishes only for a \(C^1\)
matching. We do not know how to construct bounded subsolutions with exactly
critical tail \(\vartheta A\), \(\vartheta>1\), in this way. On the other
hand, for the semilinear model \(m=1\), \(\sigma=0\), Souplet and Weissler
\cite{SouWei03} proved the existence of positive regular self-similar
solutions, hence of global solutions, whose initial value is \(\vartheta\)
times the singular steady state, for every \(\vartheta>1\) close to \(1\). The
case \(\vartheta=1\) is due to Galaktionov and V\'azquez
\cite{GalaktionovVazquez97}. This holds exactly in the range between the
Serrin and the Joseph--Lundgren exponents, and this range is optimal
\cite{SouWei03}. It suggests that, for \eqref{eq:CP}, one should expect
\(L_\ast(\ac)>A\) when the Joseph--Lundgren condition fails.
\end{remark}

\begin{conjecture}[Identification of the threshold]\label{conj:threshold}
Let \(N\ge3\), \(m>1\), \(\sigma\in(-2,0)\), and let \(p\) satisfy the
Joseph--Lundgren condition \eqref{eq:JL}. Then the threshold
\(L_\ast=L_\ast(\ac)\) of Theorem~\ref{thm:threshold} equals the amplitude of
the singular steady state:
\[
L_\ast=A=\bigl[\ac m(N-2-\ac m)\bigr]^{\frac{1}{p-m}} .
\]
Equivalently, in view of Theorem~\ref{thm:conj91}, the solution with datum
\(L\phi_{\ac}\) blows up in finite time for every \(L>A\). For \(p>\pc\) not
satisfying \eqref{eq:JL}, we expect \(L_\ast(\ac)>A\).
\end{conjecture}

\begin{conjecture}[Attraction to the singular steady state]\label{conj:attract}
Under the hypotheses of Conjecture~\ref{conj:threshold}, let
\(u_0\in\mathcal A_r\) satisfy \(u_0\le S\) and
\(\lim_{\abs{x}\to\infty}\abs{x}^{\ac}u_0(x)=A\). Then the solution is global
and converges to the singular steady state,
\[
\abs{x}^{\ac}u(t,x)\longrightarrow A\qquad(t\to\infty),
\]
locally uniformly in \(\RN\setminus\{0\}\).
\end{conjecture}
\noindent For \(L<A\), global existence and uniform boundedness of the
solution with datum \(L\phi_{\ac}\) are established in
Theorem~\ref{thm:conj91}, and the same holds for \(L=A\) when \eqref{eq:JL}
fails, by Proposition~\ref{rem:LequalA}. In all these cases, the tail bound
\(u(t,x)\le A\abs{x}^{-\ac}\) is preserved for all times by
Theorem~\ref{thm:trap}.

\begin{remark}\rm\label{rem:peaking}
The peaking solutions of \cite{ISSpeak} show that, at the critical decay, the
tail amplitude alone does not determine the outcome, even below \(A\). For
\(\sigma\in(-2,0)\) and \(p\in(\pss,p_{JL}(\sigma))\)
\cite[Theorem~2.2]{ISSpeak}, as well as in the further ranges of
\cite[Theorems~2.3, 2.4 and Corollary~2.5]{ISSpeak}, there are backward
self-similar solutions \(U(t,x)=(T-t)^{-\alpha}f\bigl(\abs{x}(T-t)^{-\beta}\bigr)\)
which blow up at \(t=T\), and whose profiles satisfy
\(\xi^{\ac}f(\xi)\to K\in(0,A)\) as \(\xi\to\infty\)
\cite[Eq.~(2.6)]{ISSpeak}. Their values at \(t=0\) are bounded, with exactly
critical tail of amplitude \(K<A\). This does not contradict
Theorem~\ref{thm:conj91} or Corollary~\ref{cor:belowS}, which require a
pointwise domination by \(\ol u_\vartheta\) or by \(\vartheta S\) on the
whole of \(\RN\), and not only at infinity. It strengthens, at
\(\gamma=\ac\), the conclusion of Proposition~\ref{prop:nodecayonly}, and it
explains why the threshold \(L_\ast\) is defined along the rays
\(L\mapsto L\phi_\gamma\). We do not discuss here whether these
self-similar solutions coincide with the solution of Assumption~\ref{ass:IL}
for the same initial value, which would require a uniqueness result for
bounded data with critical tail. On the other hand, for \(p>\pss\) and
\(0<K<A\), the unbounded initial value \(K\abs{x}^{-\ac}\) produces a global
solution, bounded for \(t>0\) and of forward self-similar form
\cite[Proposition~3.1]{ISSpeak}. This is consistent with the lower bound
\(L_\ast(\ac)\ge A\) of Corollary~\ref{cor:sandwich}.
\end{remark}

\subsection{Slowly varying corrections at the critical decay}

We now propose a refinement at the exact borderline, which connects the problem
to the theory of regular variation \cite{BGT87,Seneta76}. Write the critical
datum as
\begin{equation}\label{eq:karamata-data}
u_0(x)=\abs{x}^{-\ac}\,\ell(\abs{x})\quad\text{for }\abs{x}\ge1,\qquad
\ell\ \text{slowly varying at $+\infty$ in the sense of Karamata},
\end{equation}
so that \(\ell(\lambda r)/\ell(r)\to1\) as \(r\to\infty\) for every
\(\lambda>0\). The pure power case \(\ell\equiv L\) is \eqref{eq:tail} with
\(\gamma=\ac\). Our results already describe the two extreme regimes. If
\(\limsup_{r\to\infty}\ell(r)<A\) and \(u_0\) is small enough near the origin,
then the solution is global by Theorem~\ref{thm:conj91}. If
\(\liminf_{r\to\infty}\ell(r)>L_0\), then the solution blows up by
Theorem~\ref{thm:blowup-crit}. Both criteria depend only on the \emph{limit} of
\(\ell\). The genuinely borderline situation is \(\ell(r)\to L_\ast\), where the
limit carries no information and the outcome must be decided by the rate and
the sign of the approach.

\begin{proposition}\label{prop:karamata}
Let \(N\ge3\), \(m>1\), \(\sigma\in(-2,0)\) and \(p>\pc\), let \(L_0\) be as in
\eqref{eq:L0}, and let \(u_0\in\mathcal{A}_r\) for some \(r>r_c\) satisfy
\(u_0(x)=\abs{x}^{-\ac}\ell(\abs{x})\) for \(\abs{x}\ge1\), where
\(\ell:[1,\infty)\to(0,\infty)\) is measurable.
\begin{enumerate}[label=\textup{(\roman*)},leftmargin=2.2em]
\item If there are \(\vartheta\in(0,1)\) and \(R_0\ge1\) such that
\(\ell(r)\le\vartheta A\) for all \(r\ge R_0\) and
\(u_0(x)\le\vartheta A\,R_0^{-\ac}\) for a.e.\ \(\abs{x}\le R_0\), then
\(\Tmax(u_0)=\infty\) and
\[
u(t,x)\le R_0^{-\ac}\,\ol u_\vartheta(x/R_0)\le R_0^{-\ac}\norm{\ol u_\vartheta}_\infty
\qquad\text{for all }t\ge0\text{ and a.e. }x\in\RN,
\]
with \(\ol u_\vartheta\) as in Theorem~\ref{thm:conj91}.
\item If \(\liminf_{r\to\infty}\ell(r)>L_0\), then \(\Tmax(u_0)<\infty\).
\end{enumerate}
\end{proposition}
The slowly varying character of \(\ell\) plays no role in
Proposition~\ref{prop:karamata}, which only involves the upper and lower limits
of \(\ell\). It becomes relevant at the borderline \(\ell(r)\to L_\ast\), see
Open Problem~\ref{conj:karamata}.

\begin{proof}
(i) By Remark~\ref{rem:scaling} with \(\lambda=1/R_0\), the function
\(v(x)\eqdef R_0^{-\ac}\ol u_\vartheta(x/R_0)\) is a bounded, radially
nonincreasing stationary weak supersolution in \(\RN\). By
\eqref{eq:barrier-theta}, we have \(v(x)=\vartheta A\abs{x}^{-\ac}\) for
\(\abs{x}\ge R_0\) and \(v(x)\ge\vartheta AR_0^{-\ac}\) for \(\abs{x}\le R_0\).
For \(\abs{x}\ge R_0\), we obtain
\(u_0(x)=\abs{x}^{-\ac}\ell(\abs{x})\le\vartheta A\abs{x}^{-\ac}=v(x)\), and for
\(\abs{x}\le R_0\) we have \(u_0(x)\le\vartheta AR_0^{-\ac}\le v(x)\). Hence
\(u_0\le v\) a.e.\ in \(\RN\). Lemma~\ref{lem:comparison}, applied in its
second alternative, gives \(u(t)\le v\) on \([0,\Tmax(u_0))\), and
Assumption~\ref{ass:IL}(A) yields \(\Tmax(u_0)=\infty\).

(ii) By assumption, there are \(\kappa_1>L_0\) and \(R_0\ge1\) such that
\(\ell(r)\ge\kappa_1\) for all \(r\ge R_0\). Hence
\(\liminf_{\abs{x}\to\infty}\abs{x}^{\ac}u_0(x)\ge\kappa_1>L_0\), and
Theorem~\ref{thm:blowup-crit} gives \(\Tmax(u_0)<\infty\).
\end{proof}

\begin{openproblem}[Karamata criterion at the critical decay]\label{conj:karamata}
Let \(N\ge3\), \(m>1\), \(\sigma\in(-2,0)\), \(p>\pc\), and let \(u_0\) be as
in \eqref{eq:karamata-data} with \(\ell\) slowly varying and
\(\lim_{r\to\infty}\ell(r)=L_\ast\). Decide the alternative between global
existence and blow-up in terms of the slowly varying correction
\(\ell-L_\ast\). Concretely:
\begin{enumerate}[label=\textup{(\alph*)},leftmargin=2.2em]
\item Is it true that \(\ell(r)\le L_\ast\) for all large \(r\) implies global
existence?
\item For the model family
\(\ell_{\delta,\beta}(r)=L_\ast\bigl(1+\delta(\log r)^{-\beta}\bigr)\), with
\(\delta>0\) and \(\beta>0\), determine the set of pairs \((\delta,\beta)\) for
which the solution blows up in finite time. Heuristically, the annulus energy
computation of Lemma~\ref{lem:energy-annulus} becomes logarithmically neutral
at \(\gamma=\ac\), and the balance is then governed by an integral of the form
\(\int^{\infty}\bigl(\ell(r)-L_\ast\bigr)\,\frac{dr}{r\,\Lambda(r)}\) for a
slowly varying gauge \(\Lambda\). One expects blow-up when this integral
diverges to \(+\infty\), and global existence when it converges or diverges to
\(-\infty\).
\end{enumerate}
\end{openproblem}

We stress that Open Problem~\ref{conj:karamata} is stated as an open problem,
not as a theorem in disguise. The schematic integral criterion above is a
heuristic. It is suggested by the scale invariance of the annulus energy at
\(\gamma=\ac\) (Remark~\ref{rem:whycritical}), and by the way Karamata's
integral machinery governs threshold phenomena in related Fujita-type problems
with regularly varying coefficients, see for instance
\cite{BBM25,KumarMajdoub26,Majdoub26EJQTDE}. Its
precise formulation, and in particular the identification of the gauge
\(\Lambda\), is open.

\subsection{{The range \texorpdfstring{$\pF<p\le\pc$}{pF<p<=pc} and self-similar solutions}}

\begin{openproblem}[The range $\pF<p\le\pc$]\label{prob:pcrange}
Let \(N\ge3\), \(m>1\), \(\sigma\in(-2,0)\) and \(\pF<p\le\pc\).
Theorem~\ref{thm:sce} still identifies \(\ac\) as the second critical exponent
(its proof uses only \(p>\pF\)), and Theorem~\ref{thm:threshold}(i)--(iv)
remains valid. However, no algebraic singular steady state exists
(Proposition~\ref{prop:steady}), so that the barriers of
Theorems~\ref{thm:barrier} and \ref{thm:conj91} are not available.
Theorem~\ref{thm:allp}, proved below, shows that
\(0<\sup_{b}c_\infty(b)\le L_\ast(\ac)\le L_0\), the barriers being the global
self-similar solutions of \cite{ISS26} with algebraic decay. Determine
\(L_\ast(\ac)\) in this range, and
decide whether it coincides with the supremum of the tail constants
\(c_\infty\) of the global self-similar solutions. In the semilinear case, the
corresponding threshold role of forward self-similar solutions with algebraic
tails is analysed in \cite{Naito12,Naito20}.
\end{openproblem}

\begin{remark}\rm \label{rem:selfsimilar}
Forward self-similar solutions can provide barriers in the whole range
\(p>\pF\), including \(\pF<p\le\pc\), where no singular steady state exists.
Looking for solutions of \eqref{eq:CP} of the form
\(u(t,x)=t^{-\alpha}f\bigl(\abs{x}t^{-\beta}\bigr)\) forces
\[
\alpha(m-1)+2\beta=1,\qquad \alpha(p-1)-\sigma\beta=1 .
\]
Subtracting, we get \(\alpha(p-m)=(\sigma+2)\beta\), and therefore
\begin{equation}\label{eq:ssexp}
\frac{\alpha}{\beta}=\ac,\qquad
\beta=\frac1\theta,\qquad \alpha=\frac{\ac}{\theta},
\qquad\theta=\ac(m-1)+2>0,
\end{equation}
with \(\theta\) exactly the temporal scaling exponent of \eqref{eq:scaling}. In
particular, if the profile \(f\) has an algebraic tail,
\(f(\xi)\sim c_\infty\xi^{-\alpha/\beta}\) as \(\xi\to\infty\), then
\(\alpha=\beta\ac\) gives
\[
u(t,x)\sim c_\infty\abs{x}^{-\ac}\qquad\text{for }\abs{x}\gg t^{\beta}.
\]
Hence \emph{every} global self-similar solution with algebraic tail decays in
space at exactly the critical rate \(\ac\), with initial trace
\(c_\infty\abs{x}^{-\ac}\). Such solutions have recently been classified in
\cite{ISS26} for \(m\ge1\), \(p>\pF\) and \(\sigma>\max\{-2,-N\}\). For
\(m>1\) and \(\pF<p<p_S(\sigma)=\pss\), there is a compactly supported profile,
and there is \(b_\ast>0\) such that, for every \(b\in(0,b_\ast)\), the profile
\(f(\cdot;b)\) with \(f(0)=b\) satisfies
\(\xi^{\ac}f(\xi;b)\to c_\infty(b)\in(0,\infty)\) as \(\xi\to\infty\). For
\(p\ge\pss\), the latter property holds for every \(b>0\)
\cite[Theorem~1.2]{ISS26}. The exponents in \eqref{eq:ssexp} coincide with those
of \cite{ISS26}, since \(\theta=\bigl(\sigma(m-1)+2(p-1)\bigr)/(p-m)\). The
initial trace \(c_\infty\abs{x}^{-\ac}\) is unbounded at the origin, and it
lies outside the class \eqref{eq:class}. This obstruction disappears after a
time shift, since \(U(\tau,\cdot)\) is bounded for every \(\tau>0\). One then
only needs a comparison principle with bounded time-dependent supersolutions,
which is the evolutionary analogue of Lemma~\ref{lem:comparison}, see
Lemma~\ref{lem:comparisonU} and Lemma~\ref{lem:ssregular} below.
\end{remark}

\begin{proposition}\label{prop:selfsimilar}
Let \(N\ge3\), \(m>1\), \(\sigma\in(-2,0)\) and \(p>\pF\), and let
\begin{equation}\label{eq:U}
U(t,x)=t^{-\alpha}f\bigl(\abs{x}t^{-\beta}\bigr),\qquad
\alpha=\frac{\ac}{\theta},\quad\beta=\frac1\theta,\quad\theta=\ac(m-1)+2,
\end{equation}
where the profile \(f:[0,\infty)\to\R\) satisfies the following two
conditions.
\begin{enumerate}[label=\textup{(S\arabic*)},leftmargin=2.6em]
\item \label{eq:P1} \(f\in C([0,\infty))\cap C^2((0,\infty))\), \(f>0\) on
\([0,\infty)\), and \(f\) solves
\begin{equation}\label{eq:profile}
(f^m)''(\xi)+\frac{N-1}{\xi}\,(f^m)'(\xi)+\alpha f(\xi)+\beta\xi f'(\xi)
+\xi^{\sigma}f(\xi)^p=0,\qquad\xi>0 .
\end{equation}
\item \label{eq:P2} \(\xi^{\ac}f(\xi)\to c_\infty\in(0,\infty)\) as
\(\xi\to\infty\).
\end{enumerate}
Then \(L_\ast(\ac)\ge c_\infty\). In particular, \(L_\ast(\ac)>0\).
\end{proposition}

The proof relies on the following lemma, which shows that the time-shifted
self-similar function is an admissible supersolution in
Lemma~\ref{lem:comparisonU}. In particular, no Dirac mass appears at the
origin, and this is a consequence of the boundedness of \(f\) near \(\xi=0\).

\begin{lemma}\label{lem:ssregular}
Let \(N\ge3\), \(m>1\), \(\sigma\in(-2,0)\) and \(p>\pF\), and let \(f\) satisfy
\textup{(S1)} and \textup{(S2)}. Then the following hold.
\begin{enumerate}[label=\textup{(\roman*)},leftmargin=2.2em]
\item \(\xi^{N-1}(f^m)'(\xi)\to0\) as \(\xi\downarrow0\), and there is
\(C>0\) such that \(\abs{(f^m)'(\xi)}+\abs{f'(\xi)}\le C\xi^{\sigma+1}\) for
\(\xi\in(0,1]\).
\item For every \(\tau>0\) and \(T>0\), the function
\(U_\tau(t,x)\eqdef U(t+\tau,x)\) is continuous and positive on
\([0,T]\times\RN\), with \(U_\tau\le\tau^{-\alpha}\sup f\). It satisfies
\textup{(U1)}, and \textup{(U2)} with equality, on \([0,T)\times\RN\). More
precisely, \(\partial_tU_\tau\in L^\infty_{\mathrm{loc}}([0,T]\times\RN)\) and
\(\abs{\nabla U_\tau^m(t,x)}\le C_{R}\abs{x}^{\sigma+1}\) for
\((t,x)\in[0,T]\times B(0,R)\).
\end{enumerate}
\end{lemma}

\begin{proof}
(i) Set \(Q(\xi)\eqdef\xi^{N-1}(f^m)'(\xi)+\beta\xi^Nf(\xi)\) for \(\xi>0\).
Multiplying \eqref{eq:profile} by \(\xi^{N-1}\), we obtain
\(\bigl(\xi^{N-1}(f^m)'\bigr)'=-\xi^{N-1}\bigl(\alpha f+\beta\xi f'+\xi^{\sigma}f^p\bigr)\).
Since \((\beta\xi^Nf)'=N\beta\xi^{N-1}f+\beta\xi^Nf'\), we infer that
\[
Q'(\xi)=-(\alpha-N\beta)\,\xi^{N-1}f(\xi)-\xi^{N-1+\sigma}f(\xi)^p,\qquad\xi>0 .
\]
The right-hand side is integrable on \((0,1)\), because \(f\) is bounded on
\([0,1]\) and \(N-1+\sigma>-1\). Hence \(Q(\xi)\) has a finite limit \(Q_0\) as
\(\xi\downarrow0\). Assume that \(Q_0\ne0\). Since \(\beta\xi^Nf(\xi)\to0\),
there is \(\xi_1\in(0,1)\) such that \((f^m)'\) has the sign of \(Q_0\) and
\(\abs{(f^m)'(\xi)}\ge\abs{Q_0}\xi^{1-N}/2\) on \((0,\xi_1)\). Integrating, we get
\(\abs{f^m(\xi)-f^m(\xi_1)}\ge\frac{\abs{Q_0}}2\int_\xi^{\xi_1}r^{1-N}dr\to\infty\)
as \(\xi\downarrow0\), because \(N\ge3\). This contradicts the continuity of
\(f\) at \(0\). Therefore \(Q_0=0\) and
\[
\xi^{N-1}(f^m)'(\xi)=-\beta\xi^Nf(\xi)-\int_0^\xi\Bigl((\alpha-N\beta)r^{N-1}f(r)
+r^{N-1+\sigma}f(r)^p\Bigr)dr .
\]
In particular \(\xi^{N-1}(f^m)'(\xi)\to0\). Let \(M_1\eqdef\max_{[0,1]}f\) and
\(m_1\eqdef\min_{[0,1]}f>0\). For \(\xi\in(0,1]\), we have
\(\xi\le\xi^{\sigma+1}\), because \(\sigma+1<1\), and therefore
\[
\abs{(f^m)'(\xi)}\le\beta M_1\xi+\frac{\abs{\alpha-N\beta}M_1}{N}\,\xi
+\frac{M_1^p}{N+\sigma}\,\xi^{\sigma+1}\le C\xi^{\sigma+1}.
\]
Since \(f'=(f^m)'/(mf^{m-1})\) and \(f\ge m_1\) on \([0,1]\), the same bound holds
for \(f'\), with another constant.

(ii) By (S1) and (S2), \(f\) is continuous and positive on \([0,\infty)\) and
tends to \(0\) at infinity, so that \(0<f\le\sup f<\infty\). Hence \(U_\tau\) is
continuous and positive on \([0,T]\times\RN\), and
\(U_\tau\le\tau^{-\alpha}\sup f\). In particular,
\(U_\tau\in L^\infty((0,T)\times\RN)\cap C([0,T],L^1_{\mathrm{loc}}(\RN))\).
Write \(s=t+\tau\in[\tau,\tau+T]\) and \(\xi=\abs{x}s^{-\beta}\). Since
\(f>0\) and \(f\in C^2((0,\infty))\), we have \(f^m\in C^2((0,\infty))\), and
\(U_\tau\) is of class \(C^1\) in \(t\) and \(C^2\) in \(x\) on
\([0,T]\times(\RN\setminus\{0\})\), with
\[
\partial_tU_\tau=-s^{-\alpha-1}\bigl(\alpha f+\beta\xi f'\bigr)(\xi),\qquad
\nabla U_\tau^m=s^{-\alpha m-\beta}(f^m)'(\xi)\,\frac{x}{\abs{x}},
\]
\[
\Delta U_\tau^m=s^{-\alpha m-2\beta}\Bigl((f^m)''+\frac{N-1}{\xi}(f^m)'\Bigr)(\xi),
\qquad
\abs{x}^{\sigma}U_\tau^p=s^{\sigma\beta-\alpha p}\,\xi^{\sigma}f(\xi)^p .
\]
The exponents in \eqref{eq:U} satisfy \(\alpha(m-1)+2\beta=1\) and
\(\alpha(p-1)-\sigma\beta=1\). Indeed, the first identity is
\((\ac(m-1)+2)/\theta=1\), and the second one follows from
\(\ac(p-1)-\sigma=\ac(p-m)+\ac(m-1)-\sigma=\theta\), since \(\ac(p-m)=\sigma+2\).
Thus \(\alpha m+2\beta=\alpha p-\sigma\beta=\alpha+1\), and \eqref{eq:profile}
yields
\begin{equation}\label{eq:ssclassical}
\partial_tU_\tau=\Delta U_\tau^m+\abs{x}^{\sigma}U_\tau^p
\qquad\text{in }[0,T]\times(\RN\setminus\{0\}).
\end{equation}
Fix \(R>0\). For \((t,x)\in[0,T]\times B(0,R)\), we have
\(\xi\in[0,\Xi]\) with \(\Xi\eqdef R\tau^{-\beta}\). By (i) and the continuity
of \(f'\) and \((f^m)'\) on \([1,\Xi]\), the functions \(\xi\mapsto\xi f'(\xi)\)
and \(\xi\mapsto\xi^{-\sigma-1}(f^m)'(\xi)\) are bounded on \((0,\Xi]\). Hence
\(\partial_tU_\tau\) is bounded on \([0,T]\times B(0,R)\), and
\(\abs{\nabla U_\tau^m(t,x)}\le C_R\abs{x}^{\sigma+1}\) there. Since
\(2(\sigma+1)>-2>-N\), the function \(\abs{x}^{\sigma+1}\) is square integrable
on \(B(0,R)\).

For each \(t\), the function \(U_\tau^m(t,\cdot)\) is continuous on \(\RN\),
of class \(C^1\) on \(\RN\setminus\{0\}\), and its pointwise gradient is
locally integrable. Integrating by parts on \(\{\abs{x}>\eps\}\) against
\(\phi\in C^1_c(\RN)\), the boundary term on \(\{\abs{x}=\eps\}\) is bounded by
\(N\omega_N\eps^{N-1}\norm{\phi}_\infty\sup U_\tau^m\to0\). Hence the pointwise
gradient is the distributional gradient, and
\(U_\tau^m\in L^2(0,T;H^1_{\mathrm{loc}}(\RN))\). Since \(t\mapsto
U_\tau(t,x)\) is of class \(C^1\) on \([0,T]\) for every \(x\ne0\), the
distributional time derivative of \(U_\tau\) is the pointwise one, and
\(\partial_tU_\tau\in L^\infty_{\mathrm{loc}}([0,T]\times\RN)\). This proves
(U1).

Let \(\varphi\in C^\infty_c((0,T)\times\RN)\) and \(t\in(0,T)\). Green's formula
on \(\{\abs{x}>\eps\}\) gives
\[
\int_{\{\abs{x}>\eps\}}\nabla U_\tau^m\cdot\nabla\varphi\dx
=-\int_{\{\abs{x}>\eps\}}\Delta U_\tau^m\,\varphi\dx
+\int_{\{\abs{x}=\eps\}}\varphi\,\nabla U_\tau^m\cdot\nu\,d\Hmeas ,
\]
with \(\nu=-x/\abs{x}\). The boundary term is bounded by
\(N\omega_N\norm{\varphi}_\infty C_R\,\eps^{N+\sigma}\), which tends to \(0\)
since \(N+\sigma>0\). By \eqref{eq:ssclassical},
\(\Delta U_\tau^m=\partial_tU_\tau-\abs{x}^{\sigma}U_\tau^p\) on
\(\RN\setminus\{0\}\), and this function is locally integrable on \(\RN\),
because \(\partial_tU_\tau\) and \(U_\tau\) are locally bounded and
\(\sigma>-N\). Letting \(\eps\downarrow0\) by dominated convergence, and
integrating in \(t\), we obtain
\[
\int_0^T\!\!\int_{\RN}\Bigl(\partial_tU_\tau\,\varphi+\nabla U_\tau^m\cdot\nabla\varphi
-\abs{x}^{\sigma}U_\tau^p\varphi\Bigr)\dx\,dt=0 ,
\]
which is (U2) with equality.
\end{proof}

\begin{proof}[Proof of Proposition~\ref{prop:selfsimilar}]
Set \(M_f\eqdef\sup_{\xi\ge0}f(\xi)<\infty\), see the proof of
Lemma~\ref{lem:ssregular}(ii). Fix \(L\in(0,c_\infty)\). By \textup{(S2)},
there is \(\xi_0\ge1\) such that \(f(\xi)\ge L\xi^{-\ac}\) for \(\xi\ge\xi_0\),
and \(m_0\eqdef\min_{[0,\xi_0]}f>0\) by \textup{(S1)}. Choose \(\tau\in(0,1]\)
so small that \(\xi_0\tau^{\beta}\le1\) and \(\tau^{-\alpha}m_0\ge L\). For
\(\abs{x}\le\xi_0\tau^{\beta}\), we have
\(U(\tau,x)\ge\tau^{-\alpha}m_0\ge L\ge L\phi_{\ac}(x)\). For
\(\abs{x}\ge\xi_0\tau^{\beta}\), using \(\alpha=\beta\ac\), we get
\[
U(\tau,x)=\tau^{-\alpha}f\bigl(\abs{x}\tau^{-\beta}\bigr)
\ge\tau^{-\alpha}L\bigl(\abs{x}\tau^{-\beta}\bigr)^{-\ac}
=L\abs{x}^{-\ac}\ge L\phi_{\ac}(x).
\]
Hence
\begin{equation}\label{eq:comparaison-init}
L\phi_{\ac}(x)\le U(\tau,x)\qquad\text{for all }x\in\RN .
\end{equation}
Let \(u\) be the solution with datum \(L\phi_{\ac}\in\mathcal A_r\), \(r>r_0\).
For every \(T>0\), Lemma~\ref{lem:ssregular}(ii) shows that \(U_\tau\) satisfies
\textup{(U1)} and \textup{(U2)} on \([0,T)\times\RN\) and is positive, while
\textup{(U3)} is \eqref{eq:comparaison-init}. Lemma~\ref{lem:comparisonU} gives
\(u(t)\le U(t+\tau)\le(t+\tau)^{-\alpha}M_f\le\tau^{-\alpha}M_f\) for a.e.\
\(x\) and all \(t<\min\{T,\Tmax(L\phi_{\ac})\}\). Since \(T\) is arbitrary,
Assumption~\ref{ass:IL}(A) yields \(\Tmax(L\phi_{\ac})=\infty\). Thus
\((0,c_\infty)\subset\mathcal G_{\ac}\), and \(L_\ast(\ac)\ge c_\infty\).
\end{proof}

We now remove the positivity assumption in \textup{(S1)} for the profiles of
\cite{ISS26}. The key observation is that a nonnegative profile which is
positive near the origin and near infinity cannot vanish in between. This is
the self-similar counterpart of the persistence of positivity for porous
medium type equations: a zero of \(f\) at \(\xi_2>0\) with \(f>0\) beyond
\(\xi_2\) would correspond to a contact sphere \(\{\abs{x}=\xi_2t^{\beta}\}\)
moving outward into a region where the solution is positive. We formulate the
equation in the flux form which appears in the proof of
Lemma~\ref{lem:ssregular}, because this is the form satisfied by weak
continuations of a profile through a zero.

\begin{lemma}\label{lem:positivity}
Let \(N\ge3\), \(m>1\), \(\sigma\in(-2,0)\) and \(p>\pF\), and let \(\alpha\),
\(\beta\) be as in \eqref{eq:U}. Let \(f:[0,\infty)\to[0,\infty)\) be
continuous, with \(f(0)>0\) and \(f^m\in C^1((0,\infty))\), and set
\(Q(\xi)\eqdef\xi^{N-1}(f^m)'(\xi)+\beta\xi^Nf(\xi)\). Assume that
\begin{equation}\label{eq:Qflux}
Q(\xi)-Q(\zeta)=\int_\zeta^\xi\Bigl((N\beta-\alpha)\,r^{N-1}f(r)
-r^{N-1+\sigma}f(r)^p\Bigr)dr\qquad\text{for all }0<\zeta<\xi,
\end{equation}
and that \(f(\xi)>0\) for all sufficiently large \(\xi\). Then \(f>0\) on
\([0,\infty)\), \(f\in C^2((0,\infty))\), and \(f\) solves \eqref{eq:profile}
classically on \((0,\infty)\). Conversely, every positive solution
\(f\in C^2((0,\infty))\) of \eqref{eq:profile} satisfies \eqref{eq:Qflux}.
\end{lemma}

\begin{proof}
The converse statement is the computation of \(Q'\) carried out at the
beginning of the proof of Lemma~\ref{lem:ssregular}. Let us prove the direct
statement. The zero set \(Z\eqdef\{\xi>0:\ f(\xi)=0\}\) is closed in
\((0,\infty)\). Since \(f\) is continuous with \(f(0)>0\), and since \(f\) is
positive for large \(\xi\), there are \(0<\delta<R\) such that
\(Z\subset[\delta,R]\). Hence \(Z\) is compact. Assume by contradiction that
\(Z\ne\emptyset\), and let \(\xi_2\eqdef\max Z>0\). Then \(f(\xi_2)=0\) and
\(f>0\) on \((\xi_2,\infty)\). The function \(f^m\) is of class \(C^1\) and
nonnegative, and it attains its minimum \(0\) at the interior point \(\xi_2\).
Therefore \((f^m)'(\xi_2)=0\), and we conclude that \(Q(\xi_2)=0\).

Set \(M_2\eqdef\max_{[\xi_2,\xi_2+1]}f\) and
\(F(\xi)\eqdef\int_{\xi_2}^{\xi}f(r)\,dr\) for \(\xi\ge\xi_2\). Note that
\(F(\xi)>0\) for \(\xi>\xi_2\) and that \(F\) is nondecreasing. Since
\(N-1+\sigma>0\), \eqref{eq:Qflux} with \(\zeta=\xi_2\) gives, for
\(\xi\in[\xi_2,\xi_2+1]\),
\[
\abs{Q(\xi)}\le C_1F(\xi),\qquad
C_1\eqdef\abs{N\beta-\alpha}(\xi_2+1)^{N-1}+(\xi_2+1)^{N-1+\sigma}M_2^{p-1}.
\]
Using the definition of \(Q\), \(r\ge\xi_2\) and the monotonicity of \(F\), we
infer that
\[
(f^m)'(r)=r^{1-N}Q(r)-\beta rf(r)\le C_1\xi_2^{1-N}F(\xi)-\beta\xi_2f(r),
\qquad \xi_2<r\le\xi\le\xi_2+1 .
\]
We integrate this inequality over \((\xi_2,\xi)\). Since \(f^m(\xi_2)=0\), we
obtain
\[
f^m(\xi)\le F(\xi)\Bigl(C_1\xi_2^{1-N}(\xi-\xi_2)-\beta\xi_2\Bigr),
\qquad \xi\in(\xi_2,\xi_2+1].
\]
For \(\xi-\xi_2<\min\{1,\beta\xi_2^{N}/C_1\}\), the bracket is negative, while
\(F(\xi)>0\). Hence \(f^m(\xi)<0\), which is impossible. We conclude that
\(Z=\emptyset\), that is, \(f>0\) on \([0,\infty)\).

Since \(f>0\) and \(f^m\in C^1((0,\infty))\), we have \(f=(f^m)^{1/m}\in
C^1((0,\infty))\). By \eqref{eq:Qflux}, \(Q\) is of class \(C^1\), so that
\((f^m)'=\xi^{1-N}Q-\beta\xi f\) is of class \(C^1\). Hence \(f^m\in
C^2((0,\infty))\) and \(f\in C^2((0,\infty))\). Differentiating
\eqref{eq:Qflux} and using
\(Q'=\xi^{N-1}\bigl((f^m)''+\frac{N-1}{\xi}(f^m)'\bigr)+N\beta\xi^{N-1}f+\beta\xi^Nf'\),
we recover \eqref{eq:profile}.
\end{proof}

\begin{remark}\rm \label{rem:S1S2}
Let us state precisely how the profiles of \cite{ISS26} enter. The profile
equation \eqref{eq:profile} is \cite[Eq.~(1.7)]{ISS26}, with the same
exponents \(\alpha\) and \(\beta\), since
\(\theta=\bigl(\sigma(m-1)+2(p-1)\bigr)/(p-m)\). The standing assumptions
\cite[Eq.~(1.4)]{ISS26}, namely \(m\ge1\), \(p>\pF\) and
\(\sigma>\max\{-2,-N\}\), cover our range. For \(b\in\mathcal B_{\rm alg}\),
the profile \(f(\cdot;b)\) is the solution of \eqref{eq:profile} determined by
\(f(0;b)=b\) and the expansion \cite[Eq.~(1.10)]{ISS26} at \(\xi=0\), and it
satisfies \eqref{eq:ISStail}. In particular, \(f(\cdot;b)\) is continuous at
\(0\), positive on a neighbourhood of \(0\) and, by \eqref{eq:ISStail},
positive for all large \(\xi\). On its maximal interval of positivity
\([0,\xi_{\max})\) containing \(0\), the profile is a classical solution of
\eqref{eq:profile}. Two situations may a priori be considered. If
\(\xi_{\max}=\infty\), then \(f(\cdot;b)>0\) on \([0,\infty)\) directly. If
\(\xi_{\max}<\infty\), the limit \eqref{eq:ISStail} only makes sense for a
continuation of the profile beyond its first zero. Such a continuation is
meaningful as a nonnegative weak solution, which amounts to the flux
formulation \eqref{eq:Qflux} together with \(f^m\in C^1\), as for the contact
condition \cite[Eq.~(1.13)]{ISS26} and the discussion following
\cite[Theorem~1.2]{ISS26}. Lemma~\ref{lem:positivity} then shows that this
second situation cannot occur. In both cases, \(f(\cdot;b)\) satisfies
\textup{(S1)} and \textup{(S2)}, with \(c_\infty=c_\infty(b)\).
\end{remark}

\begin{proof}[Proof of Theorem~\ref{thm:allp}]
Let \(b\in\mathcal B_{\rm alg}\). By Remark~\ref{rem:S1S2}, the profile
\(f(\cdot;b)\) satisfies \textup{(S1)} and \textup{(S2)} with
\(c_\infty=c_\infty(b)\). Let \(u_0\in\mathcal A_r\), \(r>r_c\), satisfy
\(u_0\le U(\tau,\cdot)\) for some \(\tau>0\). By Lemma~\ref{lem:ssregular}(ii),
for every \(T>0\), the function \(U_\tau=U(\cdot+\tau,\cdot)\) is positive and
satisfies \textup{(U1)} and \textup{(U2)} on \([0,T)\times\RN\), and \textup{(U3)}
holds by assumption. Lemma~\ref{lem:comparisonU} yields
\(u(t)\le U(t+\tau)\le\tau^{-\alpha}M_b\) for
\(t<\min\{T,\Tmax(u_0)\}\). Since \(T\) is arbitrary, the blow-up alternative
of Assumption~\ref{ass:IL}(A) gives \(\Tmax(u_0)=\infty\). Next,
Proposition~\ref{prop:selfsimilar} gives \(L_\ast(\ac)\ge c_\infty(b)\), and we
take the supremum over \(b\in\mathcal B_{\rm alg}\). The upper bound
\(L_\ast(\ac)\le2^{\ac}M_0=L_0\) is Theorem~\ref{thm:threshold}(ii). Finally,
\(L_\ast(\gamma)=0\) for \(\gamma<\ac\) by Theorem~\ref{thm:threshold}(iii),
while \(L_\ast(\gamma)\ge L_\ast(\ac)>0\) for \(\gamma\ge\ac\) by the
monotonicity in Theorem~\ref{thm:threshold}(ii). This proves that the zero set
of \(L_\ast\) is \((0,\ac)\).
\end{proof}

%=====================================================================
\section{Discussion}\label{sec:discussion}
%=====================================================================
\subsection{The emerging picture}
For \eqref{eq:CP} in the regime \(m>1\), \(\sigma\in(-2,0)\), \(N\ge3\), the
critical exponents are ordered as \(1<\pG<m<\pF<\pc<\pss\), and the dynamics
split as follows. For \(1<p\le\pF\), the theory of \cite{IL1} is complete. All
solutions with integrable bounded data are global (with grow-up) for
\(p\in(1,\pG]\), and every nontrivial solution blows up in finite time for
\(\pG<p\le\pF\). For \(p>\pF\), both global and blowing-up solutions coexist,
and the separating mechanism is the spatial decay of \(u_0\). The present work
identifies the second critical exponent as
\[
\ac=\frac{\sigma+2}{p-m}=\frac{N}{r_0+1}
=\bigl(\text{decay rate of the singular steady state }S\bigr).
\]
It proves universal blow-up for \(\gamma<\ac\) and global existence of small
data for \(\gamma>\ac\). It also describes the borderline \(\gamma=\ac\) for
\(p>\pc\), where \(S\) exists. There, a finite amplitude threshold
\(L_\ast(\ac)\in[A,L_0]\) separates global existence from blow-up
(Corollary~\ref{cor:sandwich}), and the amplitude \(A\) itself is attained by
global solutions whenever the Joseph--Lundgren condition fails
(Proposition~\ref{rem:LequalA}). Moreover, all bounded data lying below
\(\vartheta S\) with \(\vartheta<1\) produce global bounded solutions
(Corollary~\ref{cor:belowS}). The whole picture is encoded in the nondecreasing
threshold function \(L_\ast\) of Theorem~\ref{thm:threshold}. For every
\(p>\pF\), its zero set contains \((0,\ac)\) and is contained in \((0,\ac]\),
and it is exactly \((0,\ac)\) (Corollary~\ref{cor:char} and
Theorem~\ref{thm:allp}). The
singular steady state is the organizing object. It produces \(\ac\), it traps
the flow from above, and, after being matched to an inner profile, it produces
bounded barriers at criticality.

\subsection{\texorpdfstring{{Explicit constants}}{Explicit constants}, and what remains}

The two constants \(A\), below which global existence is proved, and \(L_0\),
above which blow-up is proved, are explicit, see \eqref{eq:explicit} and
\eqref{eq:L0}. The lower constant \(A\) is optimal for barriers of the form
used in Theorem~\ref{thm:conj91}, since a multiple \(\vartheta S\) is a
supersolution only for \(\vartheta\le1\). The upper constant \(L_0\) is not
expected to be sharp. It may be optimized over the cut-off profile \(\zeta\)
in \eqref{eq:zeta}, which leads to the constant \(L_0^{\ast}\) of
\eqref{eq:L0star}, a variational problem of Hardy--Sobolev type on an annulus.
Conjecture~\ref{conj:threshold} asserts that \(L_\ast(\ac)=A\) under the
Joseph--Lundgren condition \eqref{eq:JL}, while Remark~\ref{rem:noglue}
explains why \(L_\ast(\ac)>A\) is to be expected when \eqref{eq:JL} fails. The
explicit barrier of Theorem~\ref{thm:barrier} only gives the weaker bound
\(\mu A\). By Remark~\ref{rem:scaling}, this loss cannot be removed by
rescaling, and it is due to the explicit inner profile, as
Theorem~\ref{thm:conj91} shows.

\subsection{On the threshold \texorpdfstring{$\pc$}{pc} and the low-dimensional cases}
The singular steady state \(S=A\abs{x}^{-\ac}\) exists only for \(p>\pc\) and
\(N\ge3\). In the range \(\pF<p\le\pc\), forward self-similar solutions
replace it, see Theorem~\ref{thm:allp}, Lemma~\ref{lem:positivity} and Open
Problem~\ref{prob:pcrange}. In dimensions \(N\in\{1,2\}\), the admissible range is
\begin{equation}\label{eq:standing}
m>1,\qquad p>1,\qquad \sigma\in\bigl(\max\{-2,-N\},0\bigr),
\end{equation}
which gives \(\sigma\in(-1,0)\) for \(N=1\). No algebraic singular steady state
of the above form exists in these dimensions. Indeed, by \eqref{eq:laplacian}
one would need \(\ac m<N-2\), which is impossible for \(N\le2\) since
\(\ac m>0\). The blow-up side of the theory, however, does not use \(N\ge3\).
The annulus computation of Lemma~\ref{lem:energy-annulus} is dimension free,
and Assumption~\ref{ass:IL}(B) only relies on the comparison principle for the
regularized problems \eqref{eq:Reta}. Moreover, the results of \cite{IL1} which
we use, namely \cite[Theorems~2.2 and 2.9, Corollary~2.7, Propositions~3.3 and
3.4]{IL1}, are stated under \eqref{eq:standing} in every dimension. The
uniqueness result \cite[Theorem~2.5]{IL1}, on which \cite[Lemma~3.6 and
Corollary~2.7]{IL1} rely, applies for \(N=2\) when \(p>\pG\), and for \(N=1\)
when \(p>1-\sigma(m-1)\). Both conditions hold for \(p>\pF\), because
\(\pF>m>\max\{\pG,1-\sigma(m-1)\}\) (recall that \(\sigma>-1\) when \(N=1\)).
Hence the proofs of Theorems~\ref{thm:slow} and \ref{thm:blowup-crit}, of
Proposition~\ref{prop:nodecayonly} and of Theorem~\ref{thm:threshold}(i)--(iii)
apply verbatim for \(N\in\{1,2\}\). What is genuinely missing in low dimensions
is the barrier, hence the global existence side at and near criticality.
Adapting the weighted estimates used in the proof of \cite[Theorem~2.5]{IL1} to
construct barriers at criticality in low dimensions is a natural next step.

\subsection{Relation to \texorpdfstring{\cite{IL2}}{[IL2]} and to the classical theory}

The non-existence threshold of \cite{IL2} concerns \(\sigma>0\), where the
potential \emph{grows} at infinity, and where sufficiently slow decay of the
datum produces non-existence of any solution, at the rate \(\sigma/(p-1)\) for
\(p\ge\pG\). The present problem with \(\sigma<0\) is complementary, and in a
precise sense softer. The rate \(\sigma/(p-1)\) becomes negative, so that the
potential decays and slow tails do not preclude local solvability. As we
prove, however, slow tails still force finite-time blow-up, and the borderline
is the amplitude-sensitive critical decay of Theorem~\ref{thm:threshold} (see
Remark~\ref{rem:henon}). Both regimes are nevertheless connected by the same
singular steady state \(S=A\abs{x}^{-\ac}\) and by the same exponent \(\pc\),
which are exhibited in \cite[\S3(c)]{IL2} and used systematically here. This
suggests a unified treatment of the two sign regimes. We also note the
structural asymmetry pointed out in Remark~\ref{rem:why-eta}. The sign
\(\sigma<0\) is what makes the singular stationary profile a supersolution of
every regularized problem, and hence what makes the comparison theory developed
in this paper work.

The whole circle of questions is the quasilinear, singular-weight extension of
the classical Lee--Ni theory \cite{LeeNi92,GuiWang95} for
\(\partial_t u=\Delta u+u^{p}\), whose critical decay \(2/(p-1)\) is recovered
from \eqref{eq:ac} in the formal limit \(m=1\), \(\sigma=0\). It also extends
the quasilinear counterparts of this theory
\cite{AndreucciDiBenedetto91,Qi98,QiLevine93,MochizukiSuzuki97,MMH00,GuoGuo01,Suzuki02,LiMu04,MuLiWang10,MMZ11,
Galaktionov94,SGKM95}. We refer to \cite{AndreucciTedeev05} for universal
bounds and blow-up rates in the Fujita range, to \cite{HMM98} for systems, and
to \cite{ZMHZ16} for a doubly singular equation with a weighted source.
Weighted semilinear versions, in which the influence of a space-dependent
coefficient on global solvability is analysed, can be found in
\cite{Pinsky97,BTW17,CIT22}, and weighted $p$-Laplacian versions in
\cite{YangYangZheng12}. The role of the Joseph--Lundgren exponent
for data near the singular steady state is analysed in
\cite{GuiNiWang92,Wang93,PolacikYanagida03,GalaktionovVazquez97,SouWei03,Naito12}
for the semilinear equation. See also \cite[Ch.~IV]{QS19} for the semilinear state of the art,
\cite{Levine90,DengLevine00} for surveys, and \cite{Fujita66} for the origin of
the subject.

\section*{Declarations}
\noindent\textbf{Funding.} No funding was received for this work.

\noindent\textbf{Conflict of interest.} The authors declare that they have no
conflict of interest.

\noindent\textbf{Data availability.} No datasets were generated or analysed in
this work. 

%=====================================================================

\end{document}